\documentclass[reqno,a4paper]{amsart}

\usepackage[english]{babel}
\usepackage{amsmath,amsthm,amssymb,amsfonts}
\usepackage{mathrsfs}

\usepackage{enumitem}
\setlist[enumerate]{leftmargin=*}

\usepackage{hyperref}

\numberwithin{equation}{section}   %%numera le equazioni sezione per sezione

\newtheorem{theorem}{Theorem}[section]
\newtheorem{lemma}[theorem]{Lemma}
\newtheorem{proposition}[theorem]{Proposition}
\newtheorem{corollary}[theorem]{Corollary}

\theoremstyle{definition}

\newcommand{\spnt}[2]{\left\lfloor #1, #2 \right\rceil}

\DeclareMathOperator{\Span}{span}
\DeclareMathOperator{\supp}{supp}
\DeclareMathOperator*{\esssup}{ess\,sup}

\newcommand{\id}{\mathrm{id}}

\newcommand{\tc}{\,:\,}

\newcommand{\RR}{\mathbb{R}}
\newcommand{\CC}{\mathbb{C}}
\newcommand{\NN}{\mathbb{N}}
\newcommand{\ZZ}{\mathbb{Z}}

\newcommand{\DimD}{d}
\newcommand{\DimK}{k}                 % dimension of the equatorial sphere

\newcommand{\DDJJ}{j_*}             % special indices
\newcommand{\DDKK}{{k_*}}    % in the proof of the elliptic spectral cluster estimate

\newcommand{\compl}{\mathsf{c}}

\newcommand{\cI}{\mathcal{I}}          % splitting the integral
\newcommand{\cS}{\mathcal{S}}          % in the proof of estimate for weight

\newcommand{\defeq}{\mathrel{:=}}
\newcommand{\eqdef}{\mathrel{=:}}

\newcommand{\Harm}{\mathcal{H}}
\newcommand{\dimHarm}{\alpha}

\newcommand{\potV}{\mathfrak{V}}     % `potential' in Grushin operator

\newcommand{\SPerm}{\mathfrak{S}}   % permutation group

\newcommand{\sfera}{\mathbb{S}}   % sphere
\newcommand{\equat}{\mathbb{E}}   % equator

\newcommand{\SO}{\mathrm{SO}}

\newcommand{\Hz}{\mathsf{H}}    % horizontal distribution
\newcommand{\fZ}{\mathcal{Z}}    % family of Hoermander vector fields

\newcommand{\opL}{\mathcal{L}}   % spherical Grushin operator
\newcommand{\opG}{\mathcal{G}}   % plane Grushin operator

\newcommand{\AbsOp}{\mathfrak{L}}   % operator in abstract theorem
\newcommand{\absDim}{\mathfrak{d}}  % dimension in abstract theorem
\newcommand{\absWeight}{\pi}

\newcommand{\dist}{\varrho}       % control distance
\newcommand{\weight}{\varpi}      % weight
\newcommand{\meas}{\sigma}        % measure on sphere
\newcommand{\Kern}{\mathcal{K}}  % integral kernel

\newcommand{\Sob}[2]{L^{#1}_{#2}}  % Sobolev norm

\newcommand{\tX}{\widetilde{X}}      % modified spherical harmonic S^d 
\newcommand{\tSpHddd}{\mathscr{X}}

\DeclareFontFamily{U}{matha}{\hyphenchar\font45}
\DeclareFontShape{U}{matha}{m}{n}{
      <5> <6> <7> <8> <9> <10> gen * matha
      <10.95> matha10 <12> <14.4> <17.28> <20.74> <24.88> matha12
      }{}
\DeclareSymbolFont{matha}{U}{matha}{m}{n}
\DeclareFontSubstitution{U}{matha}{m}{n}

\DeclareFontFamily{U}{mathx}{\hyphenchar\font45}
\DeclareFontShape{U}{mathx}{m}{n}{
      <5> <6> <7> <8> <9> <10>
      <10.95> <12> <14.4> <17.28> <20.74> <24.88>
      mathx10
      }{}
\DeclareSymbolFont{mathx}{U}{mathx}{m}{n}
\DeclareFontSubstitution{U}{mathx}{m}{n}

\DeclareMathDelimiter{\vvvert}{0}{matha}{"7E}{mathx}{"17}

\makeatletter
\def\author@andify{%
  \nxandlist {\unskip ,\penalty-1 \space\ignorespaces}%
    {\unskip {} \@@and~}%
    {\unskip \penalty-2 \space \@@and~}%
}
\makeatother

\begin{document}

\title[Ultraspherical Grushin operators]{
Weighted spectral cluster bounds\\ and a sharp multiplier theorem\\ for ultraspherical Grushin operators}
\author{Valentina Casarino}
\address{Universit\`a degli Studi di Padova\\Stradella san Nicola 3 \\I-36100 Vicenza \\ Italy}
\email{valentina.casarino@unipd.it}
\author{Paolo Ciatti}
\address{Universit\`a degli Studi di Padova\\Via Marzolo 9 \\I-35100 Padova \\ Italy}
\email{paolo.ciatti@unipd.it}
\author{Alessio Martini}
\address{School of Mathematics \\ University of Birmingham \\ Edgbaston \\ Birmingham \\ B15 2TT \\ United Kingdom}
\email{a.martini@bham.ac.uk}

\keywords{Grushin spheres, spectral multipliers,  sub-Laplacian, sub-Riemannian geometry, hyperspherical harmonics.}
\subjclass[2010]{%
33C55, %Spherical harmonics,
42B15,  % multipliers
43A85 % harmonic analysis on homogeneous spaces
 (primary);
53C17, %geometria subriemanniana
58J50  %problemi spettrali
(secondary).
}

\thanks{
The first and the second author were partially supported by GNAMPA (Project 2018
``Operatori e disuguaglianze integrali in spazi con simmetrie")
and MIUR (PRIN 2016 ``Real and Complex Manifolds: Geometry, Topology and Harmonic Analysis").
Part of this research was carried out while the third author was visiting the Dicea, Universit\`a di Padova,  Italy,  as a recipient of a ``Visiting Scientist 2019'' grant;
he gratefully thanks the Universit\`a di Padova for the support and hospitality. 
The authors are members of the Gruppo Nazionale per l'Analisi Matematica, la Probabilit\`a e le loro Applicazioni (GNAMPA) of the Istituto Nazionale di Alta Matematica (INdAM)}

\begin{abstract}
{We study degenerate elliptic operators of Grushin type on the $d$-dimensional sphere, which are singular on a $k$-dimensional sphere for some $k < d$. 
For these operators we prove a spectral multiplier theorem of Mihlin--H\"ormander type, which is optimal whenever $2k \leq d$, and a corresponding 
Bochner--Riesz summability result. The proof hinges on suitable weighted spectral cluster bounds, which in turn depend on precise estimates for ultraspherical polynomials.
}

\end{abstract}

\maketitle

\section{Introduction}
In this  paper we 
continue the study of spherical Grushin-type operators
started in \cite{CaCiaMa} 
with the case of the two-dimensional sphere.
The focus here is on a family of hypoelliptic operators $\{\opL_{\DimD,\DimK} \}_{1\leq \DimK < \DimD}$, acting on functions defined 
on  the unit sphere   $\sfera^\DimD$ in $\RR^{1+\DimD} $, i.e., on
\begin{equation}\label{eq:sfera}
\sfera^\DimD = \{ (z_0,\dots,z_\DimD) \in \RR^{1+\DimD} \tc z_0^2 + \dots + z_\DimD^2 = 1 \},
\end{equation}
 for some $d\geq 2$.
As it is well known, the groups $\SO(1+r)$ with $1 \leq r \leq \DimD$ can be naturally identified with a sequence of nested subgroups of $\SO(1+\DimD)$ and correspondingly they act on $\sfera^\DimD$ by rotations. We denote by $\Delta_r$ the (positive semidefinite) second-order differential operator on $\sfera^\DimD$ 
corresponding through this action to the Casimir operator on $\SO(1+r)$. 
The operators
 $\Delta_r$ commute pairwise and $\Delta_\DimD$ turns out to be
the Laplace--Beltrami operator on $\sfera^\DimD$. 
The operators we are interested in are defined as
\begin{equation}\label{eq:def_Grushin_op}
\opL_{\DimD,\DimK} = \Delta_\DimD - \Delta_\DimK,
\end{equation}
with $\DimK=1,\dots,\DimD-1$.
By introducing a
 suitable system of ``cylindrical coordinates''
 $(\omega,\psi)$ on $\sfera^d$, where
$\omega \in \sfera^\DimK$ and
 $\psi = (\psi_{\DimK+1},\ldots,\psi_\DimD) \in (-\pi/2,\pi/2)^{\DimD-\DimK}$ (see Section \ref{ss:cylcoords} below for details),
one can write
$\opL_{\DimD,\DimK}$ 
more explicitly as
\begin{equation}\label{eq:Ldk-fields}
\opL_{\DimD,\DimK} = \sum_{r=\DimK+1}^\DimD Y^+_r Y_r + \potV(\psi) \, \Delta_\DimK,
\end{equation}
where the $Y_r$ and 
their formal adjoints $Y_r^+$ (with respect to the standard rotation-invariant measure $\meas$ on $\sfera^d$) are vector fields 
only depending on $\psi$, to wit,
\begin{equation}\label{eq:vfield_Yr}
Y_r = \frac{1}{\cos\psi_{r+1} \cdots \cos \psi_\DimD} \frac{\partial}{\partial \psi_r},
\end{equation}
and $\potV : (-\pi/2,\pi/2)^{\DimD-\DimK} \to \RR$ is given by
\begin{equation}\label{eq:potential}
\potV(\psi)
= \frac{1}{\cos^2 \psi_{\DimK+1} \cdots \cos^2 \psi_\DimD} - 1
= \prod_{j=\DimK+1}^\DimD (1 + \tan^2 \psi_j) - 1 .
\end{equation}

Since $\potV(\psi)$ vanishes only for $\psi = 0$, 
the formulae above show that
  each  $\opL_{\DimD,\DimK}$ is elliptic away from the $\DimK$-submanifold  $\sfera^\DimK \times \{0\}$ of $\sfera^\DimD$;
the
 loss of global ellipticity is anyway compensated by the fact that each $\opL_{\DimD,\DimK}$ is hypoelliptic and satisfies subelliptic estimates,
as shown by an application of H\"ormander's theorem for sums of squares of vector fields \cite{Hor}. 
  Indeed the expression \eqref{eq:Ldk-fields} reveals the analogy of the operators $\opL_{\DimD,\DimK}$
  with
certain degenerate elliptic
 operators $\mathcal G_{\DimD,\DimK}$ on $\RR^\DimD$,
given by
\begin{equation}\label{def:Grushin}
 \mathcal G_{\DimD,\DimK}=\Delta_{x}+|x|^2\Delta_{y},
\end{equation}
 where $x,y$ are the components of a point in $ \RR^{\DimD-\DimK}_x \times \RR^\DimK_y$
 and $\Delta_{x}$, $\Delta_{y}$ denote the corresponding (positive definite) partial Laplacians.
 
 In light of \cite{Grushin, Grushin2}, the operators $\mathcal G_{\DimD,\DimK}$ are often called
 Grushin operators;
sometimes they are also called Baouendi--Grushin operators, since shortly before the papers by V.\ V.\ Grushin appeared,
M.\ S.\ Baouendi introduced a more general class of operators containing also the
 $\mathcal G_{\DimD,\DimK}$ \cite{Ba}. In these and other works (see, e.g., \cite{Franchi,RS,DM}), the coefficient $|x|^2$ in \eqref{def:Grushin} 
may be replaced by a more general function $V(x)$.
As prototypical examples of differential operators with mixed homogeneity, operators of the form \eqref{def:Grushin}
   have attracted increasing interest
in the last fifty years;
we refer to \cite{CaCiaMa}
for a brief list of the main results, focused on the field of harmonic analysis. 
More recently,
the study of Grushin-type operators began to develop also on more general manifolds than $\RR^n$, from both a geometric  and an analytic perspective \cite{BFIuno, BFI, Pesenson, 
BPS, BoL, GMP, GMP2}.

In this article, we investigate $L^p$ boundedness properties of operators of the form $F(\sqrt{\opL_{\DimD,\DimK}})$ in connection with size and smoothness properties of the spectral multiplier $F : \RR \to \CC$; here $L^p$ spaces on the sphere $\sfera^d$ are defined in terms of the spherical measure $\meas$, and the operators $F(\sqrt{\opL_{\DimD,\DimK}})$ are initially defined on $L^2(\sfera^d)$ via the Borel functional calculus for the self-adjoint operator $\opL_{\DimD,\DimK}$. 
The study of the $L^p$ boundedness of functions of Laplace-like operators is a classical and very active area of harmonic analysis, with a number of celebrated results and open questions, already in the case of the classical Laplacian in Euclidean space (think, e.g., of the Bochner--Riesz conjecture). Regarding the spherical Grushin operators $\opL_{\DimD,\DimK}$, in the case $\DimD=2$ and 
$\DimK=1$ a sharp multiplier theorem of Mihlin--H\"ormander type and a Bochner--Riesz summability result for $\opL_{\DimD,\DimK}$ were obtained in \cite{CaCiaMa}.
Here we treat the general case $\DimD \geq 2$, $1 \leq \DimK < \DimD$, and obtain the following result.

Let $\eta \in C^\infty_c((0,\infty))$ be any nontrivial cutoff, and denote by $\Sob{q}{s}(\RR)$ the $L^q$ Sobolev space of (fractional) order $s$ on $\RR$.

\begin{theorem}\label{thm:main}
Let $D = \max\{\DimD,2\DimK\}$ and $s > D/2$.
\begin{enumerate}[label=(\roman*)]
\item\label{en:main_cptsup} For all continuous functions supported in $[-1,1]$,
\[
\sup_{t>0} \|F(t\sqrt{\opL_{\DimD,\DimK}})\|_{L^1(\sfera^\DimD) \to L^1(\sfera^\DimD)} \lesssim_s \|F\|_{\Sob{2}{s}}.
\]
\item\label{en:main_mh} For all bounded Borel functions $F : \RR \to \CC$ such that $F|_{(0,\infty)}$ is continuous,
\[
\|F(\sqrt{\opL_{\DimD,\DimK}})\|_{L^1(\sfera^\DimD) \to L^{1,\infty}(\sfera^\DimD)} \lesssim_s \sup_{t\geq 0}\|F(t\cdot) \, \eta\|_{\Sob{2}{s}}.
\]
Hence, whenever the right-hand side is finite, the operator $F(\sqrt{\opL_{\DimD,\DimK}})$ is of weak type $(1,1)$ and bounded on $L^p(\sfera^\DimD)$ for all $p \in (1,\infty)$.
\end{enumerate}
\end{theorem}

Part \ref{en:main_cptsup} of the above theorem and a standard interpolation technique imply the following Bochner--Riesz summability result.

\begin{corollary}\label{cor:main}
Let $D=\max\{\DimD,2\DimK\}$ and $p \in [1,\infty]$. If $\delta > (D-1) |1/2-1/p|$, then the Bochner--Riesz means $(1-t\opL_{\DimD,\DimK})_+^\delta$ of order $\delta$ associated with $\opL_{\DimD,\DimK}$ are bounded on $L^p(\sfera^\DimD)$ uniformly in $t \in (0,\infty)$.
\end{corollary}

It is important to point out that weaker versions of the above results, involving more restrictive requirements on the smoothness parameters $s$ and $\delta$, could be readily obtained by standard techniques. Indeed the sphere $\sfera^d$, with the measure $\sigma$ and the Carnot--Carath\'eodory distance associated to $\opL_{\DimD,\DimK}$, is a doubling metric measure space of ``homogeneous dimension'' $Q = d+k$, and the operator $\opL_{\DimD,\DimK}$ satisfies Gaussian-type heat kernel bounds. As a consequence (see, e.g., \cite{He,CoS,DOS,DzS}), one would obtain the analogue of Theorem \ref{thm:main} with smoothness requirement $s>Q/2$, measured in terms of an $L^\infty$ Sobolev norm, and the corresponding result for Bochner--Riesz means would give $L^p$ boundedness only for $\delta > Q |1/p - 1/2|$. Since $Q > D > D-1$, the results in this paper yield an improvement on the standard result for all values of $\DimD$ and $\DimK$.

As a matter of fact, in the case $\DimK \leq \DimD/2$, the above multiplier theorem is sharp, in the sense that the lower bound 
$D/2$
 to the order of smoothness $s$ required in Theorem \ref{thm:main} cannot be replaced by any smaller quantity. Since $\opL_{\DimD,\DimK}$ is elliptic away from a negligible subset of $\sfera^\DimD$, and 
 $D = \DimD$ 
 is the topological dimension of $\sfera^\DimD$ when $\DimK \leq \DimD/2$, the sharpness of the above result can be seen by comparison to the Euclidean case via a transplantation technique \cite{Mi,KST}.

The fact that for subelliptic nonelliptic operators one can often obtain ``improved'' multiplier theorems, by replacing the relevant homogeneous dimension with the topological dimension in the smoothness requirement, was first noticed in the case of sub-Laplacians on Heisenberg and related groups by D. M\"uller and E. M. Stein \cite{MuSt} and independently by W. Hebisch \cite{He_heis}, and has since been verified in multiple cases. However, despite a flurry of recent progress (see, e.g., \cite{MMuGAFA,CaCiaMa,DM,MMuNG} for more detailed accounts and further references), the question whether such an improvement is always possible remains open. The results in the present paper can therefore be considered as part of a wider programme, attempting to gain an understanding of the general problem by tackling particularly significant particular cases.

In these respects, it it relevant to point out that Theorem \ref{thm:main} above can be considered as a strengthening of the multiplier theorem for the Grushin operators $\opG_{\DimD,\DimK}$ on $\RR^\DimD$ proved in \cite{MSi}: indeed a ``nonisotropic transplantation'' technique (see, e.g., \cite[Theorem 5.2]{M}) allows one to deduce from Theorem \ref{thm:main} the analogous result where $\sfera^d$ and $\opL_{\DimD,\DimK}$ are replaced by  $\RR^d$ and $\opG_{\DimD,\DimK}$.

\bigskip

The  structure of the proof of Theorem \ref{thm:main} 
broadly follows  that of the analogous result in \cite{CaCiaMa}, but additional difficulties need to be overcome here.
An especially delicate point is the proof of the ``weighted spectral cluster estimates'' stated as Propositions \ref{prp:indici_bassi_primo_piufattori} and \ref{prp:indici_bassi_primo_peso_secondo_stratopiu} below, essentially consisting in suitable weighted $L^1 \to L^2$ norm bounds for ``weighted spectral projections''
\begin{equation}\label{eq:weighted_cluster}
(\opL_{\DimD,\DimK}/\Delta_\DimD)^{\alpha/2} \chi_{[i,i+1]}(\sqrt{\opL_{\DimD,\DimK}})
\end{equation}
associated with bands of unit width of the spectrum of $\sqrt{\opL_{\DimD,\DimK}}$.
These can be thought of as subelliptic analogues of the Agmon--Avakumovi\v{c}--H\"ormander spectral cluster estimates
\begin{equation}\label{eq:hormander_cluster}
\|\chi_{[i,i+1]}(\sqrt{\Delta_{\DimD}})\|_{L^1 \to L^2}  \lesssim i^{(\DimD-1)/2}
\end{equation}
for the elliptic Laplacian $\Delta_\DimD$,
which are valid more generally when $\sqrt{\Delta_\DimD}$ is replaced with an elliptic pseudodifferential operator of order one on a compact $\DimD$-manifold \cite{HorSpec}, and are the basic building block for a sharp multiplier theorem for elliptic operators on compact manifolds and related restriction-type estimates \cite{SoggeRiesz,SoggeClusters,SeeSo,Frank}. Thanks to pseudodifferential and Fourier integral operator techniques, estimates of the form \eqref{eq:hormander_cluster} can be proved for elliptic operators in great generality, but these techniques break down when the ellipticity assumption is weakened. Nevertheless alternative ad-hoc methods may be developed in many cases, based on a detailed analysis of the spectral decomposition of the operator under consideration, often made possible by underlying symmetries.

In the case of the spherical Grushin operator $\opL_{\DimD,\DimK}$, as a consequence of its spectral decomposition in terms of joint eigenfunctions of the operators $\Delta_\DimD,\dots,\Delta_\DimK$,
the integral kernel of the ``weighted projection'' in \eqref{eq:weighted_cluster} involves sums of $(\DimD-\DimK)$-fold tensor
 products of ultraspherical polynomials.
This is a substantial difference from the case considered in \cite{CaCiaMa} (where $\DimD-\DimK=1$) 
 and requires new ideas and greater care. Section \ref{s:weighted} of this paper is devoted to the proof of these estimates.
As in \cite{CaCiaMa}, here we make fundamental use of precise estimates for ultraspherical polynomials, which are uniform in suitable ranges of indices. These estimates, which are consequences of the asymptotic approximations of  \cite{Olver75, Olver-libro, Olver, BoydDunster}, could be of independent interest,
and their derivation is presented in an auxiliary paper \cite{CCM-Jacobi}.

In the context of subelliptic operators on compact manifolds, ``weighted spectral cluster estimates'' were first obtained in the seminal work of Cowling and Sikora \cite{CoS} for a distinguished sub-Laplacian on $\mathrm{SU}(2)$, leading to a sharp multiplier theorem in that case; their technique was then applied to many different frameworks \cite{CoKS, CCMS, M, ACMMu}.
However, the general theory developed in \cite{CoS}, based on spectral cluster estimates involving a single weight function, does not seem to be directly applicable to the spherical Grushin operator $\opL_{\DimD,\DimK}$ (which, differently from the sub-Laplacian of \cite{CoS}, is not invariant under a transitive group of isometries of the underlying manifold). For this reason, here we take the opportunity to establish an ``abstract'' multiplier theorem, which 
applies to a rather general setting of self-adjoint operators on  bounded metric measure spaces, satisfying the volume doubling property,
and extends the analogous result in \cite{CoS} to the framework of a family of scale-dependent weights.

\bigskip
It would be of great interest to establish whether  Theorem \ref{thm:main}
is sharp  
when $\DimK > \DimD/2$ or alternatively improve on it.
 The corresponding question for the Grushin operators $\opG_{\DimD,\DimK}$ on $\RR^\DimD$
has been settled in  \cite{MMu}; based on that result, one may expect that Theorem \ref{thm:main} and Corollary \ref{cor:main} actually hold with $D$ replaced by $d$.
However, when the dimension $k$ 
of the singular set is larger than the codimension, 
the approach developed in this paper, which is  
 based on a ``weighted Plancherel estimate with weights on the first layer'',
does not suffice 
to obtain such result and new methods (inspired, for instance, to those in \cite{MMu} and involving the ``second layer'' as well) appear to be necessary.

\bigskip
The paper is organised as follows.
In Section \ref{s:abstracttheorem} 
we state our abstract multiplier theorem, of which
 Theorem \ref{thm:main} will be a direct consequence; 
 in order not to burden the exposition, we postpone the proof
 of the  abstract theorem
  to an appendix (Section \ref{s:appendix}). 
 In  Section
\ref{s:prelim}
we introduce the spherical Laplacians and the Grushin operators on $\sfera^d$.
A precise estimate for the sub-Riemannian distance $\dist$ associated with the Grushin operator $\opL_{\DimD,\DimK}$ is also given.
 Moreover, we introduce a system of cylindrical coordinates on $\sfera^d$ 
 which is  key to   
 our approach.
In Section \ref{s:eigenfunctions}   
we recall
 the construction of a complete system of joint eigenfunctions of $\Delta_\DimD,\dots,\Delta_\DimK$ on $\sfera^\DimD$,
in terms of which we 
explicitly write down
the spectral decomposition of the Grushin operator $\opL_{\DimD,\DimK} = \Delta_\DimD - \Delta_\DimK$.
We also prove
some Riesz-type bounds for $\opL_{\DimD,\DimK}$, and 
we state the refined estimates for ultraspherical polynomials, which are the building blocks in the joint spectral decomposition.
Section \ref{s:weighted} is devoted to the proof of the crucial  ``weighted spectral cluster estimates'' for the Grushin operators $\opL_{\DimD,\DimK}$. 
In Section \ref{s:Plancherel} we use the Riesz-type bounds and the weighted spectral cluster estimates
to prove ``weighted Plancherel-type estimates'' for the Grushin operator $\opL_{\DimD,\DimK}$.
After this preparatory work,
the proof of Theorem \ref{thm:main}, which boils down to verifying the assumptions of the abstract theorem,
concludes the section. 

\bigskip

 Throughout the paper, 
for any two nonnegative quantities $X$ and $Y$,
we use $X \lesssim Y$   or $Y\gtrsim X$  to denote the estimate $X \le CY$  for a positive constant $C$.
The symbol $X\simeq Y$ is shorthand for $X\lesssim Y$ and
$Y \gtrsim X$. We use variants such as $\lesssim_k$ or $\simeq_k$ to indicate that the implicit constants may depend on the parameter $k$.

\section{An abstract multiplier theorem}\label{s:abstracttheorem}

We state an abstract multiplier theorem, which is a refinement of \cite[Theorem 3.6]{CoS} and \cite[Theorem 3.2]{DOS}. 
The proof of our main result, Theorem \ref{thm:main}, for the operator $\opL_{\DimD,\DimK}$ will follow from this result.

As in \cite{CoS,DOS}, for all $q \in [2,\infty]$, $N \in \NN \setminus \{0\}$ and $F : \RR \to \CC$ supported in $[0,1]$, we define the norm $\|F\|_{N,q}$ by
\[
\|F\|_{N,q} = \begin{cases}
\left(\frac{1}{N} \sum_{i=1}^N \sup_{\lambda \in [(i-1)/N,i/N]} |F(\lambda)|^q\right)^{1/q} &\text{if } q < \infty,\\
\sup_{\lambda \in [0,1]} |F(\lambda)| &\text{if } q = \infty.
\end{cases}
\]
Moreover, by $\Kern_T$ we denote the integral kernel of an operator $T$.

\begin{theorem}\label{thm:abstractmult}
Let $(X,\dist)$ be a bounded metric space, equipped with a regular Borel measure $\mu$. 
Let $\AbsOp$ a nonnegative self-adjoint operator on $L^2(X)$. Let $q \in [2,\infty]$. Suppose that there exist a family of weight functions $\absWeight_r \colon X\times X \to [0,\infty)$, where $r\in (0,1]$, and a constant $\absDim \in [1,\infty)$ such that the following conditions are satisfied:
\begin{enumerate}[label=(\alph*)]
\item\label{en:abs_doubling} the doubling condition:
\[
\mu(B(x,2r)) \lesssim \mu(B(x,r)) \qquad\forall x\in X  \quad\forall r>0;
\]
\item\label{en:abs_heatkernel} heat kernel bounds: 
\[
|\Kern_{\exp(-t\AbsOp)}(x,y)| \lesssim_N \mu(B(y,t^{1/2}))^{-1} (1+\dist(x,y)/t^{1/2})^{-N}
\]
for all $N \geq 0$, for all $t \in (0,\infty)$ and $x,y \in X$;
\item\label{en:abs_weightgrowth} the growth condition:
\begin{equation}\label{eq:weightsineq22}
C^{-1} \leq \absWeight_r(x,y) \lesssim (1+\dist(x,y)/r)^{M_0}
\end{equation}
for some $M_0 \geq 0$, for all $r \in (0,1]$  and $x,y \in X$;
\item\label{en:abs_weightint} the integrability condition:
\begin{equation}\label{eq:weightsint11}
\int_X (1+\dist(x,y)/r)^{-\beta} (\absWeight_r(x,y))^{-1} \,d\mu(x) \lesssim_{\beta} \mu (B(y,r))
\end{equation}
for all $r \in (0,1]$, $\beta >\absDim$ and for all $y \in X$;
\item\label{en:abs_plancherel} weighted Plancherel-type estimates:
\begin{equation}\label{eq:abstract_plancherel}
\esssup_{y \in X} \mu(B(y,1/N)) \int_{X} \absWeight_{1/N}(x,y) \, |\Kern_{F(\sqrt{\AbsOp})}(x,y)|^2 \,d\mu (x)
\lesssim \|F(N \cdot)\|_{N,q}^2
\end{equation}
for all $N \in \NN \setminus \{0\}$ and for all bounded Borel functions $F : \RR \to \CC$ supported in $[0,N]$.
\end{enumerate}
Finally, assume that $s>\absDim/2$. Then the following hold.
\begin{enumerate}[label=(\roman*)]
\item\label{en:abs_compact} For continuous functions $F : \RR \to \CC$ supported in $[-1,1]$,
\[
\sup_{t>0} \| F(t \sqrt{\AbsOp}) \|_{L^{1}(X) \to L^{1}(X)} \lesssim_s \, \|F\|_{\Sob{q}{s}}.
\]
\item\label{en:abs_mh} For all bounded Borel functions $F : \RR \to \CC$ continuous on $(0,\infty)$,
\begin{equation}\label{eq:main_mh_wt11}
\|F(\sqrt{\AbsOp}) \|_{L^1(X) \to L^{1,\infty}(X)} \lesssim_s \, \sup_{t \geq 0} \| F(t \cdot) \, \eta \|_{\Sob{q}{s}}.
\end{equation}
Hence, whenever the right-hand side of \eqref{eq:main_mh_wt11} is finite, the operator 
$F(\sqrt{\AbsOp})$ is of weak type $(1,1) $ and bounded on $L^p(X)$ for all $p \in (1,\infty)$.
\end{enumerate}
\end{theorem}

Since the subject is replete with technicalities, which could weigh on the discussion, we defer the proof of the abstract theorem to an appendix (Section \ref{s:appendix}).

Let us just observe that Assumption \ref{en:abs_heatkernel} only requires a polynomial decay in space (of arbitrary large order) for the heat kernel; hence this assumption is weaker than the corresponding ones in \cite{DOS}, where Gaussian-type (i.e., superexponential) decay is required, and in \cite{CoS}, where finite propagation speed for the associated wave equation is required (which, under the ``on-diagonal bound'' implied by \eqref{eq:abstract_plancherel}, is equivalent to ``second order'' Gaussian-type decay \cite{S}), and matches instead the assumption in \cite{He} (see also \cite[Section 6]{M}).

Another important feature of the above result, which is crucial for the applicability to the spherical Grushin operators $\opL_{\DimD,\DimK}$ considered in this paper, is the use of a family of weight functions, where the weight $\absWeight_r$ may depend on the scale $r$ in a nontrivial way; this constitutes another important difference to \cite{CoS}, where the weights considered are effectively scalar multiples of a single weight function (compare Assumptions \ref{en:abs_weightint} and \ref{en:abs_plancherel} above with \cite[Assumptions 2.2 and 2.5]{CoS}).

The attentive reader will have noticed that it is actually enough to verify Assumptions \ref{en:abs_weightgrowth} and \ref{en:abs_weightint} for scales $r=1/N$ for $N \in \NN \setminus \{0\}$ (indeed, one can redefine $\pi_r$ as $\pi_{1/\lfloor 1/r\rfloor}$ when $1/r \notin \NN$); the slightly redundant form of the above assumptions is just due to notational convenience.

\section{Spherical Laplacians and Grushin operators}\label{s:prelim}

\subsection{The Laplace--Beltrami operator on the unit sphere}
For $\DimD \in \NN$, $\DimD\geq 1$, let $\sfera^\DimD$ denote the unit sphere in $\RR^{1+\DimD}$, as in \eqref{eq:sfera}.
The Euclidean structure on $\RR^{1+\DimD}$ induces a natural, rotation-invariant Riemannian structure on $\sfera^\DimD$. Let $\meas$ denote the corresponding Riemannian measure, and $\Delta_\DimD$ the Laplace--Beltrami operator on the unit sphere $\sfera^\DimD$ in $\RR^{1+\DimD}$. It is possible (see, e.g., \cite{Geller}) to give a more explicit expression for $\Delta_\DimD$, namely,
\begin{equation}\label{eq:Deltad_standard}
\Delta_\DimD = - \sum_{0 \leq j < r \leq \DimD} Z_{j,r}^2,
\end{equation}
where
\[
Z_{j,r} = z_j \frac{\partial}{\partial z_r} - z_r \frac{\partial}{\partial z_j}.
\]
Indeed the rotation group $\SO(1+\DimD)$ acts naturally on $\RR^{1+\DimD}$ and $\sfera^\DimD$; via this action, the vector fields $Z_{j,r}$ ($0 \leq j < r \leq \DimD$) correspond to the standard basis of the Lie algebra of $\SO(1+\DimD)$, and $\Delta_\DimD$ corresponds to the Casimir operator. The commutation relations
\begin{equation}\label{eq:commutation}
[Z_{j,r}, Z_{j',r'}] = \delta_{r,j'} Z_{j,r'} + \delta_{j,r'} Z_{r,j'} - \delta_{j,j'} Z_{r,r'} - \delta_{r,r'} Z_{j,j'}
\end{equation}
are easily checked and correspond to those of the Lie algebra of $\SO(1+\DimD)$.

\subsection{A family of commuting Laplacians and spherical Grushin operators}

By \eqref{eq:commutation}, the operator $\Delta_\DimD$ commutes with all the vector fields $Z_{j,r}$ (this corresponds to the fact that the Casimir operator is in the centre of the universal enveloping algebra of the Lie algebra of $\SO(1+\DimD)$); in particular it commutes with each of the ``partial Laplacians''
\begin{equation}\label{eq:partial_lap}
\Delta_r = -\sum_{0 \leq j < s \leq r} Z_{j,s}^2
\end{equation}
for $r = 1,\dots,\DimD$.

Assume that $\DimD\geq 2$. We now observe that, for $r=1,\dots,\DimD-1$, we can identify $\SO(1+r)$ with a subgroup of $\SO(1+\DimD)$, by associating to each $A$ in $\SO(1+r)$ the element
\[
\begin{pmatrix} A & 0 \\ 0 & I \end{pmatrix}
\]
of $\SO(1+\DimD)$. Via this identification, the operator $\Delta_r$ corresponds to the Casimir operator of $\SO(1+r)$, and therefore it commutes with all the operators $\Delta_s$ for $s=1,\dots,r$.

In conclusion, the operators $\Delta_1,\dots,\Delta_\DimD$ commute pairwise, and admit a joint spectral decomposition. In what follows we will be interested in the study of the Grushin-type operator
\begin{equation}\label{eq:grushin}
\opL_{\DimD,\DimK} = \Delta_\DimD - \Delta_\DimK = -\sum_{r=\DimK+1}^\DimD \sum_{j=0}^{r-1} Z_{j,r}^2.
\end{equation}
for $\DimK=1,\dots,\DimD-1$.

The operator $\opL_{\DimD,\DimK}$ is not uniformly elliptic: indeed it degenerates on the $\DimK$-submanifold $\equat_{\DimD,\DimK} = \sfera^\DimK \times \{0\}$ of $\sfera^\DimD$. More precisely, if $\fZ_{\DimD,\DimK} = \{ Z_{j,r} \tc \DimK+1 \leq r \leq \DimD, \, 0 \leq j < r \}$ is the family of vector fields appearing in the sum \eqref{eq:grushin}, then
it is easily checked that, for all $z \in \sfera^\DimD$,
\begin{equation}\label{eq:horizontal}
\Hz^{\DimD,\DimK}_z \defeq \Span \{ X|_z \tc X \in \fZ_{\DimD,\DimK} \} = \begin{cases}
T_z \sfera^\DimD &\text{if } z \notin \equat_{\DimD,\DimK},\\
(T_z \equat_{\DimD,\DimK})^\perp &\text{if } z \in \equat_{\DimD,\DimK}.
\end{cases}
\end{equation}
On the other hand, the commutation relations \eqref{eq:commutation} give that
\[
[Z_{j,\DimD},Z_{j',\DimD}] = - Z_{j,j'}
\]
for all $j,j'=0,\dots,\DimD-1$; in particular the vector fields in $\fZ_{\DimD,\DimK}$, together with their Lie brackets, span the tangent space of $\sfera^\DimD$ at each point. In other words, the family of vector fields $\fZ_{\DimD,\DimK}$ satisfies H\"ormander's condition and (together with the Riemannian measure $\meas$) determines a (non-equiregular) $2$-step sub-Riemannian structure on $\sfera^\DimD$ with the horizontal distribution $\Hz^{\DimD,\DimK}$ described in \eqref{eq:horizontal}. The corresponding sub-Riemannian norm on the fibres of $\Hz^{\DimD,\DimK}$ is given, for all $p \in \sfera^\DimD$ and $v \in \Hz^{\DimD,\DimK}_p$, by
\begin{equation}\label{eq:sr_metric}
|v|_{\opL_{\DimD,\DimK}} = \inf \left\{ \sqrt{\sum_{X \in \fZ_{\DimD,\DimK}} a_X^2} \tc (a_X)_X \in \RR^{\fZ_{\DimD,\DimK}}, \ v = \sum_{X \in \fZ_{\DimD,\DimK}} a_X X|_p \right\}.
\end{equation}
For more details on sub-Riemannian geometry we refer the reader to \cite{ABB,BellaicheRisler,CalinChang,Montgomery}.

\subsection{Cylindrical coordinates}\label{ss:cylcoords}
In order to study the operator $\opL_{\DimD,\DimK}$, it is useful to introduce a system of ``cylindrical coordinates'' on $\sfera^\DimD$ that will provide a particularly revealing expression for $\opL_{\DimD,\DimK}$ in a neighbourhood of the singular set $\equat_{\DimD,\DimK}$.

For all $\omega \in \sfera^{\DimD-1}$ and $\psi \in [-\pi/2, \pi/2]$, let us define the point $\spnt{\omega}{\psi} \in \sfera^\DimD$ by
\begin{equation}\label{eq:coordsfera}
\spnt{\omega}{\psi} = ((\cos\psi)\omega, \sin\psi).
\end{equation}
Away from $\psi = \pm \pi/2$, the map $(\omega,\psi) \mapsto \spnt{\omega}{\psi}$ is a diffeomorphism onto its image, which is the sphere without the two poles; so \eqref{eq:coordsfera} can be thought of as a ``system of coordinates'' on $\sfera^\DimD$, up to null sets.
In these coordinates, the spherical measure $\sigma$ on $\sfera^\DimD$ is given by
\[
d\meas(\spnt{\omega}{\psi}) = \cos^{\DimD-1}\psi \,d\psi \,d\meas_{\DimD-1}(\omega),
\]
where $\meas_{\DimD-1}$ is the spherical measure on $\sfera^{\DimD-1}$.
Moreover, the Laplace--Beltrami operator may be written in these coordinates as
\begin{equation}\label{eq:def_Delta_d}
\Delta_\DimD 
= -\frac{1}{\cos^{\DimD-1}\psi} \frac{\partial}{\partial \psi} \cos^{\DimD-1}\psi \frac{\partial}{\partial \psi}
+ \frac{1}{\cos^2\psi} \Delta_{\DimD-1} ,
\end{equation}
where $\Delta_{\DimD-1}$, given by \eqref{eq:partial_lap}, corresponds to
the Laplace--Beltrami operator on $\sfera^{\DimD-1}$ (see, e.g., \cite[\S IX.5]{Vilenkin}).

We now iterate the previous construction.
Let $\DimK \in \NN$ such that $1 \leq \DimK < \DimD$ be fixed.
Starting from \eqref{eq:coordsfera}, we can inductively define the point
\begin{equation}\label{eq:coordinates_induct}
\spnt{\omega}{\psi}
=
\spnt{\ldots\spnt{\spnt{\omega}{\psi_{\DimK+1}}}{\psi_{\DimK+2}}\ldots}{\psi_\DimD}
\end{equation}
of $\sfera^\DimD$ for all $\psi = (\psi_{\DimK+1},\ldots,\psi_\DimD) \in [-\pi/2,\pi/2]^{\DimD-\DimK}$ and $\omega \in \sfera^\DimK$; if we restrict $\psi$ to $(-\pi/2,\pi/2)^{\DimD-\DimK}$, then \eqref{eq:coordinates_induct} defines a ``system of coordinates'' for an open subset $\Omega_{\DimD,\DimK}$ of $\sfera^\DimD$ of full measure, namely,
\[
\Omega_{\DimD,\DimK} = \{ \spnt{\omega}{\psi} \tc \omega \in \sfera^\DimK, \, \psi \in (-\pi/2,\pi/2)^{\DimD-\DimK} \}.
\] 
In these coordinates, the spherical measure $\meas$ on $\sfera^\DimD$ is given by
\begin{equation}\label{eq:spherical_measure}
d\meas(\spnt{\omega}{\psi}) = 
\cos^{\DimD-1} \psi_\DimD \cdots \cos^\DimK \psi_{\DimK+1} \,d\psi_\DimD \cdots \,d\psi_{\DimK+1} \,d\meas_\DimK(\omega),
\end{equation}
where $\meas_\DimK$ is the spherical measure on $\sfera^\DimK$. Moreover, starting from \eqref{eq:def_Delta_d},
we get inductively that
\begin{multline*}
\Delta_\DimD = -\sum_{r=\DimK+1}^\DimD \frac{1}{\cos^2 \psi_{r+1} \cdots \cos^2 \psi_\DimD}
\frac{1}{\cos^{r-1} \psi_r} \frac{\partial}{\partial \psi_r} \cos^{r-1} \psi_r \frac{\partial}{\partial \psi_r} \\
+ \frac{1}{\cos^{2}\psi_{\DimK+1} \cdots \cos^{2}\psi_\DimD}  \Delta_\DimK,
\end{multline*}
where again $\Delta_\DimK$ is the operator given by \eqref{eq:partial_lap}.

In particular, the Grushin operator $\opL_{\DimD,\DimK} = \Delta_\DimD - \Delta_\DimK$ on $\sfera^\DimD$ may be written in these coordinates as
in \eqref{eq:Ldk-fields},
where the vector fields $Y_r$ 
and the function $\potV : (-\pi/2,\pi/2)^{\DimD-\DimK} \to \RR$ are defined by \eqref{eq:vfield_Yr} and \eqref{eq:potential} respectively.
Note that $\potV(\psi)$ vanishes only for $\psi = 0$, corresponding to the singular set $\equat_{\DimD,\DimK}$.
We also remark that
\begin{equation}\label{eq:potV_approx}
\frac{1}{\cos \psi_{r+1} \cdots \cos \psi_\DimD} \simeq 1, \qquad \potV(\psi) \simeq |\psi|^2
\end{equation}
for $r=\DimK+1,\dots,\DimD$, uniformly for $|\psi| \leq \varepsilon$, for any given $\varepsilon \in (0,\pi/2)$.

The formula \eqref{eq:Ldk-fields} for the sub-Laplacian corresponds to a somewhat more explicit expression for the sub-Riemannian norm \eqref{eq:sr_metric} on the fibres of the horizontal distribution, which is better written by identifying, via the ``coordinates''  \eqref{eq:coordinates_induct}, the tangent space $T_{\spnt{\omega}{\psi}} \sfera^\DimD$ with $T_\omega \sfera^\DimK \times \RR^{\DimD-\DimK}$ for all $\spnt{\omega}{\psi} \in \Omega_{\DimD,\DimK}$. Under this identification,
\begin{equation}\label{eq:horiz_expl}
\Hz^{\DimD,\DimK}_{\spnt{\omega}{\psi}} = \begin{cases}
T_\omega \sfera^\DimK \times \RR^{\DimD-\DimK} &\text{if } \psi \neq 0,\\
\{0\} \times \RR^{\DimD-\DimK} &\text{if } \psi = 0
\end{cases}
\end{equation}
and, for all $(v,w) \in \Hz^{\DimD,\DimK}_{\spnt{\omega}{\psi}}$, its sub-Riemannian norm satisfies
\begin{equation}\label{eq:subriem_norm_expl}
|(v,w)|_{\opL_{\DimD,\DimK}}^2 = \begin{cases}
\sum_{r=\DimK+1}^\DimD (\cos \psi_{r+1} \cdots \cos \psi_\DimD)^2 |w_r|^2 + \potV(\psi)^{-1} |v|^2 &\text{if } \psi \neq 0,\\
\sum_{r=\DimK+1}^\DimD (\cos \psi_{r+1} \cdots \cos \psi_\DimD)^2 |w_r|^2 &\text{if } \psi = 0,
\end{cases}
\end{equation}
where $w = (w_{\DimK+1},\dots,w_\DimD) \in \RR^{\DimD-\DimK}$ and $|v|$ is the Riemannian norm of $v \in T_\omega \sfera^\DimK$.

\subsection{The sub-Riemannian distance}
Thanks to \eqref{eq:subriem_norm_expl}, we can obtain a precise estimate for the sub-Riemannian distance $\dist$ associated with the Grushin operator $\opL_{\DimD,\DimK}$.
This is the analogue of \cite[Proposition 5.1]{RS}, that treats the case of ``flat'' Grushin operators on $\RR^n$, and \cite[Proposition 2.1]{CaCiaMa}, that treats the case of $\opL_{2,1}$ on $\sfera^2$.
In the statement below we represent the points of the sphere in the form $\spnt{\omega}{\psi}$ for $\omega \in \sfera_\DimK$, $\psi \in [-\pi/2,\pi/2]^{\DimD-\DimK}$, as in \eqref{eq:coordinates_induct}. We also denote by $\dist_{R,\sfera^\DimK}$ and $\dist_{R,\sfera^\DimD}$ the Riemannian distances on the spheres $\sfera^\DimK$ and $\sfera^\DimD$.

\begin{proposition}\label{prp:subriemannian}
Let $\varepsilon \in (0,\pi/2)$. The sub-Riemannian distance $\dist$ on $\sfera^\DimD$ associated with $\opL_{\DimD,\DimK}$ satisfies
\begin{equation}\label{eq:rho-dist}
\dist (\spnt{\omega}{\psi},\spnt{\omega'}{\psi'}) \simeq |\psi-\psi'| + \min\left\{{\dist_{R,\sfera^\DimK}(\omega,\omega')}^{1/2}, \frac{\dist_{R,\sfera^\DimK}(\omega,\omega')}{\max\{|\psi|,| \psi'|\}}\right\},
\end{equation}
if $\max\{|\psi|,|\psi'|\} \leq \varepsilon$;
if instead $\max\{|\psi|,|\psi'|\} \geq \varepsilon$, then
\begin{equation}\label{eq:rho-dist-far}
\dist (\spnt{\omega}{\psi},\spnt{\omega'}{\psi'}) 
\simeq \dist_{R,\sfera^\DimD} (\spnt{\omega}{\psi},\spnt{\omega'}{\psi'}).
\end{equation}
Consequently, the $\meas$-measure $V(\spnt{\omega}{\psi},r)$ of the $\dist$-ball centred at $\spnt{\omega}{\psi}$ with radius $r \geq 0$ satisfies 
\begin{equation}\label{eq:volume}
V(\spnt{\omega}{\psi}, r) \simeq \min\{ 1, r^\DimD \max \{r, |\psi| \}^\DimK \}.
\end{equation}
The implicit constants may depend on $\varepsilon$.
\end{proposition}
\begin{proof}
Note that the sub-Riemannian distance $\dist$ and the Riemannian distance $\dist_{R,\sfera^\DimD}$ are locally equivalent far from the singular set $\equat_{\DimD,\DimK}$: since $\Hz^{\DimD,\DimK}_p = T_p M$ for all $p \in \sfera^\DimD \setminus \equat_{\DimD,\DimK}$ (see \eqref{eq:horizontal}), and the Riemannian and sub-Riemannian inner products on $T_p M$ depend continuously on $p$, it is enough to apply \cite[Lemma 2.3]{CaCiaMa} by choosing as $M$ and $N$
 the Riemannian and sub-Riemannian $\sfera^\DimD$ respectively, and as $F$ the identity map restricted to any open subset $U$ of $\sfera$ with compact closure not intersecting $\equat$. 
 Then   \cite[Lemma 2.2]{CaCiaMa}, applied with $K = \{ (\spnt{\omega}{\psi},\spnt{\omega'}{\psi'}) \in \sfera^\DimD \times \sfera^\DimD \tc \max \{|\psi|,|\psi'|\} \geq \varepsilon \}$,
yields  \eqref{eq:rho-dist-far}.

Note now that the expression in the right-hand side of \eqref{eq:rho-dist} defines a continuous function 
$\Phi : \Omega_{\DimD,\DimK} \times \Omega_{\DimD,\DimK} \to [0,\infty)$, 
which is nondegenerate in the sense of  \cite[Lemma 2.2]{CaCiaMa}.
Hence, in order to prove the equivalence 
\eqref{eq:rho-dist}, it is enough to show that $\Phi$ and $\dist$ are locally equivalent at each point $p_0 \in \Omega_{\DimD,\DimK}$, and then apply \cite[Lemma 2.2]{CaCiaMa} with $K = \{ (\spnt{\omega}{\psi},\spnt{\omega'}{\psi'}) \in \sfera^\DimD \times \sfera^\DimD \tc \max\{|\psi|,|\psi'|\} \leq \varepsilon \}$.

Consider now the Grushin operator $\opG = \opG_{\DimD,\DimK}$
on $\RR^{\DimD-\DimK}_x \times \RR^\DimK_y$ defined in \eqref{def:Grushin}. The associated horizontal distribution $\Hz^\opG$ and sub-Riemannian metric are given by
\begin{equation}\label{eq:flat_horizontal}
\Hz^\opG_{(x,y)} = \begin{cases}
\RR^{\DimD-\DimK} \times \RR^\DimK &\text{if } x \neq 0,\\
\RR^{\DimD-\DimK} \times \{0\} &\text{if } x = 0,
\end{cases}
\qquad
|(w,v)|_{\opG}^2 = \begin{cases}
|w|^2 + |x|^{-2} |v|^2 &\text{if } x \neq 0,\\
|w|^2 &\text{if } x = 0,
\end{cases}
\end{equation}
for all $(x,y) \in \RR^{\DimD-\DimK} \times \RR^\DimK$ and $(w,v) \in \Hz^\opG_{(x,y)}$. Moreover, according to \cite[Proposition 5.1]{RS}, the associated sub-Riemannian distance $\dist_\opG$ satisfies
\begin{equation}\label{eq:flat_distance}
\dist_\opG((x,y),(x',y')) \simeq |x-x'| + \min\left\{ |y-y'|^{1/2}, \frac{|y-y'|}{\max\{|x|,|x'|\}}\right\}.
\end{equation}

Let $p_0 = \spnt{\omega_0}{\psi_0} \in \Omega_{\DimD,\DimK}$. Choose coordinates for $\sfera^\DimK$ centred at $\omega_0$, thus determining a diffeomorphism $f$ from an open neighbourhood $A$ of $0$ in $\RR^\DimK$ to a neighbourhood $f(A)$ of $p$ in $\sfera^\DimK$. By the equivalence of norms, up to shrinking $A$, we may assume that
\begin{equation}\label{eq:normequivalence_Rk}
|df_y(w)| \simeq |w|
\end{equation}
for all $y \in A$, $w \in \RR^\DimK \cong T_y \RR^\DimK$, where the norms in \eqref{eq:normequivalence_Rk} are those determined by the Riemannian structures of $\sfera^\DimK$ and $\RR^\DimK$; similarly, we may also assume that
\begin{equation}\label{eq:distequivalence_Rk}
\dist_{R,\sfera^\DimK}(f(y),f(y')) \simeq |y-y'|
\end{equation}
for all $y,y' \in A$. Let now $U = B \times A$, where $B$ is a neighbourhood of $\psi_0$ with compact closure in $(-\pi/2,\pi/2)^{\DimD-\DimK}$, and define $F : U \to \sfera^\DimD$ by $F(x,y) = \spnt{f(y)}{x}$. A comparison of \eqref{eq:horiz_expl} and \eqref{eq:subriem_norm_expl} with \eqref{eq:flat_horizontal}, taking \eqref{eq:potV_approx} and \eqref{eq:normequivalence_Rk} into account, immediately shows that \cite[Lemma 2.3]{CaCiaMa} can be applied to the map $F$ and the sub-Riemannian structures associated with $\opG$ and $\opL_{\DimD,\DimK}$; consequently, up to shrinking $U$, we obtain that
\[
\dist(F(p),F(p')) \simeq \dist_\opG(p,p') \simeq \Phi(F(p),F(p'))
\]
for all $p,p' \in U$, where the latter equivalence readily follows from \eqref{eq:flat_distance} and \eqref{eq:distequivalence_Rk}.

Finally, the estimate \eqref{eq:volume} for the volume of balls follows from \eqref{eq:spherical_measure}, \eqref{eq:rho-dist} and \eqref{eq:rho-dist-far} by considering separately the cases $|\psi|$ small and $|\psi|$ large.
\end{proof}

\section{A complete system of joint eigenfunctions}\label{s:eigenfunctions}

Let $\DimD,\DimK \in \NN$ with $1 \leq \DimK < \DimD$.
In this section we briefly recall the construction of a complete system of joint eigenfunctions of $\Delta_\DimD,\dots,\Delta_\DimK$ on $\sfera^\DimD$.
This will give in particular the spectral decomposition of the Grushin operator $\opL_{\DimD,\DimK} = \Delta_\DimD - \Delta_\DimK$.

This construction is classical and can be found in several places in the literature (see, e.g., \cite[Ch.\ IX]{Vilenkin} or \cite[Ch.\ XI]{EMOT}), where explicit formulas for spherical harmonics on spheres of arbitrary dimension are given, in terms of ultraspherical (Gegenbauer) polynomials. The discussion below is essentially meant to fix the notation that will be used later.

By the symbol $P^{(\alpha,\beta)}_j$ we shall denote the Jacobi polynomial of degree $j \in \NN$ and indices $\alpha,\beta>-1$, defined by means of Rodrigues' formula:
\begin{equation}\label{eq:jacobi}
P^{(\alpha,\beta)}_j(x)= \frac{(-1)^j}{2^j\, j!} (1-x)^{-\alpha} (1+x)^{-\beta} \left(\frac{d}{dx}\right)^j \left((1-x)^{\alpha+j} (1+x)^{\beta+j} \right)
\end{equation}
for $x \in (-1,1)$.
We recall, in particular, the symmetry relation 
\begin{equation}\label{eq:symmetry}
P_j^{(\alpha,\beta)}(x) = (-1)^j P_j^{(\beta,\alpha)}(x),
\end{equation}
for $j \in \NN$, $\alpha,\beta>-1$ and $x \in \RR$. Ultraspherical polynomials correspond to Jacobi polynomials with $\alpha=\beta$ \cite[(4.7.1)]{Szego}.

\subsection{Spectral theory of the Laplace--Beltrami operator}

We first recall some well known facts about the spectral theory of $\Delta_\DimD$ (see, e.g., \cite[Ch.\ 4]{Stein-Weiss} or \cite[Ch.\ 5]{AxBR}). The operator $\Delta_\DimD$ is essentially self-adjoint on $L^2(\sfera^\DimD)$ and has discrete spectrum: its eigenvalues are given by
\begin{equation}\label{eq:laplace_eigenvalues}
\lambda_\ell^\DimD \defeq ( \ell+(\DimD-1)/2 ) ( \ell-(\DimD-1)/2 ),
\end{equation}
where $\ell\in\NN_\DimD$, and
\begin{equation}\label{eq:Nd}
\NN_\DimD = \NN + (\DimD-1)/2.
\end{equation}
The corresponding eigenspaces, denoted by $\Harm^\ell(\sfera^\DimD)$, consist of all spherical harmonics of degree $\ell' = \ell-(\DimD-1)/2$, that is, of all restrictions to $\sfera^\DimD$ of homogeneous harmonic polynomials on $\RR^{1+\DimD}$ of degree $\ell'$;
they are finite-dimensional spaces of dimension
\begin{equation}\label{eq:dimensione}
\dimHarm_\ell(\sfera^\DimD)=\binom{\ell'+\DimD}{\ell'}-\binom{\ell'+\DimD-2}{\ell'-2} 
= \frac{2\ell'+\DimD-1}{\DimD-1} \binom{\ell'+\DimD-2}{\DimD-2}
\end{equation}
for $\ell\in\NN_\DimD$ (the last identity only makes sense when $\DimD>1$),
and in particular
\begin{equation}\label{eq:dimensione_stima}
\dimHarm_\ell(\sfera^\DimD) \simeq_\DimD \ell^{\DimD-1}
\end{equation}
(this estimate is also valid when $\DimD=1$, provided we stipulate that $0^0 = 1$).

Since $\Delta_\DimD$ is self-adjoint, its eigenspaces are mutually orthogonal, i.e.,
\[
\Harm^\ell(\sfera^\DimD) \perp \Harm^{\ell'}(\sfera^\DimD)
\]
for $\ell,\ell' \in \NN_\DimD$, $\ell \neq \ell'$. Moreover, if $E^\DimD_\ell$ is an orthonormal basis of $\Harm^\ell(\sfera^\DimD)$,
then
\begin{equation}\label{eq:SW-dim}
\sum_{Z \in E^\DimD_\ell} | Z(z) |^2 = \meas(\sfera^\DimD)^{-1} \, \dimHarm_\ell(\sfera^\DimD)
\end{equation}
for all $z \in \sfera^\DimD$ \cite[Ch.\ 4, Corollary 2.9]{Stein-Weiss}.

\subsection{Joint eigenfunctions of \texorpdfstring{$\Delta_\DimD$}{Delta d} and \texorpdfstring{$\Delta_{\DimD-1}$}{{Delta d-1}}}
We start the construction of joint eigenfunctions with the case $\DimK = \DimD-1$, and look for eigenfunctions of $\Delta_\DimD$ 
that are simultaneously eigenfunctions of $\Delta_{\DimD-1}$.

Following, e.g., \cite[\S IX.5]{Vilenkin}, one can use the expression \eqref{eq:def_Delta_d} for $\Delta_\DimD$ to solve the eigenfunction equation for $\Delta_\DimD$ via separation of variables. More precisely, we look for functions $W$ on $\sfera^\DimD$ of the form $X \otimes Z$, that is,
\[
W(\spnt{\omega}{\psi}) = X(\psi) Z(\omega)
\]
in the coordinates \eqref{eq:coordsfera}, such that
\[
\Delta_\DimD W = \lambda^\DimD_\ell W \qquad\text{and}\qquad \Delta_{\DimD-1} Z = \lambda^\DimK_m Z,
\]
for some $\ell \in \NN_\DimD$, $m \in \NN_{\DimD-1}$. This leads to a differential equation for $X$ that is solved in terms of ultraspherical polynomials. Namely, if $Z \in \Harm^m(\sfera^{\DimD-1})$ is nonzero and $\ell \geq m$, then $W = X \otimes Z$ is in $\Harm^\ell(\sfera^\DimD)$ if and only if $X$ is a multiple of
\begin{equation}\label{eq:Xlm_def_iniziale}
X_{\ell,m}^\DimD(\psi) = c_{\ell m} (\cos \psi)^{m-(\DimD-2)/2} P_{\ell-m-1/2}^{(m,m)}(\sin \psi).
\end{equation}
Here the normalization constant  $c_{\ell m}$ is chosen so that
\begin{equation}\label{eq:Xlm_def_iniziale-norm}
\int_{-\pi/2}^{\pi/2}
|X_{\ell, m}^\DimD (\psi)|^2 \,\cos^{\DimD-1}\psi \,d\psi = 1,
\end{equation}
that is, by means of \cite[(4.3.3)]{Szego},
\begin{equation}\label{eq:clm}
c_{\ell m} = \frac{\bigl[\ell \, \Gamma ( \ell-m+1/2)\, \Gamma (\ell +m+1/2)\bigr]^{1/2}}{2^m \,\Gamma( \ell+1/2)}.
\end{equation}

Define
\begin{equation}\label{eq:Id}
I_\DimD = \{ (\ell,m) \tc \ell \in \NN_\DimD, \, m \in \NN_{\DimD-1}, \,  \ell \geq m \}.
\end{equation}
Then, for all $(\ell, m)\in I_\DimD$, we obtain an injective linear map
\[
\Harm^m(\sfera^{\DimD-1}) \ni Z \mapsto X^\DimD_{\ell,m} \otimes Z \in \Harm^\ell(\sfera^\DimD) ,
\]
which is an isometry with respect to the Hilbert space structures of $L^2(\sfera^{\DimD-1})$ and $L^2(\sfera^\DimD)$, and a decomposition
\begin{equation}\label{eq:decomp-Hl}
\Harm^\ell(\sfera^\DimD) = \bigoplus_{\substack{
m \in \NN_{\DimD-1}\\ 
m \leq \ell
}}
X^\DimD_{\ell,m} \otimes \Harm^m(\sfera^{\DimD-1}) 
\end{equation}
(cf.\ \cite[p.\ 466, eq.\ (1)]{Vilenkin}). The summands in the right-hand side of \eqref{eq:decomp-Hl} are joint eigenspaces of $\Delta_\DimD$ and $\Delta_{\DimD-1}$ of eigenvalues $\lambda^\DimD_\ell$ and $\lambda^\DimK_m$ respectively; hence they are pairwise orthogonal in $L^2(\sfera^\DimD)$.

\subsection{Joint eigenfunctions of \texorpdfstring{$\Delta_\DimD,\dots,\Delta_\DimK$}{Delta d,...,Delta k}}\label{ss:dk-induction}
We go back to the general case $1 \leq \DimK < \DimD$ and we look for a complete system of joint eigenfunctions of $\Delta_\DimD,\dots,\Delta_\DimK$.

It is natural to introduce the index set
\[
J_\DimD^{(\DimK)} = \{ (\ell_\DimD,\ell_{\DimD-1},\dots,\ell_\DimK) \in \NN_\DimD \times \NN_{\DimD-1} \times \dots \times \NN_\DimK \tc \ell_\DimD \geq \ell_{\DimD-1} \geq \dots \geq \ell_\DimK \}.
\]
For all $(\ell_\DimD,\dots,\ell_\DimK) \in J_\DimD^{(\DimK)}$, let us define $X^\DimD_{\ell_\DimD,\dots,\ell_\DimK} : [-\pi/2,\pi/2]^{\DimD-\DimK} \to \RR$ by
\[
X^\DimD_{\ell_\DimD,\dots,\ell_\DimK}(\psi)
=
X^\DimD_{\ell_\DimD,\ell_{\DimD-1}}(\psi_\DimD)
\cdots
X^{\DimK+1}_{\ell_{\DimK+1},\ell_\DimK}(\psi_{\DimK+1}),
\]
where $\psi = (\psi_{\DimK+1},\dots,\psi_\DimD)$, and the functions $X^r_{\ell_r,\ell_{r-1}}$ are defined in \eqref{eq:Xlm_def_iniziale}.
Then, for all $(\ell_\DimD,\dots,\ell_\DimK) \in J_\DimD^{(\DimK)}$ and $Z \in \Harm^{\ell_\DimK}(\sfera^\DimK)$, the function $X^\DimD_{\ell_\DimD,\dots,\ell_\DimK} \otimes Z$, defined, in the coordinates \eqref{eq:coordinates_induct} on $\sfera^\DimD$, by
\begin{equation}\label{def:common_eigfncts_kd}
X^\DimD_{\ell_\DimD,\ldots,\ell_\DimK} \otimes Z : \spnt{\omega}{\psi} \mapsto X^\DimD_{\ell_\DimD,\dots,\ell_\DimK}(\psi) Z(\omega),
\end{equation}
is an eigenfunction of $\Delta_\DimD,\dots,\Delta_\DimK$ of respective eigenvalues $\lambda^\DimD_{\ell_\DimD},\dots,\lambda^\DimK_{\ell_\DimK}$. More precisely, iterating \eqref{eq:decomp-Hl}, we obtain the orthogonal direct sum decomposition
\[
\Harm^\ell(\sfera^\DimD) = \bigoplus_{\substack{
(\ell_\DimD,\dots,\ell_\DimK) \in J^{(\DimK)}_\DimD \\ \ell_\DimD = \ell
}}
X^\DimD_{\ell_\DimD,\dots,\ell_\DimK} \otimes \Harm^{\ell_\DimK}(\sfera^\DimK) .
\]

As a consequence, each function $f\in L^2(\sfera^\DimD)$ may be written as
\begin{equation}\label{eq:orth_decomp}
f = \sum_{(\ell_\DimD,\dots,\ell_\DimK) \in J_\DimD^{(\DimK)}} \sum_{Z \in E^\DimK_m} c_{\ell_\DimD,\ldots,\ell_\DimK,Z} X^\DimD_{\ell_\DimD,\ldots,\ell_\DimK} \otimes Z ,
\end{equation}
where $E^\DimK_m$ is an orthonormal basis of $\Harm^m(\sfera^\DimK)$ and
\[
c_{\ell_\DimD,\ldots,\ell_\DimK,Z} = \langle f, X^\DimD_{\ell_\DimD,\ldots,\ell_\DimK} \otimes Z \rangle.
\]
In particular, for all $(\ell_\DimD,\dots,\ell_\DimK) \in J_\DimD^{(\DimK)}$, the orthogonal projection $\pi^\DimD_{\ell_\DimD,\dots,\ell_\DimK}$ of $L^2 (\sfera^\DimD)$ onto the joint eigenspace of $\Delta_\DimD,\dots,\Delta_\DimK$ of eigenvalues $\lambda^\DimD_{\ell_\DimD},\dots,\lambda^\DimK_{\ell_\DimK}$
is given by
\begin{equation}\label{eq:def_joint_orthog_projection}
\pi^\DimD_{\ell_\DimD,\dots,\ell_\DimK} : f \mapsto \sum_{Z \in E^\DimK_m} \langle f,  X^\DimD_{\ell_\DimD,\ldots,\ell_\DimK} \otimes Z  \rangle X^\DimD_{\ell_\DimD,\ldots,\ell_\DimK} \otimes Z .
\end{equation}
Consequently, the integral kernel $K_{\ell_\DimD,\dots,\ell_\DimK}^\DimD$ of $\pi^\DimD_{\ell_\DimD,\dots,\ell_\DimK}$ is given by
\begin{equation}\label{eq:def_int_kernel_kd}
\begin{split}
K_{\ell_\DimD,\dots,\ell_\DimK}^\DimD(\spnt{\omega}{\psi},\spnt{\omega'}{\psi'})
&=
K_m^\DimK(\omega,\omega') X^\DimD_{\ell_\DimD,\ldots,\ell_\DimK}(\psi) \, X^\DimD_{\ell_\DimD,\ldots,\ell_\DimK}(\psi'),
\end{split}
\end{equation}
where
\[
K_m^\DimK(\omega,\omega') 
\defeq
\sum_{Z \in E^\DimK_m}
Z(\omega) \overline{Z(\omega')}
\]
is the integral kernel of the orthogonal projection of $L^2(\sfera^\DimK)$ onto $\Harm^m(\sfera^\DimK)$.

For all bounded Borel functions $F : \RR^{\DimD-\DimK+1} \to \CC$, we can express the operator $F(\Delta_\DimD,\dots,\Delta_\DimK)$ in the joint functional calculus of $\Delta_\DimD,\dots,\Delta_\DimK$ on $L^2(\sfera^\DimD)$ as
\begin{equation}\label{eq:jointcalculus}
F(\Delta_\DimD,\dots,\Delta_\DimK)
=
\sum_{(\ell_\DimD,\dots,\ell_\DimK) \in J_\DimD^{(\DimK)}}
F(\lambda^\DimD_{\ell_\DimD},\dots,\lambda^\DimK_{\ell_\DimK}) 
\, \pi^\DimD_{\ell_\DimD,\dots,\ell_\DimK},
\end{equation}
Correspondingly, the integral kernel $\Kern_{F(\Delta_\DimD,\dots,\Delta_\DimK)}$ of the operator $F(\Delta_\DimD,\dots,\Delta_\DimK)$ is given by
\begin{equation}\label{eq:kernelformula}
\Kern_{F(\Delta_\DimD,\dots,\Delta_\DimK)} =
\sum_{ (\ell_\DimD,\dots,\ell_\DimK) \in J_\DimD^{(\DimK)} } 
F(\lambda^\DimD_{\ell_\DimD},\dots,\lambda^\DimK_{\ell_\DimK})
\, K^\DimD_{\ell_\DimD,\dots,\ell_\DimK},
\end{equation}
and in particular,
by \eqref{eq:SW-dim} and \eqref{eq:Xlm_def_iniziale-norm}, for all $\spnt{\omega'}{\psi'} \in \sfera^\DimD$,
\begin{multline}\label{eq:kernel_l2norm}
\|\Kern_{F(\Delta_\DimD,\dots,\Delta_\DimK)}(\cdot, \spnt{\omega'}{\psi'})\|_{L^2(\sfera^\DimD)}^2 \\
=
\frac{1}{\meas_\DimK(\sfera^\DimK)}
\sum_{\substack{
(\ell_\DimD,\ldots,\ell_{\DimK}) \in J_\DimD^{(\DimK)} 
}}
\dimHarm_m(\sfera^\DimK)
|F(\lambda^\DimD_{\ell_\DimD},\dots,\lambda^\DimK_{\ell_\DimK})|^2
|X^\DimD_{\ell_\DimD,\ldots,\ell_\DimK}(\psi')|^2,
\end{multline}
where $\meas_\DimK$ is the Lebesgue measure on $\sfera^\DimK$, and $\dimHarm_m(\sfera^\DimK)$ denotes the dimension of $\Harm^m(\sfera^\DimK)$ as in \eqref{eq:dimensione}.

We note that the operators of the form \eqref{eq:jointcalculus} include those in the functional calculus of the Grushin operator $\opL_{\DimD,\DimK} = \Delta_\DimD - \Delta_\DimK$; namely,
\[
F(\opL_{\DimD,\DimK}) 
= 
\sum_{(\ell,m) \in I_\DimD^{(\DimK)}}
F(\lambda^{\DimD,\DimK}_{\ell,m}) 
\sum_{\substack{
(\ell_\DimD,\dots,\ell_\DimK) \in J_\DimD^{(\DimK)} \\
\ell_\DimD = \ell, \ \ell_\DimK = m
}}
\pi^\DimD_{\ell_\DimD,\dots,\ell_\DimK},
\]
where
\[\begin{split}
I_\DimD^{(\DimK)} &= \{ (\ell,m) \tc \ell \in \NN_\DimD, \, m \in \NN_{\DimK}, \ \ell\geq m + (\DimD - \DimK)/2 \} \\
&= \{ (\ell,m) \tc \exists (\ell_\DimD,\dots,\ell_\DimK) \in J_\DimD^{(\DimK)} \tc \ell_\DimD = \ell,\ \ell_\DimK = m \}
\end{split}\]
and, for all $(\ell,m) \in I_\DimD^{(\DimK)}$,
\begin{equation}\label{eq:grushin_eigenvalue}
\lambda^{\DimD,\DimK}_{\ell,m} = \lambda^\DimD_\ell - \lambda^\DimK_m.
\end{equation}

\subsection{Riesz-type bounds}\label{s:Riesz}
In this section we prove certain weighted $L^2$ bounds involving the joint functional calculus of $\Delta_\DimD,\dots,\Delta_\DimK$, which, in combination with the weighted spectral cluster estimates in Section \ref{s:weighted} below, play a fundamental role in satisfying the assumptions on the weight in the abstract theorem and proving our main result. A somewhat similar estimate was obtained in \cite[Lemma 2.5]{CaCiaMa} in the case $\DimD=2$ and $\DimK = 1$. Differently from \cite{CaCiaMa}, the estimate in Proposition \ref{prp:Riesz} below is proved for arbitrarily large powers of the weight; this prevents us from using the elementary ``quadratic form majorization'' method exploited in the previous paper, and requires a more careful analysis, based on the explicit eigenfunction expansion developed in the previous sections.

For later use, it is convenient to reparametrise the functions $X_{\ell,m}^\DimD$ defined in \eqref{eq:Xlm_def_iniziale}: namely, we introduce the functions $\tX_{\ell,m}^\DimD :[-1,1] \to \RR$ defined by
\begin{equation}\label{eq:Xtilde}
\tX_{\ell, m}^\DimD (x)
= c_{\ell m} (1-x^2)^{m/2-(\DimD-2)/4} P_{\ell-m-1/2}^{(m,m)}(x),
\end{equation}
where $(\ell,m) \in I_\DimD$ and $c_{\ell m}$ is given by \eqref{eq:clm}.

Let $t_{\DimD,\DimD} : \sfera^\DimD \to [-1,1]$ denote the restriction of the projection map $(z_0,\dots,z_\DimD) \mapsto z_\DimD$. Inductively, for $r=2,\dots,\DimD-1$, we define $t_{\DimD,r} : \sfera^\DimD \to [-1,1]$ by setting, for all $(\omega,\psi) \in \sfera^{\DimD-1} \times [-\pi/2,\pi/2]$,
\[
t_{\DimD,r}(\spnt{\omega}{\psi}) = \begin{cases}
t_{\DimD-1,r}(\omega) &\text{if } |\psi| < \pi/2,\\
0 &\text{otherwise}.
\end{cases}
\]
In particular, for all $(\omega,\psi) \in \sfera^\DimK \times (-\pi/2,\pi/2)^{\DimD-\DimK}$ and $r = \DimK+1,\dots,\DimD$,
\[
t_{\DimD,r}(\spnt{\omega}{\psi}) = \sin\psi_r,
\]
where $\psi=(\psi_{\DimK+1},\dots,\psi_\DimD)$.
Finally, we set, for $1 \leq \DimK < \DimD$,
\begin{equation}\label{eq:taudk}
\tau_{\DimD,\DimK} = \sum_{r=\DimK+1}^\DimD |t_{\DimD,r}|.
\end{equation}

\begin{proposition}\label{prp:Riesz}
Let $1 \leq \DimK < \DimD$. For all $N \in [0,\infty)$ and all $f \in L^2(\sfera^\DimD)$ such that $f \perp \ker\Delta_{\DimK+1}$,
\begin{equation}\label{eq:riesz_veryfinal}
\| \tau_{\DimD,\DimK}^N f \|_{L^2(\sfera^\DimD)} \lesssim_N \| (\opL_{\DimD,\DimK}/\Delta_{\DimK+1})^{N/2} f \|_{L^2(\sfera^\DimD)}.
\end{equation}
\end{proposition}
\begin{proof}
By interpolation, it is enough to prove the estimate in the case $N \in \NN$.

Let us first prove the inequality in the case $\DimK=\DimD-1$. From \eqref{eq:Xtilde} and known identities for Jacobi polynomials \cite[\S 10.9, eqs.\ (4) and (13), pp.\ 174--175]{EMOT}, one easily deduces that
\begin{equation}\label{eq:riesz_identity}
x \tX_{\ell, m}^\DimD (x) = \alpha_{\ell,m} \tX_{\ell+1,m}^\DimD(x) + \alpha_{\ell-1,m} \tX_{\ell-1,m}^\DimD(x)
\end{equation}
for all $(\ell,m) \in I_\DimD$, where
\[
\alpha_{\ell,m} = \sqrt{\frac{(\ell-m+1/2)(\ell+m+1/2)}{4\ell(\ell+1)}}.
\]
We remark that, in the case $(\ell-1,m) \notin I_\DimD$, the condition $(\ell,m) \in I_\DimD$ forces $\ell-m-1/2 = 0$ and $\alpha_{\ell-1,m} = 0$; in other words, the term with $\tX_{\ell-1,m}$ in the right-hand side of \eqref{eq:riesz_identity} appears only when $(\ell-1,m) \in I_\DimD$ too. On the other hand, if $(\ell,m) \in I_\DimD$, then $\alpha_{\ell,m} \simeq \sqrt{(\ell-m)(\ell+m)}/\ell$.
Consequently, by iterating \eqref{eq:riesz_identity}, we easily obtain, for all $N \in \NN$ and $(\ell,m) \in I_\DimD$,
\begin{equation}\label{eq:riesz_iterata}
x^N \tX_{\ell, m}^\DimD (x) = \sum_{j=0}^N \alpha^{N,j}_{\ell,m} \tX_{\ell-N+2j,m}^\DimD(x),
\end{equation}
where
\[
\alpha^{N,j}_{\ell,m} \simeq_N \begin{cases}
((\ell-m)(\ell+m)/\ell^2)^{N/2} &\text{if } (\ell-N+2j,m) \in I_\DimD,\\
0 &\text{otherwise.}
\end{cases}
\]

Let now $f \in L^2(\sfera^\DimD)$; then we can write  (see \eqref{eq:orth_decomp})
\[
f = \sum_{(\ell,m) \in I_\DimD} \sum_{Z \in E^{\DimD-1}_m} a_{\ell,m,Z} X_{\ell,m}^\DimD \otimes Z, 
\]
where, for all $m \in \NN_{\DimD-1}$, $E^{\DimD-1}_m$ is an orthonormal system of eigenfunctions of $\Delta_{\DimD-1}$ on $\sfera^{\DimD-1}$ of eigenvalue $\lambda^{\DimD-1}_m$.
Then from \eqref{eq:riesz_iterata} we deduce
\[
t_{\DimD,\DimD}^N f = \sum_{j=0}^N \sum_{(\ell,m) \in I_\DimD} \sum_{Z \in E_{\DimD-1,m}} a_{\ell,m,Z} \, \alpha^{N,j}_{\ell,m} X_{\ell-N+2j,m}^\DimD \otimes Z
\]
and consequently, by the orthogonality properties of the $X_{\ell,m}\otimes Z$ (see Section \ref{s:eigenfunctions}),
\begin{equation}\label{eq:preliminar_riesz}
\| t_{\DimD,\DimD}^N f \|_{L^2(\sfera^\DimD)}^2 \lesssim_N \sum_{(\ell,m) \in I_\DimD} \left(\frac{(\ell-m)(\ell+m)}{\ell^2}\right)^{N}  \sum_{Z \in E_{\DimD-1,m}} a_{\ell,m,Z}^2.
\end{equation}
Recall that $\Delta_\DimD (X^\DimD_{\ell,m} \otimes Z) = \lambda^\DimD_\ell (X^\DimD_{\ell,m} \otimes Z)$ and $\Delta_{\DimD-1} (X^\DimD_{\ell,m} \otimes Z) = \lambda^{\DimD-1}_m (X^\DimD_{\ell,m} \otimes Z)$, where $\lambda_\ell^\DimD = \ell^2-((\DimD-1)/2)^2$ and $\lambda_m^{\DimD-1} = m^2-((\DimD-2)/2)^2$; hence
\begin{equation}\label{eq:spectral_ineq}
\lambda_\ell^\DimD \geq \lambda_m^{\DimD-1} \text{ for all } (\ell,m) \in I_\DimD
\end{equation}
and
\[
\lambda_\ell^\DimD = \lambda_m^{\DimD-1} \text{ if and only if } \lambda_\ell^\DimD = 0.
\]
In particular, for all $(\ell,m) \in I_\DimD$,
\begin{equation}\label{eq:riesz_development}
\lambda_\ell^\DimD \simeq \ell^2, \qquad \lambda_\ell^\DimD - \lambda_m^{\DimD-1} \simeq \ell^2-m^2 \qquad\text{whenever }\lambda_\ell^\DimD \neq 0.
\end{equation}
If $f \perp \ker \Delta_\DimD$, then the coefficients $a_{\ell,m,Z}$ in \eqref{eq:preliminar_riesz} vanish unless $\lambda^\DimD_\ell \neq 0$, and from \eqref{eq:riesz_development} we deduce
\begin{equation}\label{eq:riesz_codim1}
\| t_{\DimD,\DimD}^N f \|_{L^2(\sfera^\DimD)} \lesssim_N \| ((\Delta_\DimD-\Delta_{\DimD-1})/\Delta_\DimD)^{N/2} f \|_{L^2(\sfera^\DimD)},
\end{equation}
which is \eqref{eq:riesz_veryfinal} in the case $\DimK=\DimD-1$.

Let now $2 \leq r \leq \DimD$. By the discussion in Section \ref{s:prelim}, up to null sets we can identify $\sfera^\DimD$ with $\sfera^r \times [-\pi/2,\pi/2]^{\DimD-r}$ with coordinates $(\omega,\psi)$ and measure $\cos\psi_\DimD^{\DimD-1} \cdots \cos\psi_{r+1}^r \,d\psi \,d\omega$. Consequently the space $L^2(\sfera^\DimD)$ is the Hilbert tensor product of the spaces $L^2(\sfera^r)$ and $L^2([-\pi/2,\pi/2]^{\DimD-r}, \cos\psi_\DimD^{\DimD-1} \cdots \cos\psi_{r+1}^r \,d\psi)$. Hence the inequality \eqref{eq:riesz_codim1}, applied with $\DimD=r$, yields a corresponding inequality on the sphere $\sfera^\DimD$, namely
\begin{equation}\label{eq:riesz_codim1_bis}
\| t_{\DimD,r}^N f \|_{L^2(\sfera^\DimD)} \lesssim_N \| ((\Delta_r-\Delta_{r-1})/\Delta_r)^{N/2} f \|_{L^2(\sfera^\DimD)}
\end{equation}
for all $f \perp \ker \Delta_r$.
On the other hand, from \eqref{eq:spectral_ineq} we deduce that $\Delta_{r_1} \geq \Delta_{r_2}$ spectrally whenever $1 \leq r_1 \leq r_2 \leq \DimD$ (recall that the $\Delta_r$ have a joint spectral decomposition, see Section \ref{s:eigenfunctions}), so \eqref{eq:riesz_codim1_bis} implies
\begin{equation}\label{eq:riesz_codim1_tris}
\| t_{\DimD,r}^N f \|_{L^2(\sfera^\DimD)} \lesssim_N \| ((\Delta_\DimD-\Delta_\DimK)/\Delta_{\DimK+1})^{N/2} f \|_{L^2(\sfera^\DimD)}
\end{equation}
whenever $\DimK < r \leq \DimD$ and $f \perp \ker \Delta_{\DimK+1}$. The desired inequality \eqref{eq:riesz_veryfinal} then follows by summing the inequalities \eqref{eq:riesz_codim1_tris} for $r=\DimK+1,\dots,\DimD$.
\end{proof}

\subsection{Estimates for ultraspherical polynomials}\label{s:bounds}

In this section we present a number of estimates for the functions $X_{\ell,m}^\DimD$ (or rather, their reparametrisations $\tX_{\ell,m}^\DimD$ from \eqref{eq:Xtilde}), which will play a crucial role in the subsequent developments. 

We first state some basic uniform bounds that follow from the previous discussion (see \eqref{eq:SW-dim} and \eqref{eq:decomp-Hl}).
In the statement below, we convene that $0^0 = 1$.

\begin{proposition}\label{prp:stime-classiche}
Let $\DimD \in \NN$, $\DimD \geq 2$.
\begin{enumerate}[label=(\roman*)]
\item\label{en:stime_somma}
For all $\ell \in \NN_\DimD$ and $x \in [-1,1]$,
\[
\sum_{\substack{
m \in \NN_{\DimD-1} \\ \, m \leq \ell
}} m^{\DimD-2} |\tX^\DimD_{\ell,m}(x)|^2 \lesssim_\DimD \, \ell^{\DimD-1} .
\]
\item\label{en:stime_unif}
$\|\tX_{\ell,m}^\DimD\|_\infty \lesssim_\DimD \, \ell^{(\DimD-1)/2}/ m^{(\DimD-2)/2}$ for all $(\ell,m) \in I_\DimD$.
\end{enumerate}
\end{proposition}

More refined pointwise estimates can be derived from 
asymptotic approximations of ultraspherical polynomials in terms of Hermite polynomials and Bessel functions, obtained in works of Olver \cite{Olver} and Boyd and Dunster \cite{BoydDunster} in the regimes $m \geq \epsilon \ell$ and $m \leq \epsilon \ell$ respectively, where $\epsilon \in (0,1)$.
Here and subsequently, for all $\ell,m \in \RR$ with $\ell \neq 0$ and $0 \leq m \leq \ell$,
$a_{\ell,m}$ and $b_{\ell,m}$ will denote the numbers in $[0,1]$ defined by
\begin{equation}\label{eq:def-bldm}
b_{\ell,m} = \frac{m}{\ell}
\end{equation}
and
\begin{equation}\label{eq:def-aldm}
a_{\ell,m}^2 =1-b_{\ell,m}^2=\frac{(\ell-m)(\ell+m)}{\ell^2}.
\end{equation}
The points $\pm a_{\ell,m} \in [-1,1]$ play the role of ``transition points'' for the functions $\tX^\DimD_{\ell,m}$ in the estimates that follow.

\begin{theorem}\label{thm:est_olver_boyddunster}
Let $\DimD \in \NN$, $\DimD \geq 2$.
Let $\epsilon \in (0,1)$. There exists $c \in (0,1)$ such that, for all $(\ell,m) \in I_\DimD$, if $m \geq \epsilon \ell$ then
\begin{equation}\label{eq:est_olver}
|\tX_{\ell,m}^\DimD(x)| \lesssim_{\DimD,\epsilon} \begin{cases}
(\ell^{-1} + |x^2-a_{\ell,m}^2|)^{-1/4} &\text{for all $x\in[-1,1]$,}\\
|x|^{-1/2} \exp(-c\ell x^2) &\text{for $|x| \geq 2 a_{\ell,m}$,}
\end{cases}
\end{equation}
while, if $m \leq \epsilon \ell$, then
\begin{equation}\label{eq:est_boyddunster}
|\tX^\DimD_{\ell,m}(x)|
\lesssim_{\DimD,\epsilon} \begin{cases}
y^{-(\DimD-2)/2} \left( \frac{(1+m)^{4/3}}{\ell^{2}} + |y^2-b_{\ell,m}^2|\right)^{-1/4} & \text{for all $x \in[-1,1]$,}\\
\ell^{(\DimD-1)/2} \, \, 2^{-m} &\text{if $y \leq b_{\ell,m}/(2e)$,}
\end{cases}
\end{equation}
where $y = \sqrt{1-x^2}$.
\end{theorem}

In the case $d=2$, the derivation of the estimates in Theorem \ref{thm:est_olver_boyddunster} from the asymptotic approximations in \cite{Olver,BoydDunster} is presented in \cite[Section 3]{CaCiaMa}; a number of variations and new ideas are required in the general case $d\geq 2$, and we refer to \cite{CCM-Jacobi} for a complete proof (indeed, in \cite{CCM-Jacobi} a stronger decay is proved in the regime $m \geq \epsilon \ell$ for $|x| \geq 2 a_{\ell,m}$ than the one given in \eqref{eq:est_olver}). Here we only remark that combining the above estimates yields the following bound.

\begin{corollary}\label{cor:est_combinata}
Let $\DimD \in \NN$, $\DimD \geq 2$. There exists $c \in (0,\infty)$ such that, for all $(\ell,m) \in I_\DimD$ and $x \in [-1,1]$,
\begin{equation}\label{eq:est_combinata}
|\tX^\DimD_{\ell,m}(x)|
\lesssim_\DimD \begin{cases}
y^{-(\DimD-2)/2} \left( \frac{1+m}{\ell^{2}} + |y^2-b_{\ell,m}^2|\right)^{-1/4} & \text{for all $x \in[-1,1]$,}\\
\ell^{(\DimD-1)/2} \,  \exp(-cm) &\text{if $y \leq b_{\ell,m}/(2e)$,}
\end{cases}
\end{equation}
where $y = \sqrt{1-x^2}$.
\end{corollary}
\begin{proof}
Let $\epsilon \in (0,1)$ be a parameter to be fixed later. If $m \leq \epsilon \ell$, the desired estimates immediately follow from \eqref{eq:est_boyddunster}, by taking any $c \leq \log 2$ (indeed, note that $(1+m)^{4/3} \geq 1+m$).

On the other hand, for $m \geq \epsilon \ell$, we may apply the estimates \eqref{eq:est_olver}. Note that $m \simeq \ell \gtrsim 1$ in this range, so $1/\ell \simeq (1+m)/\ell$; moreover $|x^2-a_{\ell,m}|^2 = |y^2-b_{\ell,m}|^2$ and $y \leq 1$, so the first estimate in \eqref{eq:est_combinata} immediately follows from the first estimate in \eqref{eq:est_olver}.

Assume now that $y \leq b_{\ell,m}/(2e)$. Since $b_{\ell,m} \geq \epsilon$ in this range, $a_{\ell,m}^2/(1-\epsilon^2) \leq 1$. Consequently
\[
x^2 = 1-y^2 \geq 1-\frac{b_{\ell,m}^2}{4e^2} = \frac{(4e^2-1) + a_{\ell,m}^2}{4e^2} \geq \min\left\{ \frac{4e^2-1}{4e^2}, \frac{1-(\epsilon/(2e))^2}{1-\epsilon^2} a_{\ell,m}^2 \right\}.
\]
This shows that, on the one side, $|x| \gtrsim 1$; on the other side, if $\epsilon \in (0,1)$ is chosen sufficiently large, then $|x| \geq 2 a_{\ell,m}$. Therefore we can apply the second estimate in \eqref{eq:est_olver} and obtain that
\[
|\tX_{\ell,m}^\DimD(x)| \lesssim \exp(-c' \ell)
\]
for a suitable constant $c' \in (0,\infty)$. Since $\ell \simeq m$ in this range, this clearly implies the second estimate in \eqref{eq:est_combinata} for an appropriate choice of $c$.
\end{proof}

\section{Weighted spectral cluster estimates}\label{s:weighted}

As a consequence of the estimates in Section \ref{s:bounds}, we obtain ``weighted spectral cluster estimates'' for the Grushin operators $\opL_{\DimD,\DimK}$.

Fix $\DimD,\DimK \in \NN$ with $1 \leq \DimK < \DimD$. 
For $(\ell,m)\in I_\DimD^{(\DimK)}$ and $x = (x_\DimD,x_{\DimD-1},\dots,x_{\DimK+1}) \in [-1,1]^{\DimD-\DimK}$, define
\begin{equation}\label{eq:def_tSpHddd-multiple}
\tSpHddd_{\ell,m}^{\DimD,\DimK}(x) =
\sum_{\substack{
(\ell_\DimD,\dots,\ell_\DimK) \in J_\DimD^{(\DimK)} \\
\ell_\DimD = \ell, \ \ell_\DimK = m
}}
\bigl|\tX_{\ell_\DimD,\ell_{\DimD-1}}^\DimD(x_\DimD) \bigr|^2
\ldots
\bigl|\tX_{\ell_{\DimK+1},\ell_{\DimK}}^{\DimK+1}(x_{\DimK+1}) \bigr|^2,
\end{equation}
where $\tX_{r,s}^\DimD$ has been defined in \eqref{eq:Xtilde}. We are interested in bounds for suitable weighted sums of the $\tSpHddd_{\ell,m}^{\DimD,\DimK}$ for indices $\ell,m$ such that the eigenvalue $\sqrt{\lambda_{\ell,m}^{\DimD,\DimK}}$ of $\sqrt{\opL_{\DimD,\DimK}}$ ranges in an interval of unit length (whence the name ``spectral cluster''). The bounds that we obtain are different in nature according to whether $m \leq \epsilon \ell$ or $m \geq \epsilon \ell$ for some fixed $\epsilon \in (0,1)$, and are presented as separate statements. We remark that, in the case $m \leq \epsilon \ell$, the eigenvalue $\lambda_{\ell,m}^{\DimD,\DimK}$ of $\opL_{\DimD,\DimK}$ is comparable with the eigenvalue $\lambda^\DimD_\ell$ of $\Delta_\DimD$; consequently, the range $m \leq \epsilon \ell$ will be referred to as the ``elliptic regime'', while the range $m \geq \epsilon \ell$ will be called the ``subelliptic regime''.

\begin{proposition}[subelliptic regime]\label{prp:indici_bassi_primo_piufattori}
Let $\epsilon \in (0,1)$ and $\DimD \geq 2$. Fix $1 \leq \DimK \leq \DimD-1$.
Then, for all $i \in \NN \setminus \{0\}$, $\alpha \in [0,\DimK/2)$, and $x\in[-1,1]^{\DimD-\DimK}$,
\begin{equation}\label{eq:indici_bassi_primo_peso_secondo_strato_KK}
 \sum_{\substack{
(\ell,m)\in I_\DimD^{(\DimK)} \\
m \geq \epsilon \ell \\
\lambda_{\ell,m}^{\DimD,\DimK} \in [i^2,(i+1)^2]
}}
\dimHarm_m(\sfera^{\DimK}) \ 
\tSpHddd_{\ell,m}^{\DimD,\DimK}(x)
\Bigl[ \sqrt{\lambda_{\ell,m}^{\DimD,\DimK}}/\ell \Bigr]^{2\alpha} 
\lesssim_\epsilon \, i^{\DimD-1} \ \min \{i, 1/|x|\}^{\DimK-2\alpha}  .
\end{equation}
\end{proposition}

\begin{proposition}[elliptic regime]\label{prp:indici_bassi_primo_peso_secondo_stratopiu}
Let $\epsilon \in (0,1)$ and $\DimD \geq 2$. Fix $1 \leq \DimK \leq \DimD-1$.
Then, for all $i \in \NN \setminus \{0\}$ and $x \in [-1,1]^{\DimD-\DimK}$,
\begin{equation}\label{eq:indici_bassi_primo_peso_secondo_strato_dk}
\sum_{\substack{
(\ell,m)\in I_\DimD^{(\DimK)} \\  
m \leq \epsilon \ell \\
\lambda_{\ell,m}^{\DimD,\DimK} \in[i^2, (i+1)^2]
}}
\dimHarm_m(\sfera^\DimK) \
\tSpHddd_{\ell,m}^{\DimD,\DimK}(x) 
\lesssim_\epsilon \, i^{\DimD-1} ,
\end{equation}
where $\tSpHddd_{\ell,m}^{\DimD,\DimK}$ was defined in \eqref{eq:def_tSpHddd-multiple}.
\end{proposition}

Analogous estimates are proved in \cite[Section 4]{CaCiaMa} in the case $\DimD=2$ and $\DimK=1$; in that case, each of the products in \eqref{eq:def_tSpHddd-multiple} reduces to a single factor. Treating the general case, with multiple factors, presents substantial additional difficulties, as one may appreciate from the discussion below.

The following lemma will be repeatedly used in the proofs that follow, in combination with \cite[Lemma 4.1]{CaCiaMa} and \cite[Lemma 5.7]{DM}, to justify the passage from a sum to the corresponding integral.

\begin{lemma}\label{lem:der_Xi}
For $a,t \in \RR$, $s \in (0,\infty)$, define
\begin{equation}\label{eq:def_Xi}
\Xi(a,s,t) = (s+|a-t|)^{-1/2}.
\end{equation}
Let $\kappa \in [1,\infty)$. Let $\Omega \subseteq \RR^n$, and $\alpha_j : \Omega \to (0,\infty)$, $\beta_j: \Omega \to \RR$ be such that
\[
|\nabla \alpha_j|, |\nabla \beta_j| \leq \kappa \alpha_j
\]
for $j=1,\dots,N$.
Define $\tilde \Xi(y,x) = \prod_{j=1}^N \Xi(y_j,\alpha_j(x),\beta_j(x))$ for $y \in \RR^n$ and $x \in \Omega$. Then, for all $y \in \RR^N$ and $x \in \Omega$,
\[
|\nabla_x \tilde \Xi(y,x)| \leq N \kappa \, \tilde\Xi(y,x).
\]
\end{lemma}
\begin{proof}
By the Leibniz rule, it is enough to consider the case $N=1$. Set $\alpha = \alpha_1$, $\beta = \beta_1$. Define $X(a,s,t) = s+|a-t|$ and $\tilde X(y,x) = X(y,\alpha(x),\beta(x))$. 
Note now that
\[
X(a,s,t) \geq s, \qquad |\partial_s X(a,s,t)|,|\partial_t X(a,s,t)| \leq 1,
\]
whence, by the chain rule,
\[
|\nabla_x \tilde X(y,x)| \leq |\nabla_x \alpha(x)| + |\nabla_x \beta(x)| \leq 2\kappa \alpha(x) \leq 2\kappa \tilde X(y,x)
\]
and
\[
\frac{|\nabla_x \tilde\Xi(y,x)|}{\tilde \Xi(y,x)} = \frac{1}{2} \frac{|\nabla_x \tilde X(y,x)|}{\tilde X(y,x)} \leq \kappa,
\]
as desired.
\end{proof}

\subsection{The subelliptic regime}

Here we prove Proposition \ref{prp:indici_bassi_primo_piufattori}. To this aim, we first present a couple of lemmas that will allow us to perform a particularly useful change of variables in the proof.

\begin{lemma}\label{lem:determinant}
Let $w \in \RR^n$ and define the matrix $M(w) = (m_{i,j}(w))_{i,j=1}^n$ by
\[
m_{i,j}(w) = \begin{cases}
1 &\text{if } i=j,\\
w_i &\text{if } i>j,\\
-w_i &\text{if } i<j.
\end{cases}
\]
Then
\[
\det M(w) = \sum_{\substack{S \subseteq \{1,\dots,n\} \\ |S| \text{ even}}} \prod_{j\in S} w_j.
\]
\end{lemma}
\begin{proof}
Observe that $m_{i,j}(w) = \delta_{i,j} + \rho_{i,j} w_j$, where
\[
\rho_{j,s} = \begin{cases}
1 &\text{if } s<j,\\
-1 &\text{if } s>j,\\
0 &\text{if } s=j.
\end{cases}
\]
Consequently, if $\SPerm_n$ denotes the group of permutations of the set $\{1,\dots,n\}$ and $\epsilon(\sigma)$ denotes the signature of the permutation $\sigma$, then
\[\begin{split}
\det M(w)
&= \sum_{\sigma \in \SPerm_n} \epsilon(\sigma) \prod_{j=1}^n m_{i,j}(w)  \\
&= \sum_{\sigma \in \SPerm_n} \epsilon(\sigma) \prod_{j \tc \sigma(j) \neq j} \rho_{j,\sigma(j)} w_j \\
&= \sum_{S \subseteq \{1,\dots,n\}}
\left(\prod_{j \in S} w_j \right)
\sum_{\substack{\sigma \in \SPerm_n \\ \sigma|_{S^\compl} = \id }}
\epsilon(\sigma) \prod_{j \in S} \rho_{j,\sigma(j)} \\
&= \sum_{S \subseteq \{1,\dots,n\}}
\left(\prod_{j \in S} w_j \right)
\det(\rho_{l,m})_{l,m=1}^{|S|},
\end{split}\]
where $S^\compl = \{1,\dots,n\} \setminus S$. We note that $(\rho_{l,m})_{l,m=1}^{|S|}$ is a skewsymmetric matrix, so its determinant vanishes when $|S|$ is odd; if $|S|$ is even, instead, its determinant is the square of its pfaffian, and using the Laplace-type expansion for pfaffians (see, e.g., \cite[\S III.5, p.\ 142]{Artin}) one can see inductively that the determinant is $1$.
\end{proof}

\begin{lemma}\label{lem:jacobiandeterminant}
Let $\Omega = \{ v \in \RR^n \tc \hat v_j \neq -1 \text{ for all } j =1,\dots,n \}$, where
\[
\hat v_j = \sum_{r=j+1}^n v_r - \sum_{r=1}^{j-1} v_r. 
\]
Let $v \mapsto w$ be the map from $\Omega$ to $\RR^n$
defined by
\[
w_j = \frac{v_j}{1+\hat v_j}
\]
for $j=1,\dots,n$.
Then
\[
\det (\partial_{v_s} w_j)_{j,s=1}^n =
\left( \prod_{j =1}^n \frac{1}{1+\hat v_j} \right)
\sum_{\substack{
S \subseteq \{1,\dots,n\} \\
|S| \text{ even}
}}
\prod_{j \in S} \frac{v_j}{1+\hat v_j}.
\]
Moreover, for all $\epsilon \in (0,1)$, the map $v \mapsto w$ is injective when restricted to
\[
\Omega_\epsilon \defeq \left\{ v \in \RR^n \tc v_j \geq 0 \ \forall j=1,\dots,n, \, \sum_j v_j \leq \epsilon \right\}.
\]
\end{lemma}
\begin{proof}
From the definition it is immediate that
\[
\partial_{v_s} w_j = \frac{1}{1+\hat v_j} m_{j,s}(w),
\]
where $M(w) = \{ m_{j,s}(w) \}_{j,s=1}^n$ is the matrix defined in Lemma \ref{lem:determinant},
so
\[
\det (\partial_{v_s} w_j)_{j,s=1}^n = \left(\prod_{j=1}^n \frac{1}{1+\hat v_j} \right) \det M(w),
\]
and the desired expression for the determinant follows from Lemma \ref{lem:determinant}.

Note that, if $v \in \Omega_\epsilon$, $0 \leq v_j, |\hat v_j| \leq \sum_j v_j \leq \epsilon < 1$, so $1+\hat v_j > 0$ and 
$\Omega_\epsilon \subseteq \Omega$. In addition, the equations $w_j = v_j/(1+\hat v_j)$ are equivalent to $v_j - w_j \hat v_j = w_j$, that is,
\[
M(w) v = w.
\]
Since $w_j = v_j/(1+\hat v_j) \geq 0$, from Lemma \ref{lem:determinant} it follows that $\det M(w) \geq 1$, so the matrix $M(w)$ is invertible and the above equation is equivalent to $v = M(w)^{-1} w$; in other words, if $v \in \Omega_\epsilon$, then $v$ is uniquely determined by its image $w$ via the map $v \mapsto w$, that is, the map restricted to $\Omega_\epsilon$ is injective.
\end{proof}

\begin{proof}[Proof of Proposition \ref{prp:indici_bassi_primo_piufattori}]
We start by observing that, for all $(\ell,\DimK) \in I_\DimD^{(\DimK)}$, if we assume $\epsilon \ell \leq m$,
then, for all $(\ell_\DimD,\dots,\ell_\DimK) \in J_\DimD^{(\DimK)}$ with $\ell_\DimD = \ell$ and $\ell_\DimK = m$,
\begin{equation}\label{eq:inequalities_epsilon_ell}
\epsilon \ell_{j+1} \leq \ell_j, \qquad j \in \{\DimK,\dots,\DimD-1\},
\end{equation}
and in particular
\begin{equation}\label{eq:equiv_ellj}
\ell_j \simeq \ell \gtrsim 1,
\quad
\text{for all } j \in \{\DimK,\dots,\DimD\}.
\end{equation}
We also note that
\[
\lambda_{\ell,m}^{\DimD,\DimK} + (\DimD+\DimK-2)(\DimD-\DimK)/4 = \ell^2 - m^2.
\]
Thus \eqref{eq:indici_bassi_primo_peso_secondo_strato_KK} will follow from 
\begin{multline}\label{eq:indici_alti_mod_new_2fattori-I-vv}
\sum_{\substack{
(\ell_\DimD,\dots,\ell_\DimK) \in J_\DimD^{(\DimK)} \\
\epsilon \ell_\DimD \leq \ell_\DimK \\
\ell_\DimD^2 - \ell_\DimK^2 \in [i^2,(i+1)^2]
}}
\ell_\DimD^{\DimK-1-2\alpha}
\bigl|\tX_{\ell_\DimD,\ell_{\DimD-1}}^\DimD(x_\DimD) \bigr|^2
\dots
\bigl|\tX_{\ell_{\DimK+1},\ell_{\DimK}}^{\DimK+1} (x_{\DimK+1}) \bigr|^2
\\
\lesssim_\epsilon \ i^{\DimD-1-2\alpha} \min\{i,|x|^{-1}\}^{\DimK-2\alpha} ,
\end{multline}
since it will suffice to apply \eqref{eq:indici_alti_mod_new_2fattori-I-vv} $h \defeq \lceil (\DimD+\DimK-2)(\DimD-\DimK)/4 \rceil$ times, with $i$ replaced by $i,i+1,\ldots,i+h-1$, respectively. Due to \eqref{eq:symmetry}, we may restrict without loss of generality to $x \in [0,1]^{\DimD-\DimK}$. In addition, for each fixed $i$, the sum in the left-hand side of \eqref{eq:indici_alti_mod_new_2fattori-I-vv} is finite, since $\ell_\DimD - \ell_\DimK \gtrsim 1$ and therefore
\[
\ell_\DimD \leq \ell_\DimD + \ell_\DimK \lesssim \ell_\DimD^2 - \ell_\DimK^2 \leq (i+1)^2;
\]
the boundedness of the functions $\tX^{\DimD-j+1}_{\ell_{\DimD-j+1},\ell_{\DimD-j}}$ (see Proposition \ref{prp:stime-classiche}\ref{en:stime_unif}) then shows that the estimate \eqref{eq:indici_alti_mod_new_2fattori-I-vv} is trivially true for each fixed $i$ (with a constant dependent on $i$), and therefore it is enough to prove it for $i$ sufficiently large.

It is convenient to reindex the above sum. Let us set
\[
p=\ell_\DimD+\ell_\DimK, \quad
q_j=\ell_{\DimD-j+1}-\ell_{\DimD-j},
\qquad \text{for all } j \in \{1,\dots,\DimD-\DimK\}.
\]
and let us write
\[
Q \defeq q_1 + \dots + q_{\DimD-\DimK}.
\]
Then the condition $(\ell_\DimD,\dots,\ell_\DimK) \in J_\DimD^{(\DimK)}$ implies
\[
q_1,\dots,q_{\DimD-\DimK} \in \NN+1/2, \qquad p \in \NN + (\DimD+\DimK-2)/2, \qquad p \geq Q + \DimK-1.
\]
Moreover
\[
\ell_\DimD-\ell_\DimK = Q, \qquad \ell_\DimD^2 - \ell_\DimK^2 = p Q
\]
and, for $j=1,\dots,\DimD-\DimK$,
\begin{equation}\label{eq:sum-of-two-consecutive-ind}
\ell_{\DimD-j+1}+\ell_{\DimD-j} = p + \hat q_j, \qquad\text{where } \hat q_j \defeq \sum_{r=j+1}^{\DimD-\DimK} q_r - \sum_{r=1}^{j-1} q_r;
\end{equation}
in particular
\begin{equation}\label{eq:def_adlsj}
a_{\ell_{\DimD-j+1},\ell_{\DimD-j}}^2 = 1-\frac{\ell_{\DimD-j}^2}{\ell_{\DimD-j+1}^2} = \frac{4 q_j (p+\hat q_j)}{(p+\hat q_j + q_j)^2} \simeq \frac{q_j}{p}.
\end{equation}
The condition $\epsilon \ell_\DimD \leq \ell_\DimK$ is then equivalent to
\[
Q \leq \bar\epsilon^4 p,
\]
where $\bar\epsilon = \left(\frac{1-\epsilon}{1+\epsilon}\right)^{1/4} \in (0,1)$, and implies, by \eqref{eq:inequalities_epsilon_ell},
\begin{equation}\label{eq:inequalities_epsilon_pqj}
q_j \leq \bar\epsilon^4 (p+\hat q_j)
\end{equation}
for $j=1,\dots,\DimD-\DimK$.
As previously discussed, it will be enough to prove the estimate for $i$ sufficiently large; in the following we will assume that
\[
1+1/i \leq \bar\epsilon^{-1}.
\]
Under this assumption on $i$, if $pQ \in [i^2,(i+1)^2]$, then
\[
Q \leq \bar\epsilon^2 \sqrt{pQ} \leq \bar\epsilon^2 (i+1) \leq \bar\epsilon i.
\]

\medskip

Let us first consider the range 
\begin{equation}\label{eq:ineq_xinfty-multiple}
|x|_\infty \geq 2 \max_{j \in \{1,\ldots,\DimD-\DimK\}} a_{\ell_{\DimD-j+1},\ell_{\DimD-j}},
\end{equation}
where 
\begin{equation}\label{eq:xinfty-multiple}
|x|_\infty = \max_{j \in \{1,\ldots,\DimD-\DimK\}} |x_{\DimD-j+1}| .
\end{equation}
In light of \eqref{eq:est_olver}, the inequalities
\begin{equation}\label{eq:est_olver-multiple}
|\tX_{\ell_{\DimD-j+1},\ell_{\DimD-j}}^{\DimD-j+1}(x_{\DimD-j+1}) |
\lesssim \ell_{\DimD-j+1}^{1/4}
\simeq p^{1/4}
\end{equation}
hold for all $j \in \{1,\ldots,\DimD-\DimK\}$.
Moreover, for one of the quantities 
\[
|\tX_{\ell_\DimD,\ell_{\DimD-1}}^\DimD(x_\DimD)|, \, \dots, \, |\tX_{\ell_{\DimK+1},\ell_{\DimK}}^{\DimK+1}(x_{\DimK+1})|
\]
the better bound $|x|^{-1/2} \exp(-cp |x|^2) $ holds for some $c>0$, thanks to the second estimate in \eqref{eq:est_olver} and to \eqref{eq:equiv_ellj}.
As a consequence, we obtain
\[\begin{split}
\bigl|\tX_{\ell_\DimD,\ell_{v-1}}^\DimD(x_\DimD) \bigr|^2
\cdots
\bigl|\tX_{\ell_{\DimK+1},\ell_{\DimK}}^{\DimK+1} (x_{\DimK+1}) \bigr|^2
&\lesssim p^{(\DimD-\DimK-1)/2}
|x|^{-1}
 \exp(-2cp |x|^2) 
\\
&\lesssim_N |x|^{-(\DimD-\DimK)}
(p |x|^2)^{-N} 
\end{split}\]
for arbitrarily large $N \in \NN$.
Note then that the conditions \eqref{eq:ineq_xinfty-multiple} and \eqref{eq:def_adlsj} imply
\[
|x|^2\gtrsim Q/p,
\]
which, together with $\ell_\DimD^2 - \ell_\DimK^2 = p Q \in [i^2, (i+1)^2]$, yields
\[
i|x| \gtrsim Q \gtrsim 1.
\]
Then
\[\begin{split}
&\sum_{\substack{
(\ell_\DimD,\dots,\ell_\DimK) \in J_\DimD^{(\DimK)} \\
\epsilon \ell_\DimD \leq \ell_\DimK \\
\ell_\DimD^2 - \ell_\DimK^2 \in [i^2,(i+1)^2] \\
|x|_\infty \geq 2 \max_j a_{\ell_{\DimD-j+1},\ell_{\DimD-j}}
}}
\ell_\DimD^{\DimK-1-2\alpha}
\bigl|\tX_{\ell_\DimD,\ell_{\DimD-1}}^\DimD(x_\DimD) \bigr|^2
\cdots
\bigl|\tX_{\ell_{\DimK+1},\ell_{\DimK}}^{\DimK+1}(x_{\DimK+1}) \bigr|^2
\\
&\lesssim 
|x|^{-(\DimD-\DimK)-2N}
\sum_{\substack{
Q \leq \bar\epsilon^4 p \\
pQ \in[i^2, (i+1)^2]
\\
|x|^2\gtrsim Q/p 
}}
p^{\DimK-1-2\alpha-N}
\\
&\lesssim 
|x|^{-(\DimD-\DimK)-2N}
\sum_{\substack{
Q \lesssim i|x|
}}
\sum_{
p \in[i^2/Q, (i+1)^2/Q]}
p^{\DimK-1-2\alpha-N}
\\
&\lesssim 
|x|^{-(\DimD-\DimK)-2N}
\sum_{\substack{
Q \lesssim i|x|
}}
(i/Q)
(i^2/Q)^{\DimK-1-2\alpha-N} 
\\
&=
i^{2\DimK-1-4\alpha-2N} 
|x|^{-(\DimD-\DimK)-2N} 
\sum_{\substack{
Q \lesssim i|x|
}}
Q^{N-\DimK+2\alpha} 
\\
&\lesssim 
i^{2\DimK-1-4\alpha-2N} 
|x|^{-(\DimD-\DimK)-2N} 
( i|x|)^{N+\DimD-2\DimK+2\alpha} 
\\
&=
i^{\DimD-1-2\alpha} 
|x|^{-\DimK+2\alpha} \
(i|x|)^{-N}
\lesssim
i^{\DimD-1-2\alpha} \min\{i,|x|^{-1}\}^{\DimK-2\alpha} ,
\end{split}\]
since $i |x| \gtrsim 1$ and $\DimK - 2\alpha > 0$, provided $N$ is large enough.
Note that, in estimating the sum in $p$, we used the fact that the interval $[i^2/Q,(i+1)^2/Q]$ has length $(2i+1)/Q \simeq i/Q \gtrsim 1$.

\medskip

Let us now discuss the range
\begin{equation}\label{eq:ineq_secondcase-multiple}
|x|_\infty \leq 2 \max_{j \in \{1,\ldots,\DimD-\DimK\}} a_{\ell_{\DimD-j+1},\ell_{\DimD-j}} .
\end{equation}
We first note that \eqref{eq:ineq_secondcase-multiple} implies
\[
|x|_\infty^2 \lesssim {Q}/{p},
\]
which, combined with $p Q \in [i^2, (i+1)^2]$ and $Q \in \NN + (\DimD-\DimK)/2$, implies
\[
p \lesssim {i/ |x|_\infty} \qquad\text{and}\qquad Q \gtrsim \max\{1, i |x|_\infty\}.
\]
Note that, by \eqref{eq:def_adlsj}, for all $j=1,\dots,\DimD-\DimK$,
\begin{equation}\label{eq:trans_phi}
(a_{\ell_{\DimD-j+1},\ell_{\DimD-j}})^2= \varphi(q_j/(p + \hat q_j)),
\end{equation}
where $\varphi(w) = 4w/(1+w)^2$. 
Note that the map $\varphi : [0,1] \to [0,1]$ is an increasing bijection, 
such that $w \leq \varphi(w) \leq 4w$; its derivative is given by
$\varphi'(w) = 4 \frac{1-w}{(1+w)^3}$
and vanishes only at $w=1$.
As a consequence, setting $\bar x_j = \sqrt{\varphi^{-1}(x_j^2)}$, with $j \in \{1,\ldots,\DimD-\DimK\}$, 
one has $\bar x_j \simeq |x_j|$; moreover, in light of \eqref{eq:trans_phi} and \eqref{eq:inequalities_epsilon_pqj},
\[
\bigl| x_{\DimD-j+1}^2 -(a_{\ell_{\DimD-j+1},\ell_{\DimD-j}})^2 \bigr|
\simeq
\bigl| \bar x_{\DimD-j+1}^2- q_j/(p + \hat q_j) \bigr|,
\]
uniformly for $x \in [0,1]^{\DimD-\DimK}$.
In particular, in this range, by \eqref{eq:est_olver},
\[
\bigl|\tX_{\ell_{\DimD-j+1},\ell_{\DimD-j}}^{\DimD-j+1}(x_{\DimD-j+1}) \bigr|^2 
\lesssim \Xi(\bar x_{\DimD-j+1}^2, 1/p, q_j/ ( p + \hat q_j))
\]
for all $j = 1,\dots,\DimD-\DimK$, where $\Xi$ is defined as in \eqref{eq:def_Xi}.
\\
Then
\[\begin{split}
&\sum_{\substack{
(\ell_\DimD,\dots,\ell_\DimK) \in J_\DimD^{(\DimK)} \\
\epsilon \ell_\DimD \leq \ell_\DimK \\
\ell_\DimD^2 - \ell_\DimK^2 \in [i^2,(i+1)^2] \\
|x|_\infty \leq 2 \max_j a_{\ell_{\DimD-j+1},\ell_{\DimD-j}}
}}
\ell_\DimD^{\DimK-1-2\alpha}
\bigl|\tX_{\ell_\DimD,\ell_{\DimD-1}}^\DimD (x_\DimD) \bigr|^2
\dots
\bigl|\tX_{\ell_{\DimK+1},\ell_{\DimK}}^{\DimK+1} (x_{\DimK+1}) \bigr|^2
\\
&\lesssim 
\sum_{\substack{
Q \leq \bar\epsilon^4 p \\
pQ \in [i^2,(i+1)^2] \\
|x|^2 \lesssim Q/p
}}
p^{\DimK-1-2\alpha}
\prod_{j=1}^{\DimD-\DimK} \Xi(\bar x_{\DimD-j+1}^2, 1/p,  q_j/ ( p + \hat q_j)) 
\\
&
\lesssim 
\sum_{\substack{
\max\{1,i|x|\} \lesssim Q \leq \bar\epsilon i
}}
\left(\frac{i^2}{Q}\right)^{\DimK-1-2\alpha}
\sum_{p \in [i^2/Q,(i+1)^2/Q]}
\vec\Xi(\bar x,p,q) ,
\end{split}\]
where $\bar x = (\bar x_\DimD,\dots,\bar x_{\DimK+1})$, $q = (q_1,\dots,q_{\DimD-\DimK})$ and
\[
\vec\Xi(\bar x,p,q) = 
\prod_{j=1}^{\DimD-\DimK} \Xi(\bar x_{\DimD-j+1}^2, 1/p,  q_j/(p + \hat q_j)) .
\]
It is easily seen that
\[
|\partial_p (1/p)|, |\partial_p (q_j/ (p + \hat q_j))| \lesssim 1/p
\]
for all $j=1,\dots,\DimD-\DimK$, on the range of summation (note that $q_j + |\hat q_j| \leq Q \leq \bar\epsilon^2 i^2/Q \leq \bar\epsilon^2 p$ and $\bar\epsilon < 1$, whence $p + \hat q_j\simeq p \gtrsim q_j \gtrsim 1$). Thus, by Lemma \ref{lem:der_Xi},
\[
|\partial_p \vec\Xi(\bar x,p,q)| \lesssim \vec\Xi(\bar x,p,q).
\]
Moreover the interval $[i^2/Q,(i+1)^2/Q]$ has length $(2i+1)/Q \simeq i/Q \gtrsim 1$.
Hence, by \cite[Lemma 4.1]{CaCiaMa},
\[\begin{split}
&\sum_{\substack{
\max\{1,i|x|\} \lesssim Q \leq \bar\epsilon i
}}
\left(\frac{i^2}{Q}\right)^{\DimK-1-2\alpha}
\sum_{p \in [i^2/Q,(i+1)^2/Q]}
\vec\Xi(\bar x,p,q) 
\\
&\lesssim 
\sum_{\substack{
\max\{1,i|x|\} \lesssim Q \leq \bar\epsilon i
}}
\left(\frac{i^2}{Q}\right)^{\DimK-1-2\alpha}
\int_{i^2/Q}^{(i+1)^2/Q}
\vec\Xi(\bar x,p,q) 
\,dp \\
&\simeq
i^{2\DimK-1-4\alpha}
\int_{i}^{i+1}
\sum_{\substack{
\max\{1,i|x|\} \lesssim Q \leq \bar\epsilon i
}}
\hat\Xi(\bar x,u,q)
 \,du,
\end{split}\]
where the change of variables $p = u^2/Q$ was used, and
\[\begin{split}
\hat\Xi(\bar x,u,q) &= Q^{2\alpha-\DimK} \vec\Xi(\bar x,u^2/Q,q) \\
&= Q^{2\alpha-\DimK} \prod_{j=1}^{\DimD-\DimK} \Xi(\bar x_{\DimD-j+1}^2, Q/u^2, q_j Q/(u^2+\hat q_j Q)).
\end{split}\]
It is easily checked that
\[
|\nabla_q (Q/u^2)|, |\nabla_q (q_j Q/(u^2 + \hat q_j Q))|\lesssim Q/u^2
\]
for all $j=1,\dots,\DimD-\DimK$,
on the range of summation (note that $|\hat q_j| Q \leq Q^2 
\leq \bar\epsilon^2 i^2 \leq \bar\epsilon^2 u^2$ 
and $\bar\epsilon < 1$, so $u^2 + \hat q_j Q \simeq u^2$), and therefore, by Lemma \ref{lem:der_Xi} and the Leibniz rule,
\[
|
\nabla_q \hat\Xi (\bar x,u,q)
| \lesssim \hat\Xi(\bar x,u,q).
\]
Hence, by \cite[Lemma 5.7]{DM},
\[\begin{split}
&i^{2\DimK-1-4\alpha}
\int_{i}^{i+1}
\sum_{\substack{
\max\{1,i|x|\} \lesssim Q \leq \bar\epsilon i
}}
\hat\Xi(\bar x,u,q)
 \,du
\\
&\lesssim
i^{2\DimK-1-4\alpha}
\int_i^{i+1}
\int_{\substack{
\max\{1,i|x|\} \lesssim Q \leq \bar\epsilon i
}}
\hat\Xi(\bar x,u,q)
\,dq \,du 
\\
&\simeq
\iint_{\substack{
\max\{1,i|x|\} \lesssim Q \leq \bar\epsilon i \\
pQ \in [i^2,(i+1)^2]
}}
\vec\Xi(\bar x,p,q) 
\, p^{\DimK-1-2\alpha}
\,dq \,dp \\
&\lesssim
\iint_{\substack{
\max\{i^{-1},|x|\}^2 \lesssim V \leq \bar\epsilon^2 \\
p^2 V \in [i^2,(i+1)^2] 
}}
\vec\Xi(\bar x,p,pv) \,p^{\DimD-1-2\alpha}
\,dp \,dv 
\end{split}\]
where the change of variables $q_j = p v_j$ was used, and $V = \sum_{j=1}^{\DimD-\DimK} v_j$
(note that $Q \leq \bar\epsilon i$ and $pQ \geq i^2$ implies $V = Q^2/(pQ) \leq \bar\epsilon^2$). Now,
\[\begin{split}
&\iint_{\substack{
\max\{i^{-1},|x|\}^2 \lesssim V \leq \bar\epsilon^2 \\
p^2 V \in [i^2,(i+1)^2] 
}}
\vec\Xi(\bar x,p,pv) \,p^{\DimD-1-2\alpha}
\,dp \,dv \\
&
\lesssim
\iint_{\substack{
\max\{i^{-1},|x|\}^2 \lesssim V \leq \bar\epsilon^2 \\
p^2 V \in [i^2,(i+1)^2] 
}}
\left(\frac{i}{\sqrt{V}}\right)^{\DimD-1-2\alpha} 
\prod_{j=1}^{\DimD-\DimK} \left|\bar x_{\DimD-j+1}^2 - \frac{v_j}{1+\hat v_j}\right|^{-1/2} 
\,dp \,dv
\\
&\lesssim
i^{\DimD-1-2\alpha}
\int_{\substack{
\max\{i^{-1},|x|\}^2 \lesssim V \leq \bar\epsilon^2 
}}
V^{-(\DimD-2\alpha)/2} 
\prod_{j=1}^{\DimD-\DimK} \left|\bar x_{\DimD-j+1}^2 - \frac{v_j}{1+\hat v_j}\right|^{-1/2} 
\,dv ,
\end{split}\]
where
$\hat v_j = \sum_{r=j+1}^{\DimD-\DimK} v_r - \sum_{r=1}^{j-1} v_j$, and the fact that the interval $[i/\sqrt{V},(i+1)/\sqrt{V}]$ has length $V^{-1/2}$ was used.
We can now use the change of variables
\[
w_j = \frac{v_j}{1+\hat v_j},
\]
observing that $w_j \simeq v_j$ for all $j \in \{1,\ldots,\DimD-\DimK\}$, and (see Lemma \ref{lem:jacobiandeterminant} below)
\begin{equation}\label{eq:jacobianvw}
\det (\partial_{v_s} w_j)_{j,s=1,\dots,\DimD-\DimK}
=
\left(\prod_{j=1}^{\DimD-\DimK}  \frac{1}{1+\hat v_j}\right)
\sum_{\substack{
S \subseteq \{1,\dots,\DimD-\DimK\} \\
|S| \text{ even}
}}
\prod_{j\in S} \frac{v_j}{1+\hat v_j}
\simeq 1
\end{equation}
on the domain of integration (here we use the fact  that $v_j,|\hat v_j| \in[0,\bar\epsilon^2]$ for all $j \in \{1,\ldots,\DimD-\DimK\}$ and $\bar\epsilon < 1$), so
\[\begin{split}
&
i^{\DimD-1-2\alpha}
\int_{\substack{
\max\{i^{-1},|x|\}^2 \lesssim V \leq \bar\epsilon^2 
}}
V^{-(\DimD-2\alpha)/2} 
\prod_{j=1}^{\DimD-\DimK} \left|\bar x_{\DimD-j+1}^2 - \frac{v_j}{1+\hat v_j}\right|^{-1/2} 
\,dv 
\\
&\simeq
i^{\DimD-1-2\alpha}
\int_{\substack{
\max\{i^{-1},|x|\}^2 \lesssim |w| 
}}
|w|^{-(\DimD-2\alpha)/2} 
\prod_{j=1}^{\DimD-\DimK} \left|\bar x_{\DimD-j+1}^2 - w_j \right|^{-1/2} 
\,dw .
\end{split}\]

In order to conclude, it is enough to bound the last integral with a multiple of $\min\{i,|x|^{-1}\}^{\DimK-2\alpha}$.
To do this, it is convenient to split the domain of integration according to whether $w_j$ is larger or smaller than $2\bar x_{\DimD-j+1}^2$ for each $j=1,\dots,\DimD-\DimK$, and according to which $j$ corresponds to the maximum component $w_j$ of $w$. In other words,
\begin{multline*}
\int_{\substack{
\max\{i^{-1},|x|\}^2 \lesssim |w| 
}}
|w|^{-(\DimD-2\alpha)/2} 
\prod_{j=1}^{\DimD-\DimK} \left|\bar x_{\DimD-j+1}^2 - w_j \right|^{-1/2} 
\,dw 
\\
\leq
\sum_{J \subseteq \{1,\dots,\DimD-\DimK\}}
\sum_{\DDJJ = 1}^{\DimD-\DimK}
\int_{\substack{
\max\{i^{-1},|x|\}^2 \lesssim |w| \\
w_{\DDJJ} = \max_j w_j \\
w_j \geq 2 \bar x_{\DimD-j+1}^2 \ \forall j \in J \\
w_j \leq 2 \bar x_{\DimD-j+1}^2 \ \forall j \in J^\compl
}}
|w|^{-(\DimD-2\alpha)/2} 
\prod_{j=1}^{\DimD-\DimK} \left|\bar x_{\DimD-j+1}^2 - w_j \right|^{-1/2} 
\,dw ,
\end{multline*}
where $J^\compl = \{1,\dots,\DimD-\DimK\} \setminus J$.
We estimate separately each summand, depending on the choice of $\DDJJ \in \{1,\dots,\DimD-\DimK\}$ and $J \subseteq \{1,\dots,\DimD-\DimK\}$, noting that, in the respective domain of integration, $|\bar x_{\DimD-j+1}^2 - w_j |^{-1/2} \simeq w_j^{-1/2}$ for all $j \in J$.

Suppose first that $\DDJJ \in J$, and set $J' = J \setminus \{\DDJJ\}$. Then
\[\begin{split}
&
\int_{\substack{
\max\{i^{-1},|x|\}^2 \lesssim |w| \\
w_{\DDJJ} = \max_j w_j \\
w_j \geq 2 \bar x_{\DimD-j+1}^2 \ \forall j \in J \\
w_j \leq 2 \bar x_{\DimD-j+1}^2 \ \forall j \in J^\compl
}}
|w|^{-(\DimD-2\alpha)/2} 
\prod_{j=1}^{\DimD-\DimK} \left|\bar x_{\DimD-j+1}^2 - w_j \right|^{-1/2} 
\,dw 
\\
&\lesssim 
\int_{\max\{i^{-1},|x|\}^2 \lesssim w_{\DDJJ}} w_{\DDJJ}^{-(\DimD-2\alpha)/2-1/2} \left(\prod_{j\in J'} \int_{w_j \leq w_{\DDJJ}} w_j^{-1/2} \,dw_j \right) \,dw_{\DDJJ}
\\
&\qquad\qquad \times
\left(
\prod_{j\in J^\compl}
\int_{w_j \leq 2\bar x_{\DimD-j+1}^2} \left|\bar x_{\DimD-j+1}^2 - w_j \right|^{-1/2} \,dw_j
\right)
\\
&\lesssim 
\left( \prod_{j\in J^\compl} |x_{\DimD-j+1}| \right)
\int_{\max\{i^{-1},|x|\}^2 \lesssim w_{\DDJJ}} w_{\DDJJ}^{-(\DimD-|J|-2\alpha)/2-1} \,dw_{\DDJJ} \\
&\lesssim 
|x|^{|J^\compl|}
\max\{i^{-1},|x|\}^{|J|-\DimD+2\alpha}
\leq 
\max\{i^{-1},|x|\}^{-\DimK+2\alpha} = \min\{i,|x|^{-1}\}^{\DimK-2\alpha},
\end{split}\]
which is the desired estimate. Here we used that $\DimD-|J|-2\alpha \geq \DimK-2\alpha > 0$.

Suppose instead that $\DDJJ \notin J$. In this range, $|x|^2 \lesssim \max\{i^{-1},|x|\}^2 \lesssim |w| \simeq w_{\DDJJ} \lesssim x_{\DimD-\DDJJ+1}^2 \leq |x|^2$, whence $w_{\DDJJ} \simeq |w|\simeq \max\{i^{-1},|x|\}^2$. So
\[\begin{split}
&
\int_{\substack{
\max\{i^{-1},|x|\}^2 \lesssim |w| \\
w_{\DDJJ} = \max_j w_j \\
w_j \geq 2 \bar x_{\DimD-j+1}^2 \ \forall j \in J \\
w_j \leq 2 \bar x_{\DimD-j+1}^2 \ \forall j \in J^\compl
}}
|w|^{-(\DimD-2\alpha)/2} \prod_{j=1}^{\DimD-\DimK} \left|\bar x_{\DimD-j+1}^2 - w_j \right|^{-1/2} \,dw 
\\
&\lesssim 
\max\{i^{-1},|x|\}^{-(\DimD-2\alpha)}
\left( \prod_{j\in J^\compl} \int_{w_j \leq 2 \bar x_{\DimD-j+1}^2} \left|\bar x_{\DimD-j+1}^2 - w_j \right|^{-1/2} \,dw_j \right) \\
&\qquad\qquad\times
\left(
\prod_{j \in J} \int_{w_j \lesssim \max\{i^{-1},|x|\}^2} w_j^{-1/2} \,dw_j \right) \\
&\lesssim
\max\{i^{-1},|x|\}^{-(\DimD-2\alpha)+|J|} |x|^{|J^\compl|}
\leq
\min\{i,|x|^{-1}\}^{\DimK-2\alpha},
\end{split}\]
and we are done.
\end{proof}

\subsection{The elliptic regime}

We now discuss the proof of Proposition \ref{prp:indici_bassi_primo_peso_secondo_stratopiu}. We first observe that a straightforward iteration of Proposition \ref{prp:stime-classiche}\ref{en:stime_somma} yields the following estimate.

\begin{lemma}\label{lem:somma_ripetuta}
Fix $\DimD \in\NN$, $\DimD \geq 2$, and $s \in \NN$, $1 \leq s \leq \DimD-1$.
For all $\ell_\DimD\in\NN_\DimD$ and all $(x_\DimD,\dots,x_{s+1}) \in [-1,1]^{\DimD-s}$,
\[
\sum_{\substack{(\ell_{\DimD-1},\dots,\ell_s) \in J_{\DimD-1}^{(s)} \\
\ell_{\DimD-1} \leq \ell_\DimD
}}
\ell_s^{s-1}
\bigl|\tX_{\ell_\DimD,\ell_{\DimD-1}}^\DimD(x_\DimD) \bigr|^2
\ldots
\bigl|\tX_{\ell_{s+1}, \ell_{s}}^{s+1}(x_{s+1}) \bigr|^2
\lesssim \ell_\DimD^{\DimD-1}.
\]
\end{lemma}

\begin{proof}[Proof of Proposition \ref{prp:indici_bassi_primo_peso_secondo_stratopiu}]
Due to the symmetry property of Jacobi polynomials \eqref{eq:symmetry}, we may restrict to $x \in [0,1]^{\DimD-\DimK}$.
Since
\[
\lambda_{\ell,m}^{\DimD,\DimK} + \frac{\DimD^2-\DimK^2+2(\DimK-\DimD)}{4} = \ell^2 - m^2
\]
for $(\ell,m) \in I_\DimD^{(\DimK)}$, it suffices to prove the estimate
\begin{equation}\label{eq:indici_bassi_primo_peso_secondo_strato_ausiliaria_dk}
\sum_{\substack{
(\ell_\DimD,\dots,\ell_\DimK) \in J_\DimD^{(\DimK)} \\
\ell_\DimK \leq \epsilon \ell_\DimD\\
\ell_\DimD^2 - \ell_\DimK^2 \in [i^2,(i+1)^2]}}
\ell_\DimK^{\DimK-1}
\bigl|\tX_{\ell_\DimD,\ell_{\DimD-1}}^\DimD(x_\DimD) \bigr|^2
\ldots
\bigl|\tX_{\ell_{\DimK+1},\ell_{\DimK}}^{\DimK+1}(x_{\DimK+1}) \bigr|^2
\lesssim_{\epsilon} i^{\DimD-1}
\end{equation}
for all $x \in [0,1]^{\DimD-\DimK}$ and $i \in \NN \setminus \{0\}$.
Indeed, the condition $\lambda_{\ell_\DimD,\ell_\DimK}^{\DimD,\DimK} \in [i^2,(i+1)^2]$ implies $\lambda_{\ell_\DimD,\ell_\DimK}^{\DimD,\DimK} +
(\DimD^2-\DimK^2+2(\DimK-\DimD))/{4} \in [i^2,(i+\DimD-1)^2]$, so \eqref{eq:indici_bassi_primo_peso_secondo_strato_dk} follows by applying \eqref{eq:indici_bassi_primo_peso_secondo_strato_ausiliaria_dk} $h \defeq \lceil (\DimD+\DimK-2)(\DimD-\DimK)/4 \rceil$ times, in correspondence of $i$, $i+1$, $\ldots$, $i+h-1$, respectively.
Note also that, since $\ell_\DimK \leq \epsilon \ell_\DimD$, we have $\lambda_{\ell,m}^{\DimD,\DimK} = \lambda_{\ell_\DimD,\ell_\DimK}^{\DimD,\DimK} \simeq \ell_\DimD^2$, whence
\[
\ell_\DimD \simeq i.
\]

We first consider the terms in the sum with $\ell_\DimK = 0$ (observe that this may happen only for $\DimK=1$). The condition $\ell_\DimD^2 \in [i^2,(i+1)^2]$ uniquely determines the value of $\ell_\DimD$. Using the estimate in Proposition \ref{prp:stime-classiche}\ref{en:stime_unif} to bound $\tX_{\ell_{\DimK+1},0}^{\DimK+1}(x_{\DimK+1})$ in the left-hand side of \eqref{eq:indici_bassi_primo_peso_secondo_strato_ausiliaria_dk} and then applying Lemma \ref{lem:somma_ripetuta}, we obtain
\[\begin{split}
&\sum_{\substack{
(\ell_\DimD,\dots,\ell_{\DimK+1}) \in J_\DimD^{(\DimK+1)} \\
\ell_\DimD^2 \in [i^2,(i+1)^2]}}
\bigl|\tX_{\ell_\DimD,\ell_{\DimD-1}}^\DimD(x_\DimD) \bigr|^2
\ldots
\bigl|\tX_{\ell_{\DimK+2},\ell_{\DimK+1}}^{\DimK+2}(x_{\DimK+2}) \bigr|^2 
\bigl|\tX_{\ell_{\DimK+1},0}^{\DimK+1}(x_{\DimK+1}) \bigr|^2 \\
&\lesssim  \sum_{\substack{
\ell_\DimD \in \NN_\DimD \\ 
\ell_\DimD^2 \in [i^2,(i+1)^2]
}}
\sum_{\substack{
(\ell_{\DimD-1},\dots,\ell_{\DimK+1}) \in J_{\DimD-1}^{(\DimK+1)} \\
{\ell_{\DimD-1} \leq \ell_\DimD}
}}
\ell_{\DimK+1}^{\DimK}
\bigl|\tX_{\ell_\DimD,\ell_{\DimD-1}}^\DimD(x_\DimD) \bigr|^2
\ldots 
\bigl|\tX_{\ell_{\DimK+2},\ell_{\DimK+1}}^{\DimK+2}(x_{\DimK+2}) \bigr|^2
\\
&\lesssim
\sum_{\substack{
\ell_\DimD \in \NN_\DimD \\ 
\ell_\DimD^2 \in [i^2,(i+1)^2]
}}
\ell_\DimD^{\DimD-1} \\
&\lesssim 
i^{\DimD-1}.
\end{split}\]
In what follows, we shall therefore assume $\ell_\DimK > 0$.

\medskip

Define $y_j \defeq \sqrt{1-x_j^2}$ for $j=\DimK+1,\dots,\DimD$.
Let us first consider the case where
\[
y_{\DDJJ} \le b_{\ell_{\DDJJ}, \ell_{\DDJJ-1}}/(2e) \qquad\text{for some } \DDJJ \in \{\DimK+1,\ldots,\DimD\}.
\]
By \eqref{eq:est_combinata}, in this case,
\[
|\tX_{\ell_{\DDJJ}, \ell_{\DDJJ-1}}^{\DDJJ}(x_{\DDJJ})|^2 
\lesssim \ell^{\DDJJ-1}_{\DDJJ} \ 
e^{-2c \ell_{\DDJJ-1}},
\]
for a suitable $c \in (0,\infty)$, whence
\[\begin{split}
&\sum_{\substack{
(\ell_\DimD,\dots,\ell_{\DimK}) \in J_\DimD^{(\DimK)} \\
0 < \ell_\DimK \leq \epsilon \ell_\DimD\\
\ell_\DimD^2 - \ell_\DimK^2 \in [i^2,(i+1)^2] \\
y_\DimK \leq b_{\ell_{\DDJJ},\ell_{\DDJJ-1}}/(2e)
}}
\ell_\DimK^{\DimK-1}
\bigl|\tX_{\ell_\DimD,\ell_{\DimD-1}}^\DimD(x_\DimD) \bigr|^2
\ldots
\bigl|\tX_{\ell_{\DimK+1},\ell_{\DimK}}^{\DimK+1}(x_{\DimK+1}) \bigr|^2 \\
&\lesssim
\sum_{\ell_{\DDJJ-1}\in\NN_{\DDJJ-1} }
e^{-2c \ell_{\DDJJ-1}} \\
&\qquad\times
\sum_{\substack{
(\ell_{\DDJJ-2},\dots,\ell_{\DimK}) \in J_{\DDJJ-2}^{(\DimK)} \\
\ell_{\DDJJ-2} \leq \ell_{\DDJJ-1} \\
\ell_\DimK \lesssim i
}}
\ell_\DimK^{\DimK-1}
\bigl|\tX_{\ell_{\DDJJ-1}, \ell_{\DDJJ-2}}^{\DDJJ-1}(x_{\DDJJ-1}) \bigr|^2
\ldots
\bigl|\tX_{\ell_{\DimK+1},\ell_{\DimK}}^{\DimK+1}(x_{\DimK+1}) \bigr|^2
\\
&\qquad\times
\sum_{\substack{
(\ell_\DimD,\dots,\ell_{\DDJJ}) \in J_\DimD^{(\DDJJ)} \\
\ell_{\DDJJ} \geq \ell_{\DDJJ-1} \\
\ell_\DimD^2 \in[\ell_\DimK^2+i^2, \ell_\DimK^2+(i+1)^2]
}}
\bigl| \tX_{\ell_\DimD,\ell_{\DimD-1}}^\DimD (x_\DimD) \bigr|^2
\ldots
\bigl|\tX_{\ell_{\DDJJ+1}, \ell_{\DDJJ}}^{J+1} (x_{\DDJJ+1}) \bigr|^2
\ell_{\DDJJ}^{\DDJJ-1},
\end{split}\]
where we also extended the sum in $\ell_{\DDJJ-1}$ to all $\NN_{\DDJJ-1}$.
As already observed, due to the condition $\ell_\DimD^2 \in [\ell_\DimK^2+i^2,\ell_\DimK^2+(i+1)^2]$, for a fixed $\ell_\DimK \lesssim i$ the sum over $\ell_\DimD$ essentially contains only one term and $\ell_\DimD \simeq i$.
Thus, by applying Lemma \ref{lem:somma_ripetuta} first to the sum over $\ell_{\DDJJ},\cdots,\ell_{\DimD-1}$ and then to the sum over $\ell_{\DDJJ}-2,\dots,\ell_\DimK$, we get
\[\begin{split}
&\sum_{\substack{
(\ell_\DimD,\dots,\ell_{\DimK}) \in J_\DimD^{(\DimK)} \\
0 < \ell_\DimK \leq \epsilon \ell_\DimD \\
\ell_\DimD^2 - \ell_\DimK^2 \in [i^2,(i+1)^2] \\
y_\DimK \leq b_{\ell_{\DDJJ},\ell_{\DDJJ-1}}/(2e)}}
\ell_\DimK^{\DimK-1}
\bigl|\tX_{\ell_\DimD,\ell_{\DimD-1}}^\DimD(x_\DimD) \bigr|^2
\ldots
\bigl|\tX_{\ell_{\DimK+1},\ell_\DimK}^{\DimK+1}(x_{\DimK+1}) \bigr|^2 \\
&\lesssim
i^{\DimD-1} \sum_{\ell_{\DDJJ-1}\in\NN_{\DDJJ-1}}
e^{-2c \ell_{\DDJJ-1}} \\
&\qquad\times
\sum_{\substack{
(\ell_{\DDJJ-2},\dots,\ell_{\DimK}) \in J_{\DDJJ-2}^{(\DimK)} \\
\ell_{\DDJJ-2} \leq \ell_{\DDJJ-1}
}}
\ell_\DimK^{\DimK-1}
\bigl|\tX_{\ell_{\DDJJ-1}, \ell_{\DDJJ-2}}^{\DDJJ-1}(x_{\DDJJ-1}) \bigr|^2
\ldots
\bigl|\tX_{\ell_{\DimK+1},\ell_{\DimK}}^{\DimK+1}(x_{\DimK+1}) \bigr|^2
\\
&\lesssim
i^{\DimD-1} \sum_{\ell_{\DDJJ-1}\in\NN_{\DDJJ-1} }
e^{-2c \ell_{\DDJJ-1}} \ell_{\DDJJ-1} \lesssim i^{\DimD-1}.
\end{split}\]

\medskip
From here on, we shall assume
\[
y_j > b_{\ell_j, \ell_{j-1}}/(2e) \qquad\text{for all } j \in \{\DimK+1,\ldots,\DimD\}.
\]
Note that this implies $y_j > 0$ for $j=\DimK+1,\dots,\DimD$.
We call $\DDKK$ the smallest index in $\{\DimK,\ldots,\DimD\}$ for which one has
\[
b_{\ell_j,\ell_{j-1}}/(2e) < y_j \leq 2 b_{\ell_j,\ell_{j-1}} \qquad \text{for all } j > \DDKK.
\]
Note that the above inequality implies that
\[
\ell_{j-1} \simeq \ell_j y_j \qquad \text{for all } j > \DDKK,
\]
and moreover, by Corollary \ref{cor:est_combinata},
\begin{equation}\label{eq:altre}
\bigl|\tX_{\ell_j,\ell_{j-1}}^j (x_j) \bigr|^2
\lesssim
y^{-(j-2)}_j \ \Xi(y_j^2, \ell_{j-1}/\ell_j^2, \ell_{j-1}^2/\ell_j^2)
\end{equation}
for $\DDKK < j \leq \DimD$, where $\Xi$ was defined in \eqref{eq:def_Xi}.

Assume first that $\DDKK > \DimK$. Then $y_\DDKK > 2b_{\ell_\DDKK,\ell_{\DDKK-1}}$, that is,
\[
\ell_{\DDKK-1} < \frac{1}{2} y_\DDKK \ell_\DDKK.
\]
whence, by Corollary \ref{cor:est_combinata},
\[
|\tX_{\ell_\DDKK,\ell_{\DDKK-1}}^\DDKK (x_\DDKK)|^2 \lesssim \ y_\DDKK^{-(\DDKK-1)} .
\]
Hence
\[\begin{split}
&\sum_{\substack{
(\ell_\DimD,\dots,\ell_{\DimK}) \in J_\DimD^{(\DimK)} \\
0 < \ell_\DimK \leq \epsilon \ell_\DimD\\
\ell_\DimD^2 - \ell_\DimK^2 \in [i^2,(i+1)^2] \\
b_{\ell_j,\ell_{j-1}}/(2e) < y_j \leq 2 b_{\ell_j,\ell_{j-1}} \forall j>\DDKK \\
y_\DDKK > 2 b_{\ell_\DDKK,\ell_{\DDKK-1}} 
}}
\ell_\DimK^{\DimK-1}
\bigl|\tX_{\ell_\DimD,\ell_{\DimD-1}}^\DimD(x_\DimD) \bigr|^2
\ldots
\bigl|\tX_{\ell_{\DimK+1},\ell_\DimK}^{\DimK+1} (x_{\DimK+1}) \bigr|^2 \\
&\lesssim 
\sum_{\substack{
(\ell_{\DDKK-1},\dots,\ell_\DimK) \in J_{\DDKK-1}^{(\DimK)} \\
\ell_{\DDKK-1} \lesssim i y_\DDKK y_{\DDKK+1} \cdots y_\DimD
}}
\ell_\DimK^{\DimK-1} 
\bigl|\tX_{\ell_{\DDKK-1},\ell_{\DDKK-2}}^\DimD(x_{\DDKK-1}) \bigr|^2
\ldots
\bigl|\tX_{\ell_{\DimK+1},\ell_\DimK}^{\DimK+1} (x_{\DimK+1}) \bigr|^2 \\
&\quad\times 
\sum_{\substack{
(\ell_\DimD,\dots,\ell_\DDKK) \in J_\DimD^{(\DDKK)} \\
\ell_{j-1} \simeq \ell_j y_j \ \forall j > \DDKK \\
\ell_\DimD \in [\sqrt{i^2+\ell_\DimK^2},\sqrt{(i+1)^2+\ell_\DimK^2}]
}} 
y_\DDKK^{-(\DDKK-1)}
\prod_{j=\DDKK+1}^\DimD y_j^{-(j-2)} \Xi(y_j^2,\ell_{j-1}/\ell_j^2,\ell_{j-1}^2/\ell_j^2).
\end{split}\]
Note now that, for $j=\DDKK+1,\dots,\DimD$,
\[
\nabla_{(\ell_\DimD,\dots,\ell_\DDKK)} (\ell_{j-1}/\ell_j^2), \nabla_{(\ell_\DimD,\dots,\ell_\DDKK)} (\ell_{j-1}^2/\ell_j^2), \lesssim \ell_{j-1}/\ell_j^2
\]
in the range of summation; moreover the interval $[\sqrt{i^2+\ell_\DimK^2},\sqrt{(i+1)^2+\ell_\DimK^2}]$ has length $\simeq 1$ and its endpoints are $\simeq i$, because $\ell_\DimK \lesssim i$. Hence, in view of Lemma \ref{lem:der_Xi}, we can apply \cite[Lemma 5.7]{DM} to the inner sum and obtain
\[\begin{split}
&\sum_{\substack{
(\ell_\DimD,\dots,\ell_\DDKK) \in J_\DimD^{(\DDKK)} \\
\ell_{j-1} \simeq \ell_j y_j \ \forall j>\DDKK \\
\ell_\DimD \in [\sqrt{i^2+\ell_\DimK^2},\sqrt{(i+1)^2+\ell_\DimK^2}]
}} 
y_\DDKK^{-(\DDKK-1)}
\prod_{j=\DDKK+1}^\DimD y_j^{-(j-2)} \Xi(y_j^2,\ell_{j-1}/\ell_j^2,\ell_{j-1}^2/\ell_j^2) \\
&\lesssim
\left( y_\DDKK^{-1} \prod_{j=\DDKK}^\DimD y_j^{-(j-2)} \right)
\int_{\substack{
\ell_\DimD \in [\sqrt{i^2+\ell_\DimK^2},\sqrt{(i+1)^2+\ell_\DimK^2}] \\
\ell_{j-1} \simeq \ell_j y_j \ \forall j>\DDKK 
}}
\prod_{j=\DDKK+1}^\DimD \left|y_j^2 - \frac{\ell_{j-1}^2}{\ell_j^2} \right|^{-1/2} \,d\ell_\DDKK \cdots \,d\ell_\DimD.
\end{split}\]
The change of variables $t_{j-1} = \ell_{j-1}/(\ell_j y_j)$, $j=\DDKK+1,\dots,\DimD$, then gives
\[\begin{split}
&\int_{\substack{
\ell_\DimD \in [\sqrt{i^2+\ell_\DimK^2},\sqrt{(i+1)^2+\ell_\DimK^2}] \\
\ell_{j-1} \simeq \ell_j y_j \ \forall j>\DDKK 
}}
\prod_{j=\DDKK+1}^\DimD \left|y_j^2 - \frac{\ell_{j-1}^2}{\ell_j^2} \right|^{-1/2} \,d\ell_\DDKK \cdots \,d\ell_\DimD \\
&\lesssim \int_{\ell_\DimD \in [\sqrt{i^2+\ell_\DimK^2},\sqrt{(i+1)^2+\ell_\DimK^2}]} 
\int_{t_\DDKK,\dots,t_{\DimD-1} \simeq 1} \prod_{j=\DDKK+1}^\DimD \frac{y_{j+1} \cdots y_\DimD \, i}{|1-t_{j-1}^2|^{1/2}}  \,dt_\DDKK \cdots \,dt_{\DimD-1} \,d\ell_\DimD \\
&\simeq \prod_{j=\DDKK+1}^\DimD (y_{j+1} \cdots y_\DimD \, i) = i^{\DimD-\DDKK} \prod_{j=\DDKK+2}^\DimD y_j^{j-\DDKK-1},
\end{split}\]
whence
\begin{multline*}
\sum_{\substack{
(\ell_\DimD,\dots,\ell_\DDKK) \in J_\DimD^{(\DDKK)} \\
\ell_{j-1} \simeq \ell_j y_j \ \forall j>\DDKK \\
\ell_\DimD \in [\sqrt{i^2+\ell_\DimK^2},\sqrt{(i+1)^2+\ell_\DimK^2}]
}} 
y_\DDKK^{-(\DDKK-1)}
\prod_{j=\DDKK+1}^\DimD y_j^{-(j-2)} \Xi(y_j^2,\ell_{j-1}/\ell_j^2,\ell_{j-1}^2/\ell_j^2) \\
\lesssim i^{\DimD-\DDKK} (y_\DDKK \cdots y_\DimD)^{1-\DDKK}
\end{multline*}
and
\[\begin{split}
&\sum_{\substack{
(\ell_\DimD,\dots,\ell_{\DimK}) \in J_\DimD^{(\DimK)} \\
0 < \ell_\DimK \leq \epsilon \ell_\DimD \\
\ell_\DimD^2 - \ell_\DimK^2 \in [i^2,(i+1)^2] \\
b_{\ell_j,\ell_{j-1}}/(2e) < y_j \leq 2 b_{\ell_j,\ell_{j-1}} \forall j>\DDKK \\
y_\DDKK > 2 b_{\ell_\DDKK,\ell_{\DDKK-1}} 
}}
\ell_\DimK^{\DimK-1}
\bigl|\tX_{\ell_\DimD,\ell_{\DimD-1}}^\DimD(x_\DimD) \bigr|^2
\ldots
\bigl|\tX_{\ell_{\DimK+1},\ell_\DimK}^{\DimK+1}(x_{\DimK+1}) \bigr|^2 \\
&\lesssim 
i^{\DimD-\DDKK} (y_\DDKK \cdots y_\DimD)^{1-\DDKK} \\
&\qquad\times
\sum_{\substack{
(\ell_{\DDKK-1},\dots,\ell_\DimK) \in J_{\DDKK-1}^{(\DimK)} \\
\ell_{\DDKK-1} \lesssim i y_\DDKK y_{\DDKK+1} \cdots y_\DimD
}}
\ell_\DimK^{\DimK-1} 
\bigl|\tX_{\ell_{\DDKK-1},\ell_{\DDKK-2}}^{\DDKK-1} (x_{\DDKK-1}) \bigr|^2
\ldots
\bigl|\tX_{\ell_{\DimK+1},\ell_\DimK}^{\DimK+1} (x_{\DimK+1}) \bigr|^2 \\
&\lesssim 
i^{\DimD-\DDKK} (y_\DDKK \cdots y_\DimD)^{1-\DDKK}
\sum_{\ell_{\DDKK-1} \lesssim i y_\DDKK y_{\DDKK+1} \cdots y_\DimD} \ell_{\DDKK-1}^{\DDKK-2} \lesssim i^{\DimD-1},
\end{split}\]
where Lemma \ref{lem:somma_ripetuta} was applied to the sum in $(\ell_\DDKK,\dots,\ell_\DimK)$ and the fact that $\DDKK \geq \DimK+1\geq 2$ was used.

\medskip

We now consider the case $\DDKK=\DimK$.
Here, by \eqref{eq:altre},
\[\begin{split}
&\sum_{\substack{
(\ell_\DimD,\dots,\ell_{\DimK}) \in J_\DimD^{(\DimK)} \\
0 < \ell_\DimK \leq \epsilon \ell_\DimD \\
\ell_\DimD^2 - \ell_\DimK^2 \in [i^2,(i+1)^2] \\
b_{\ell_j,\ell_{j-1}}/(2e) < y_j \leq 2 b_{\ell_j,\ell_{j-1}} \forall j>\DimK 
}}
\ell_\DimK^{\DimK-1}
\bigl|\tX_{\ell_\DimD,\ell_{\DimD-1}}^\DimD(x_\DimD) \bigr|^2
\ldots
\bigl|\tX_{\ell_{\DimK+1},\ell_\DimK}^{\DimK+1}(x_{\DimK+1}) \bigr|^2 \\
&\lesssim 
\sum_{\substack{
(\ell_\DimD,\dots,\ell_\DimK) \in J_\DimD^{(\DimK)} \\
\ell_{j-1} \simeq \ell_j y_j \ \forall j>\DimK \\
\ell_\DimD \in [\sqrt{i^2+\ell_\DimK^2},\sqrt{(i+1)^2+\ell_\DimK^2}]
}} 
\ell_\DimK^{\DimK-1} 
\prod_{j=\DimK+1}^\DimD y_j^{-(j-2)} \Xi(y_j^2,\ell_{j-1}/\ell_j^2,\ell_{j-1}^2/\ell_j^2) \\
&\lesssim 
i^{\DimK-1} 
\left( \prod_{j=\DimK+1}^\DimD y_j^{\DimK+1-j} \right)
\sum_{\substack{
(\ell_{\DimD-1},\dots,\ell_\DimK) \in J_{\DimD-1}^{(\DimK)} \\
\ell_j \simeq y_{j+1} \cdots y_\DimD i \ \forall \DimK \leq j < \DimD 
}}
\prod_{j=\DimK+1}^{\DimD-1} \Xi(y_j^2,\ell_{j-1}/\ell_j^2,\ell_{j-1}^2/\ell_j^2)
\\
&\quad\times
\sum_{\substack{
\ell_\DimD \in \NN_\DimD \\
\ell_\DimD \in [\sqrt{i^2+\ell_\DimK^2},\sqrt{(i+1)^2+\ell_\DimK^2}]
}}
\Xi(y_\DimD^2,\ell_{\DimD-1}/\ell_\DimD^2,\ell_{\DimD-1}^2/\ell_\DimD^2) \\
&\lesssim 
i^{\DimK-1} 
\left( \prod_{j=\DimK+1}^\DimD y_j^{\DimK+1-j} \right)
\sum_{\substack{
(\ell_{\DimD-1},\dots,\ell_\DimK) \in J_{\DimD-1}^{(\DimK)} \\
\ell_j \simeq y_{j+1} \cdots y_\DimD i \ \forall \DimK \leq j < \DimD 
}}
\prod_{j=\DimK+1}^{\DimD-1} \Xi(y_j^2,\ell_{j-1}/\ell_j^2,\ell_{j-1}^2/\ell_j^2)
\\
&\quad\times
\int_{
\ell_\DimD \in [\sqrt{i^2+\ell_\DimK^2},\sqrt{(i+1)^2+\ell_\DimK^2}]
}
\Xi(y_\DimD^2,\ell_{\DimD-1}/\ell_\DimD^2,\ell_{\DimD-1}^2/\ell_\DimD^2) \,d\ell_\DimD,
\end{split}\]
where the last inequality follows from \cite[Lemma 4.1]{CaCiaMa} together with Lemma \ref{lem:der_Xi}, the fact that
\[
|\partial_{\ell_\DimD} (\ell_{\DimD-1}/\ell_\DimD^2)|, |\partial_{\ell_\DimD} (\ell_{\DimD-1}^2/\ell_\DimD^2)| \lesssim \ell_{\DimD-1}/\ell_\DimD^2
\]
in the range of summation and the fact that (since $\ell_\DimK \lesssim i$) the length of the interval $[\sqrt{\ell_\DimK^2+i^2},\sqrt{\ell_\DimK^2+(i+1)^2}]$ is $\simeq 1$.

The change of variables $u = \sqrt{\ell_\DimD^2 - \ell_\DimK^2}$ in the inner integral then gives
\[\begin{split}
&\sum_{\substack{
(\ell_\DimD,\dots,\ell_{\DDKK}) \in J_\DimD^{(\DimK)} \\
0 < \ell_\DimK \leq \epsilon \ell_\DimD \\
\ell_\DimD^2 - \ell_\DimK^2 \in [i^2,(i+1)^2] \\
b_{\ell_j,\ell_{j-1}}/(2e) < y_j \leq 2 b_{\ell_j,\ell_{j-1}} \forall j>\DimK 
}}
\ell_\DimK^{\DimK-1}
\bigl|\tX_{\ell_\DimD,\ell_{\DimD-1}}^\DimD(x_\DimD) \bigr|^2
\ldots
\bigl|\tX_{\ell_{\DimK+1},\ell_\DimK}^{\DimK+1}(x_{\DimK+1}) \bigr|^2 \\
&\lesssim 
i^{\DimK-1} 
\left( \prod_{j=\DimK+1}^\DimD y_j^{\DimK+1-j} \right)
\sum_{\substack{
(\ell_{\DimD-1},\dots,\ell_\DimK) \in J_{\DimD-1}^{(\DimK)} \\
\ell_j \simeq y_{j+1} \cdots y_\DimD i \ \forall \DimK \leq j <\DimD 
}}
\prod_{j=\DimK+1}^{\DimD-1} \Xi(y_j^2,\ell_{j-1}/\ell_j^2,\ell_{j-1}^2/\ell_j^2)
\\
&\quad\times
\int_i^{i+1}
\Xi(y_\DimD^2,\ell_{\DimD-1}/(u^2+\ell_\DimK^2),\ell_{\DimD-1}^2/(u^2+\ell_\DimK^2)) \,du \\
&=
i^{\DimK-1} 
\left( \prod_{j=\DimK+1}^\DimD y_j^{\DimK+1-j} \right)
\int_i^{i+1} 
\sum_{\substack{
(\ell_{\DimD-1},\dots,\ell_\DimK) \in J_{\DimD-1}^{(\DimK)} \\
\ell_j \simeq y_{j+1} \cdots y_\DimD i \ \forall \DimK \leq j < \DimD 
}}
\tilde\Xi(u,\vec y,\vec\ell) \,du,
\end{split}\]
where $\vec \ell = (\ell_{\DimD-1},\dots,\ell_\DimK)$, $\vec y = (y_\DimD,\dots,y_{\DimK+1})$, and
\[
\tilde\Xi(u,\vec y,\vec\ell)
=
\Xi(y_\DimD^2,\ell_{\DimD-1}/(u^2+\ell_\DimK^2),\ell_{\DimD-1}^2/(u^2+\ell_\DimK^2))
\prod_{j=\DimK+1}^{\DimD-1} \Xi(y_j^2,\ell_{j-1}/\ell_j^2,\ell_{j-1}^2/\ell_j^2).
\]

We now observe that, since $u\in[i,i+1]$,
\[
|\nabla_{\vec\ell} \ell_{\DimD-1}/(u^2+\ell_\DimK^2)|, |\nabla_{\vec\ell} \ell_{\DimD-1}^2/(u^2+\ell_\DimK^2)| \lesssim \ell_{\DimD-1}/(u^2+\ell_\DimK^2)
\]
and
\[
|\nabla_{\vec\ell} \ell_{j-1}/\ell_j^2|, |\nabla_{\vec\ell} \ell_{j-1}^2/\ell_j^2| \lesssim \ell_{j-1}/\ell_j^2 \qquad\text{for } j=\DimK+1,\dots,\DimD-1
\]
on the range of summation. Thanks to Lemma \ref{lem:der_Xi}, we can apply \cite[Lemma 5.7]{DM} to majorize the inner sum  with the corresponding integral and obtain that
\[\begin{split}
&\sum_{\substack{
(\ell_\DimD,\dots,\ell_{\DimK}) \in J_\DimD^{(\DimK)} \\
0<\ell_\DimK \leq \epsilon \ell_\DimD\\
\ell_\DimD^2 - \ell_\DimK^2 \in [i^2,(i+1)^2] \\
b_{\ell_j,\ell_{j-1}}/(2e) < y_j \leq 2 b_{\ell_j,\ell_{j-1}} \forall j>\DimK 
}}
\ell_\DimK^{\DimK-1}
\bigl|\tX_{\ell_\DimD,\ell_{\DimD-1}}^\DimD(x_\DimD) \bigr|^2
\ldots
\bigl|\tX_{\ell_{\DimK+1},\ell_\DimK}^{\DimK+1}(x_{\DimK+1}) \bigr|^2 \\
&\lesssim 
i^{\DimK-1} 
\left( \prod_{j=\DimK+1}^\DimD y_j^{\DimK+1-j} \right)
\int_i^{i+1} 
\int_{\substack{
\ell_j \simeq y_{j+1} \cdots y_\DimD i \ \forall \DimK \leq j < \DimD 
}}
\tilde\Xi(u,\vec y,\vec\ell) \,d\ell_\DimK \cdots \,d\ell_{\DimD-1} \,du \\
&\lesssim
i^{\DimK-1}  
\left( \prod_{j=\DimK+1}^\DimD y_j^{\DimK+1-j} \right)
\int_i^{i+1} 
\int_{\substack{
\ell_j \simeq y_{j+1} \cdots y_\DimD i \ \forall \DimK \leq j < \DimD 
}}
\left| y_\DimD^2 - \frac{\ell_{\DimD-1}^2}{u^2+\ell_\DimK^2} \right|^{-1/2} \\
&\quad\times \prod_{j=\DimK+1}^{\DimD-1} \left| y_j^2 - \frac{\ell_{j-1}^2}{\ell_j^2} \right|^{-1/2}
\,d\ell_\DimK \cdots \,d\ell_{\DimD-1} \,du,
\end{split}\]
The change of variables $\ell_j = u y_{j+1} \cdots y_\DimD \tau_j$, $j=\DimK,\dots,\DimD-1$, then gives
\[\begin{split}
&\sum_{\substack{
(\ell_\DimD,\dots,\ell_{\DimK}) \in J_\DimD^{(\DimK)} \\
0 < \ell_\DimK \leq \epsilon \ell_\DimD \\
\ell_\DimD^2 - \ell_\DimK^2 \in [i^2,(i+1)^2] \\
b_{\ell_j,\ell_{j-1}}/(2e) < y_j \leq 2 b_{\ell_j,\ell_{j-1}} \forall j>\DimK 
}}
\ell_\DimK^{\DimK-1}
\bigl|\tX_{\ell_\DimD,\ell_{\DimD-1}}^\DimD(x_\DimD) \bigr|^2
\ldots
\bigl|\tX_{\ell_{\DimK+1},\ell_\DimK}^{\DimK+1}(x_{\DimK+1}) \bigr|^2 \\
&\lesssim 
i^{\DimD-1}
\left( \prod_{j=\DimK+1}^\DimD y_j \right)
\int_i^{i+1} 
\int_{\substack{
\tau_\DimK,\dots,\tau_{\DimD-1} \simeq 1 
}}
\left| y_\DimD^2 - \frac{y_\DimD^2 \tau_{\DimD-1}^2}{1+y_{\DimK+1}^2 \cdots y_\DimD^2 \tau_\DimK^2} \right|^{-1/2} \\
&\quad\times \prod_{j=\DimK+1}^{\DimD-1} \left| y_j^2 - \frac{y_j^2 \tau_{j-1}^2}{\tau_j^2} \right|^{-1/2}
\,d\tau_\DimK \cdots \,d\tau_{\DimD-1} \,du \\
&=
i^{\DimD-1}
\int_{\substack{
\tau_\DimK,\dots,\tau_{\DimD-1} \simeq 1 
}}
\left| 1 - \frac{\tau_{\DimD-1}^2}{1+y_{\DimK+1}^2 \cdots y_\DimD^2 \tau_\DimK^2} \right|^{-1/2}  \\
&\qquad\times
\prod_{j=\DimK+1}^{\DimD-1} \left| 1 - \frac{\tau_{j-1}^2}{\tau_j^2} \right|^{-1/2}
\,d\tau_\DimK \cdots \,d\tau_{\DimD-1} \\
&\simeq
i^{\DimD-1}
\int_{\substack{
\tau_\DimK,\dots,\tau_{\DimD-1} \simeq 1 
}}
\left| \frac{1}{\tau_{\DimD-1}^2} + \frac{y_{\DimK+1}^2 \cdots y_\DimD^2 \tau_\DimK^2}{\tau_{\DimD-1}^2} - 1 \right|^{-1/2} \\
&\qquad\times
\prod_{j=\DimK+1}^{\DimD-1} \left| 1 - \frac{\tau_{j-1}^2}{\tau_j^2} \right|^{-1/2}
\,d\tau_\DimK \cdots \,d\tau_{\DimD-1} .
\end{split}\]
Finally, the change of variables 
\[
t_{\DimD-1} = \frac{1}{\tau_{\DimD-1}}, \qquad t_j = \frac{\tau_j}{\tau_{j+1}} \quad\text{for } j=\DimK,\dots,\DimD-2
\]
yields
\[\begin{split}
&\sum_{\substack{
(\ell_\DimD,\dots,\ell_{\DimK}) \in J_\DimD^{(\DimK)} \\
0 < \ell_\DimK \leq \epsilon \ell_\DimD \\
\ell_\DimD^2 - \ell_\DimK^2 \in [i^2,(i+1)^2] \\
b_{\ell_j,\ell_{j-1}}/(2e) < y_j \leq 2 b_{\ell_j,\ell_{j-1}} \forall j>\DimK 
}}
\ell_\DimK^{\DimK-1}
\bigl|\tX_{\ell_\DimD,\ell_{\DimD-1}}^\DimD(x_\DimD) \bigr|^2
\ldots
\bigl|\tX_{\ell_{\DimK+1},\ell_\DimK}^{\DimK+1} (x_{\DimK+1}) \bigr|^2 \\
&\lesssim 
i^{\DimD-1}
\int_{\substack{
t_\DimK,\dots,t_{\DimD-1} \simeq 1 
}}
\left| t_{\DimD-1}^2 + y_{\DimK+1}^2 \cdots y_\DimD^2 t_\DimK^2 \cdots t_{\DimD-2}^2 - 1 \right|^{-1/2}  \\
&\qquad\times
\prod_{j=\DimK+1}^{\DimD-1} \left| 1 - t_{j-1}^2 \right|^{-1/2}
\,dt_\DimK \cdots \,dt_{\DimD-1} \\
&\lesssim 
i^{\DimD-1}
\int_{\substack{
t_\DimK,\dots,t_{\DimD-2} \simeq 1 
}}
\prod_{j=\DimK+1}^{\DimD-1} \left| 1 - t_{j-1}^2 \right|^{-1/2}
\int_{|v|\lesssim 1} |v|^{-1/2} \,dv \,dt_\DimK \cdots \,dt_{\DimD-2}
\lesssim i^{\DimD-1},
\end{split}\]
and we are done.
\end{proof}

\section{The multiplier theorem}\label{s:Plancherel}
\subsection{The weighted Plancherel-type estimate}

By means of the previous estimates we shall prove a ``weighted Plancherel-type estimate'' for the Grushin operator $\opL_{\DimD,\DimK}$.  

For all $r \in (0,\infty)$, we define the weight $\weight_r :  \sfera^\DimD \times \sfera^\DimD \to [0,\infty)$ by
\begin{equation}\label{eq:weight}
\weight_r(\spnt{\omega}{\psi},\spnt{\omega'}{\psi'}) = \frac{|\psi|}{\max\{r,|\psi'|\}}.
\end{equation}

\begin{proposition}\label{prp:weighted_plancherel}
Let $\alpha \in [0,\DimK/2)$ and $N \in \NN \setminus \{0\}$. For all Borel functions $F : \RR \to \CC$ supported in $[0,N]$, and all $z' \in \sfera^\DimD$,
\[
\| (1+\weight_{N^{-1}}(\cdot,z'))^\alpha \, \Kern_{F(\sqrt{\opL_{\DimD,\DimK}})}(\cdot,z') \|_{L^2(\sfera)} 
\lesssim_\alpha V(z',N^{-1})^{-1/2} \|F(N \cdot)\|_{N,2} .
\]
\end{proposition}
\begin{proof}
We shall prove the apparently weaker estimate
\begin{equation}\label{eq:weighted_reduced}
\| \weight_{N^{-1}}(\cdot,z')^\alpha \, \Kern_{F(\sqrt{\opL_{\DimD,\DimK}})}(\cdot,z') \|_{L^2(\sfera)} 
\lesssim_\alpha V(z',N^{-1})^{-1/2} \|F(N \cdot)\|_{N,2}
\end{equation}
for all $z' \in \sfera^\DimD$.
Proposition \ref{prp:weighted_plancherel} follows by combining the estimate \eqref{eq:weighted_reduced} with the analogous one where $\alpha = 0$.

Following \eqref{eq:kernelformula}, we can decompose
\[\begin{split}
\Kern_{F(\sqrt{\opL_{\DimD,\DimK}})}
&= 
\sum_{\substack{
(\ell_\DimD,\dots,\ell_\DimK) \in J_\DimD^{(\DimK)} \\
}} 
F\left(\sqrt{\lambda_{\ell_\DimD,\ell_\DimK}^{\DimD,\DimK}}\right)   
 \,K_{\ell_\DimD,\dots,\ell_\DimK}^\DimD
 \\
&= 
\sum_{\substack{
(\ell_\DimD,\dots,\ell_\DimK) \in J_\DimD^{(\DimK)} \\
\ell_\DimK \leq \epsilon\ell_\DimD
}} 
+
\sum_{\substack{
(\ell_\DimD,\dots,\ell_\DimK) \in J_\DimD^{(\DimK)} \\
\ell_\DimK > \epsilon\ell_\DimD
}} 
\eqdef
K_1 + K_2,
\end{split}\]
where $\epsilon = \max\{1/2,(\DimK-1)/(\DimD-1)\} \in (0,1)$ and $\lambda^{\DimD,\DimK}_{\ell_\DimD,\ell_\DimK}$ is given by \eqref{eq:grushin_eigenvalue}.

We note that, due to the choice of $\epsilon$, for all $(\ell_\DimD,\dots,\ell_\DimK) \in J_\DimD^{(\DimK)}$ with $\ell_\DimK > \epsilon \ell_\DimD$,
\begin{equation}\label{eq:equivalence_ell}
\lambda^{\DimK+1}_{\ell_{\DimK+1}} \simeq \ell_{\DimK+1}^2 \simeq \ell_\DimD^2
\end{equation}
(see \eqref{eq:laplace_eigenvalues}). In particular, $K_2(\cdot,z') \perp \ker(\Delta_{\DimK+1})$, and moreover
\[
K_2(\cdot,z') = \opL_{\DimD,\DimK}^{-\alpha/2} \Delta_{\DimK+1}^{\alpha/2} K_{2,\alpha} (\cdot,z'),
\]
where
\[\begin{split}
K_{2,\alpha}
&=\sum_{\substack{
(\ell_\DimD,\dots,\ell_\DimK) \in J_\DimD^{(\DimK)} \\
b_{\ell_\DimD,\ell_\DimK} > \epsilon
}}
({\lambda_{\ell_\DimD,\ell_\DimK}^{\DimD,\DimK}}/\lambda_{\ell_{\DimK+1}}^{\DimK+1})^{\alpha/2}
F\left(\sqrt{\lambda_{\ell_\DimD,\ell_\DimK}^{\DimD,\DimK}}\right) 
 \,K_{\ell_\DimD,\dots,\ell_\DimK}^\DimD .
\end{split}\]
Hence
\[\begin{split}
&\| \weight_{N^{-1}}(\cdot,\spnt{\omega'}{\psi'} )^\alpha 
\Kern_{F(\sqrt{\opL_{\DimD,\DimK}})}(\cdot, \spnt{\omega'}{\psi'}) \|_{L^2(\sfera^\DimD)} \\
&\leq
\| \weight_{N^{-1}}(\cdot, \spnt{\omega'}{\psi'})^\alpha K_1(\cdot, \spnt{\omega'}{\psi'}) \|_{L^2(\sfera^\DimD)} +
\| \weight_{N^{-1}}(\cdot, \spnt{\omega'}{\psi'})^\alpha K_2(\cdot, \spnt{\omega'}{\psi'}) \|_{L^2(\sfera^\DimD)} \\
&\lesssim
\min\{N,|\psi'|^{-1}\}^\alpha \left[ \| K_1(\cdot,\spnt{\omega'}{\psi'}) \|_{L^2(\sfera^\DimD)} +
\| \tau_{\DimD,\DimK}^\alpha K_2(\cdot,\spnt{\omega'}{\psi'}) \|_{L^2(\sfera^\DimD)} \right] \\
&\leq
\min\{N,|\psi'|^{-1}\}^\alpha \left[ \| K_1(\cdot, \spnt{\omega'}{\psi'}) \|_{L^2(\sfera^\DimD)} +
\| K_{2,\alpha}(\cdot,\spnt{\omega'}{\psi'} ) \|_{L^2(\sfera^\DimD)} \right]
\end{split}\]
where $\tau_{\DimD,\DimK}$ is the function defined in \eqref{eq:taudk} and in the last step the Riesz-type estimate of Proposition \ref{prp:Riesz} was used.
In light of \eqref{eq:volume}, the estimate \eqref{eq:weighted_reduced} will follow from
\begin{align}
\| K_1(\cdot,\spnt{\omega'}{\psi'} ) \|_{L^2(\sfera^2)}^2 &\lesssim_\alpha N^\DimD \min\{ N, |\psi'|^{-1} \}^{\DimK-2\alpha} \|F(N\cdot)\|_{N,2}^2, \label{eq:est_bassi_orig}\\
\| K_{2,\alpha}(\cdot,\spnt{\omega'}{\psi'}) \|_{L^2(\sfera^\DimD)}^2 &\lesssim_\alpha N^\DimD \min\{ N, |\psi'|^{-1} \}^{\DimK-2\alpha} \|F(N\cdot)\|_{N,2}^2 \label{eq:est_alti}.
\end{align}
In fact, instead of \eqref{eq:est_bassi_orig}, we shall prove the stronger estimate
\begin{equation}\label{eq:est_bassi}
\| K_1(\cdot,\spnt{\omega'}{\psi'} ) \|_{L^2(\sfera^d)}^2 \lesssim N^\DimD \|F(N\cdot)\|_{N,2}^2.
\end{equation}

In view of \eqref{eq:kernel_l2norm} and \eqref{eq:equivalence_ell}, 
we rewrite
 \eqref{eq:est_alti} and \eqref{eq:est_bassi} as
\begin{multline*}
\sum_{\substack{
(\ell_\DimD,\dots,\ell_\DimK) \in J_\DimD^{(\DimK)} \\
\ell_\DimK > \epsilon\ell_\DimD
}}
\dimHarm_{\ell_\DimK}(\sfera^\DimK)
(\lambda_{\ell_\DimD,\ell_\DimK}^{\DimD,\DimK}/\ell_\DimD^2)^{\alpha}
\left|F\left(\sqrt{\lambda_{\ell_\DimD,\ell_\DimK}^{\DimD,\DimK}}\right)\right|^2
|X^\DimD_{\ell_\DimD,\dots,\ell_\DimK}(\psi)|^2
\\
\lesssim_\alpha N^{\DimD-1} \min\{ N, |\psi'|^{-1} \}^{\DimK-2\alpha}
\sum_{i=1}^N \sup_{\lambda \in [i-1,i]} |F(\lambda)|^2
\end{multline*}
and
\[
\sum_{\substack{(\ell_\DimD,\dots,\ell_\DimK) \in J_\DimD^{(\DimK)} \\
\ell_\DimK \leq \epsilon\ell_\DimD
}}
\dimHarm_{\ell_\DimK}(\sfera^\DimK)
\left|F\left(\sqrt{\lambda_{\ell_\DimD,\ell_\DimK}^{\DimD,\DimK}}\right)\right|^2 
|X^\DimD_{\ell_\DimD,\ldots,\ell_\DimK} (\psi)|^2
\lesssim_\alpha N^{\DimD-1} 
\sum_{i=1}^N \sup_{\lambda \in [i-1,i]} |F(\lambda)|^2.
\]

So it is enough to prove that
\[
\sum_{\substack{
(\ell_\DimD,\dots,\ell_\DimK) \in J_\DimD^{(\DimK)} \\
\ell_\DimK > \epsilon \ell_\DimD \\
\lambda_{\ell_\DimD,\ell_\DimK}^{\DimD,\DimK} \in [(i-1)^2,i^2] 
}}
\dimHarm_{\ell_\DimK}(\sfera^\DimK)
({\lambda_{\ell_\DimD,\ell_\DimK}^{\DimD,\DimK}}/\ell_\DimD^2)^{\alpha}
|X^\DimD_{\ell_\DimD, \ldots, \ell_\DimK}(\psi)|^2
\lesssim_\alpha N^{\DimD-1} \min\{ N, |\psi'|^{-1} \}^{\DimK-2\alpha}
\]
and
\[
\sum_{\substack{
(\ell_\DimD,\dots,\ell_\DimK) \in J^{(\DimK)}_\DimD \\ 
\ell_\DimK \leq \epsilon\ell_\DimD \\
\lambda_{\ell_\DimD,\ell_\DimK}^{\DimD,\DimK} \in [(i-1)^2,i^2]
}} 
\dimHarm_{\ell_\DimK}(\sfera^\DimK)
|X^\DimD_{\ell_\DimD,\ldots,\ell_\DimK}(\psi)|^2
\lesssim N^{\DimD-1}
\]
for $i=1,\dots,N$. For $i=1$ it is easy to verify the above estimates, since each of the sums contains at most two summands, with $(\ell_\DimD-(\DimD-1)/2,\ell_\DimK-(\DimK-1)/2) \in \{(0,0),(1,1)\}$, and the functions $X_{\ell_\DimD,\ldots,\ell_\DimK}$ are bounded. For $i=2,\dots,N$, these estimates follow from Propositions \ref{prp:indici_bassi_primo_piufattori} and \ref{prp:indici_bassi_primo_peso_secondo_stratopiu}, applied with $m=\ell_\DimK$ and $\ell=\ell_\DimD$.
\end{proof}

\subsection{Properties of the weight}\label{s:weight_lemma}

We shall need some properties of the weights $\weight_r: \sfera^\DimD \times \sfera^\DimD \to [0,\infty)$ defined in \eqref{eq:weight}.
The following lemma extends \cite[Lemma 5.1]{CaCiaMa}, where only the case $\DimD=2$, $\DimK=1$ was treated. We refer to \cite[Lemma 12]{MSi} and \cite[Lemma 4.1]{M} for analogous results.

\begin{lemma}\label{lem:weight}
For all $r>0$ and $\alpha,\beta \geq 0$ such that $\alpha+\beta > \DimD+\DimK$ and $\alpha < \min\{\DimD-\DimK,\DimK\} $, and for all $z' \in \sfera^\DimD$,
\begin{equation}\label{eq:weightsint}
\int_{\sfera^\DimD} (1+\dist(z,z')/r)^{-\beta} (1+\weight_r(z,z'))^{-\alpha} \,d\meas(z) \lesssim_{\alpha,\beta} V(z',r).
\end{equation}
Moreover
\begin{equation}\label{eq:weightsineq}
1+ \weight_r(z,z') \lesssim (1+\dist(z,z')/r)
\end{equation}
for all $r>0$ and $z,z' \in \sfera^\DimD$.
\end{lemma}
\begin{proof}

Due to the compactness of $\sfera^\DimD$, both \eqref{eq:weightsint} and \eqref{eq:weightsineq} are obvious
for $r \geq 1$. In the following we assume therefore that $r<1$.

To prove \eqref{eq:weightsineq}, we observe that, for all $\spnt{\omega}{\psi}, \spnt{\omega'}{\psi'} \in \sfera^\DimD$,
\begin{equation}\label{eq:ineq_weight_dist}
1+\frac{|\psi|}{\max\{r,|\psi'|\}}
\simeq 1 + \frac{|\psi-\psi'|}{\max\{r,|\psi'|\}}
\lesssim 1 + \dist(\spnt{\omega}{\psi},\spnt{\omega'}{\psi'})/r.
\end{equation}
The last inequality follows immediately from \eqref{eq:rho-dist} in the case $\max\{|\psi|,|\psi'|\} < \pi/4$, and it is trivial when $|\psi'| > \pi/8$ (since $|\psi|/\max\{r,|\psi'|\} \lesssim 1$ in that case); in the remaining case ($|\psi| \geq \pi/4$ and $|\psi'| \leq \pi/8$), the points $\spnt{\omega}{\psi}$ and $\spnt{\omega'}{\psi'}$ belong to disjoint compact subsets of $\sfera^d$, whence
\begin{equation}\label{eq:punti_distanti}
\dist(\spnt{\omega}{\psi},\spnt{\omega'}{\psi'}) \simeq 1 \simeq |\psi-\psi'|
\end{equation}
and the desired inequality follows.

In order to prove \eqref{eq:weightsint}, we fix $z' = \spnt{\omega'}{\psi'} \in \sfera^\DimD$ and split the integral in the left-hand side of \eqref{eq:weightsint} into the sum $\sum_{j=0}^3 \cI_j$, where
\[
\cI_j=
\int_{\cS_j} (1+\dist(z,z')/r)^{-\beta} (1+\weight_r(z,z'))^{-\alpha} \,d\meas(z)
\]
and
\begin{align*}
\cS_0&={\left\{\spnt{\omega}{\psi}\in \sfera^\DimD \tc {\max\{|\psi|,| \psi'|\}} \geq \pi/4 \right\}},\\
\cS_1&=\left\{\spnt{\omega}{\psi} \in \sfera^\DimD \setminus \cS_0 \tc  {\rho_{R,\sfera^\DimK}(\omega,\omega')}^{1/2}
 \leq \frac{\rho_{R,\sfera^\DimK}(\omega,\omega')}{\max\{|\psi|,| \psi'|\}}\right\},\\
\cS_2&=\left\{\spnt{\omega}{\psi} \in \sfera^\DimD \setminus (\cS_0 \cup \cS_1) \tc |\psi'| \leq |\psi|/2 \right\},\\
\cS_3&=\left\{\spnt{\omega}{\psi} \in \sfera^\DimD \setminus (\cS_0 \cup \cS_1) \tc |\psi|/2 < |\psi'| \right\}.
\end{align*}

We first estimate $\cI_0$.
In the case $|\psi'| > \pi/8$, we use \eqref{eq:rho-dist-far} to conclude that
\[
\cI_0
\lesssim \int_{\sfera^\DimD} (1+\dist_R(z,z')/r)^{-\beta} \,d\meas(z) 
\lesssim r^\DimD \simeq V(z',r),
\]
since $r< 1$ and $\beta>\DimD$ (cf.\ \cite[Lemma 4.4]{DOS}). In the case $|\psi'| \leq \pi/8$, instead, $\dist(z,z') \simeq |\psi| \simeq 1$ by \eqref{eq:punti_distanti} for all $z \in \cS_0$, and
\[
\cI_0 \simeq r^\beta \max\{r,|\psi'|\}^\alpha 
= r^\DimD \max\{r,|\psi'|\}^\DimK \frac{r^{\beta-\DimD}}{\max\{r,|\psi|\}^{\DimK-\alpha}} \lesssim V(z',r),
\]
by \eqref{eq:volume}, since $\beta - \DimD > \DimK - \alpha>0$.

In order to estimate $\cI_1$, we decompose $\beta = \beta_1 + \beta_2$, with $\beta_1 > \DimD-\DimK-\alpha$ and $\beta_2 > 2\DimK$.
Thus \eqref{eq:rho-dist} and \eqref{eq:ineq_weight_dist} imply
\[\begin{split}
\cI_1
&\simeq
\int_{\cS_1} (1+\dist(z,z')/r)^{-\beta}  \left(1+\frac{|\psi-\psi'|}{\max\{r,|\psi'|\}}\right)^{-\alpha} \,d\meas(z) 
\\
&\leq
(\max\{r,|\psi'|\}/r)^{\alpha} \int_{ \cS_1} (1+\dist(z,z')/r)^{-\beta}  \left(1+{|\psi-\psi'|/r} \right)^{-\alpha} \,d\meas(z) \\
&\lesssim
(\max\{r,|\psi'|\}/r)^{\alpha}
\int_{\cS_1} (1+{\dist_{R,\sfera^\DimK}(\omega,\omega')}^{1/2}/r)^{-\beta_2} (1+|\psi-\psi'|/r)^{-\alpha-\beta_1} \,d\meas(\spnt{\omega}{\psi}) \\
&\lesssim
(\max\{r,|\psi'|\}/r)^{\alpha} \int_{\sfera^\DimK} (1+{\dist_{R,\sfera^\DimK}(\omega,\omega')}/r^2)^{-\beta_2/2} \,d\omega \\
&\qquad\times
\int_{[-\pi/4,\pi/4]^{\DimD-\DimK}} \left(1+|\psi-\psi'|/r \right)^{-\alpha-\beta_1} \,d\psi_\DimD \ldots \,d\psi_{\DimK+1} \\
&\lesssim 
( \max\{r,|\psi'|\}/r )^{\alpha} r^{2\DimK} r^{\DimD-\DimK}
= r^\DimD \max\{r,|\psi'|\}^\DimK (r/\max\{r,|\psi'|\})^{\DimK-\alpha}
\lesssim V(z',r),
\end{split}\]
since $\beta_2/2 > \DimK$ and $\alpha < \DimK$.

In order to estimate $\cI_2$, instead, we write $\beta = \tilde\beta_1 + \tilde\beta_2$, with $\tilde\beta_1 > \DimD - \alpha$ and $\tilde\beta_2 > \DimK$, so, again by \eqref{eq:rho-dist},
\[\begin{split}
\cI_2
&\simeq
\int_{\cS_2} \left(1+\frac{|\psi-\psi'|}{r} + \frac{\dist_{R,\sfera^\DimK}(\omega,\omega')}{r\max\{|\psi|,|\psi'|\}}\right)^{-\beta}
\,\left(1+\frac{|\psi|}{\max\{r,|\psi'|\}}\right)^{-\alpha} \,d\meas(\spnt{\omega}{\psi}) \\
&\lesssim
\int_{2|\psi'| \leq |\psi| \leq \pi/4} \left(1+\frac{|\psi|}{r}\right)^{-\tilde\beta_1} \left(1+\frac{|\psi|}{\max\{r,|\psi'|\}}\right)^{-\alpha} \\
&\qquad\times
\int_{\sfera^\DimK} \left(1+\frac{\dist_{R,\sfera^\DimK}(\omega,\omega')}{r|\psi|}\right)^{-\tilde\beta_2} \,d\omega \,d\psi_\DimD \ldots \,d\psi_{\DimK+1}  \\
&\lesssim
(\max\{r,|\psi'|\}/r)^\alpha \int_{[-\pi/4,\pi/4]^{\DimD-\DimK}} \left(1+\frac{|\psi|}{r}\right)^{-\tilde\beta_1-\alpha} (r|\psi|)^\DimK \,d\psi_\DimD \ldots \,d\psi_{\DimK+1}  \\
&\lesssim 
(\max\{r,|\psi'|\}/r)^\alpha r^{\DimD+\DimK}
= r^\DimD \max\{r,|\psi'|\}^\DimK (r/\max\{r,|\psi'|\})^{\DimK-\alpha}
\lesssim V(z',r)
\end{split}\]
where we used the fact that $\max\{|\psi|,| \psi'|\} \simeq |\psi-\psi'| \simeq |\psi|$ on $\cS_2$.

Finally, to estimate $\cI_3$, we decompose $\beta = \tilde\beta_1 + \tilde\beta_2$ as above and get
\[\begin{split}
\cI_3
&\lesssim
\int_{\cS_3} \left(1+\frac{|\psi-\psi'|}{r}\right)^{-\tilde\beta_1} \left(1+\frac{\rho_{R,\sfera^\DimK}(\omega,\omega')}{r|\psi'|}\right)^{-\tilde\beta_2} \,d\meas(\spnt{\omega}{\psi}) \\
&\lesssim
(r |\psi'|)^\DimK \int_{[-\pi/4,\pi/4]^{\DimD-\DimK}} \left(1+\frac{|\psi-\psi'|}{r}\right)^{-\tilde\beta_1} \,d\psi
\lesssim r^\DimD |\psi'|^\DimK  
\lesssim V(z',r),
\end{split}\]
where we used the fact that $\max\{|\psi|,|\psi'|\} \simeq |\psi'|$ on $\cS_3$.
\end{proof}

\subsection{Proof of the main result}
The previous estimates finally allow us to verify the assumptions of the abstract theorem in Section \ref{s:abstracttheorem} and prove our multiplier theorem for the Grushin operators $\opL_{\DimD,\DimK}$.

\begin{proof}[Proof of Theorem \ref{thm:main}]
Let $\alpha \in [0,\min\{\DimD-\DimK,\DimK\})$. We apply Theorem \ref{thm:abstractmult} with $(X,\dist,\mu) = (\sfera^\DimD,\dist,\meas)$, $\AbsOp = \opL_{\DimD,\DimK}$, $q=2$, $\absDim = \DimD+\DimK - \alpha$, $\pi_r = (1+\weight_r)^\alpha$. Note that the assumptions \ref{en:abs_doubling} and \ref{en:abs_heatkernel} easily follow from \cite{DzS}; as a matter of fact, \ref{en:abs_doubling} also follows from Proposition \ref{prp:subriemannian}, and \ref{en:abs_heatkernel} could be derived from Proposition \ref{prp:weighted_plancherel} via the results of \cite{Melrose,S} (cf.\ the discussion in \cite{CaCiaMa}). Moreover, the assumptions \ref{en:abs_weightgrowth} and \ref{en:abs_weightint} are proved in Lemma \ref{lem:weight}, while the assumption \ref{en:abs_plancherel} is proved in Proposition \ref{prp:weighted_plancherel}. By choosing $\alpha$ sufficiently close to $\min\{\DimD-\DimK,\DimK\}$, we can make $\absDim = \DimD+\DimK-\alpha$ arbitrarily close to $D = \max\{\DimD,2\DimK\}$, and the desired results follow.
\end{proof}

\section{Appendix. Proof of the abstract multiplier theorem}\label{s:appendix}
\begin{proof}[Proof of Theorem \ref{thm:abstractmult}]
Similarly as in \cite{M}, for all $r \in (0,\infty)$, $\beta \in [0,\infty)$, $p \in [1,\infty]$ and $K : X \times X \to \CC$, 
we define the norm $\vvvert K \vvvert_{p,\beta,r}$ 
\[
\vvvert K \vvvert_{p,\beta,r}
= \esssup_{z' \in X} \mu(B(z',r))^{1/p'} \| (1+\dist(\cdot,z')/r)^\beta \, K(\cdot,z') \|_{L^p(X)}, 
\]
where $p'=p/(p-1)$ is the conjugate exponent to $p$; if $r \in (0,1]$, we also define the norm  $\vvvert K \vvvert_{p,\beta,r}^*$ by 
\[
\vvvert K \vvvert_{p,\beta,r}^*
= \esssup_{z' \in X} \mu(B(z',r))^{1/p'} \| (1+\dist(\cdot,z')/r)^\beta \, \absWeight_r(\cdot,z') \, K(\cdot,z') \|_{L^p(X)}.
\]

Due to the doubling condition and the heat kernel bounds, we can apply \cite[Theorem 6.1]{M} to obtain that, for all $\epsilon > 0$, all $\beta \geq 0$, all $R \in (0,\infty)$ and all $F : \RR \to \CC$ supported in $[-R^2,R^2]$,
\begin{align}
\vvvert \Kern_{F(\AbsOp)} \vvvert_{2,\beta,R^{-1}} &\lesssim_{\beta,\epsilon} \|F(R^2 \cdot)\|_{\Sob{\infty}{\beta+\epsilon}}, \label{eq:nswl2bd} \\
\| F(\AbsOp) \|_{L^1(X) \to L^1(X)} &\lesssim_{\epsilon} \|F(R^2 \cdot)\|_{\Sob{\infty}{Q/2+\epsilon}}  \label{eq:nsl1bd},
\end{align}
where $Q$ is the doubling dimension of $(X,\dist,\mu)$. It is worth noting that, since $\absWeight_r \gtrsim 1$ by \eqref{eq:weightsineq22}, the estimate \eqref{eq:weightsint11} trivially holds for all $\beta > Q$, $r>0$ and $y \in X$ \cite[Lemma 4.4]{DOS}; so it is not restrictive to assume in what follows that $\absDim \leq Q$.

Set $A_t = \exp(-t^2\AbsOp)$ if $t \in [0,\infty)$ and $A_t = 0$ if $t = \infty$. From \eqref{eq:nswl2bd} we deduce that, for all $t \in [0,\infty]$, all $\epsilon > 0$, all $\beta \geq 0$, all $R \in (0,\infty)$ and all $F : \RR \to \CC$ supported in $[R/16,R]$,
\[
\vvvert \Kern_{F(\sqrt{\AbsOp}) (1-A_t)} \vvvert_{2,\beta,R^{-1}} \lesssim_{\beta,\epsilon} \|F(R \cdot)\|_{\Sob{\infty}{\beta+\epsilon}} \min \{1,(Rt)^2 \}.
\]

Let $\xi \in C_c((-1/16,1/16))$ be nonnegative with
\[
\int_R \xi(t) \,dt = 1 \qquad\text{and}\qquad \int_\RR t^k \xi(t) \,dt = 0 \text{ for }  k=1,\dots,2Q+2.
\]
(cf.\ \cite[eq.\ (18)]{M}). Then by Young's inequality we obtain that, for all $t \in [0,\infty]$, all $\epsilon > 0$, all $\beta \geq 0$, all $R \in [1,\infty)$ and all $F : \RR \to \CC$ supported in $[R/8,7R/8]$,
\[
\vvvert \Kern_{(\xi* F)(\sqrt{\AbsOp})(1-A_t)} \vvvert_{2,\beta,R^{-1}} \lesssim_{\beta,\epsilon} \|F(R \cdot)\|_{\Sob{\infty}{\beta+\epsilon}} \min \{1,(Rt)^2 \}.
\]
In particular, by \eqref{eq:weightsineq22} and Sobolev's embedding, for all $t \in [0,\infty]$, all $\epsilon > 0$, all $\beta \geq 0$, all $N \in \NN \setminus \{0\}$ and all $F : \RR \to \CC$ supported in $[N/8,7N/8]$,
\begin{equation}\label{eq:we1}
\vvvert \Kern_{(\xi*F)(\sqrt{\AbsOp})(1-A_t)} \vvvert_{2,\beta,N^{-1}}^* \lesssim_{\beta,\epsilon} \|F(N \cdot)\|_{\Sob{q}{\beta+M_0+1/q+\epsilon}} \min \{1,(Nt)^2 \}.
\end{equation}

On the other hand, by \eqref{eq:abstract_plancherel}, for all $t \in [0,\infty]$, all $N \in \NN \setminus \{0\}$ and all $F : \RR \to \CC$ supported in $[N/16,N]$,
\[
\vvvert \Kern_{F(\sqrt{\AbsOp}) (1-A_t)} \vvvert_{2,0,N^{-1}}^* \lesssim \|F(N \cdot)\|_{N,q} \min \{1,(Nt)^2 \}.
\]
Hence, by \cite[eq.\ (4.9)]{DOS}, for all $t \in [0,\infty]$, all $N \in \NN \setminus \{0\}$ and all $F : \RR \to \CC$ supported in $[N/8,7N/8]$,
\begin{equation}\label{eq:we2}
\vvvert \Kern_{(\xi*F)(\sqrt{\AbsOp})(1-A_t)} \vvvert_{2,0,N^{-1}}^* \lesssim \|F(N \cdot)\|_{L^q} \min \{1,(Nt)^2 \}.
\end{equation}

Interpolation of \eqref{eq:we1} and \eqref{eq:we2} gives that, for all $t \in [0,\infty]$, all $\epsilon > 0$, all $\beta \geq 0$, all $N \in \NN \setminus \{0\}$ and all $F : \RR \to \CC$ supported in $[N/4,3N/4]$,
\[
\vvvert \Kern_{(\xi*F)(\sqrt{\AbsOp})(1-A_t)} \vvvert_{2,\beta,N^{-1}}^* \lesssim_{\beta,\epsilon} \|F(N \cdot)\|_{\Sob{q}{\beta+\epsilon}} \min \{1,(Nt)^2 \}.
\]
By \eqref{eq:weightsint11} and H\"older's inequality we then deduce that, for all $r \in [0,\infty)$, all $t \in [0,\infty]$, all $s > \absDim/2$, all $\epsilon \in [0,s-\absDim/2)$,  all $N \in \NN \setminus \{0\}$ and all $F : \RR \to \CC$ supported in $[N/4,3N/4]$,
\begin{equation}\label{eq:mollified_l1}
\begin{split}
\esssup_{z' \in X} &\int_{X \setminus B(z',r)} |\Kern_{(\xi*F)(\sqrt{\AbsOp})(1-A_t)}(z,z')| \,d\mu(z) \\
 &\leq (1+Nr)^{-\epsilon} \vvvert \Kern_{(\xi*F)(\sqrt{\AbsOp})(1-A_t)} \vvvert_{1,\epsilon,N^{-1}} \\
 &\lesssim_{s,\epsilon} (1+Nr)^{-\epsilon} \vvvert \Kern_{(\xi*F)(\sqrt{\AbsOp})(1-A_t)} \vvvert_{2,\beta,N^{-1}}^* \\
 &\lesssim_{s,\epsilon} (1+Nr)^{-\epsilon} \|F(N \cdot)\|_{\Sob{q}{s}} \min \{1,(Nt)^2 \},
\end{split}
\end{equation}
where $\beta \in (\absDim/2+\epsilon,s)$.

On the other hand, if $D$ is the $\dist$-diameter of $X$, by \eqref{eq:weightsint11}, H\"older's inequality, \eqref{eq:abstract_plancherel} and \cite[Proposition 4.6]{DOS},  for all $s > \absDim/2$, all $\epsilon \in [0,\min\{s-\absDim/2,\absDim/2\})$, all  $N \in \NN \setminus \{0\}$ and all $F : \RR \to \CC$ supported in $[N/4,3N/4]$,
\begin{equation}\label{eq:nomollified_l1}
\begin{split}
\| (F-\xi*F)(\sqrt{\AbsOp}) \|_{1 \to 1} &= \vvvert \Kern_{(F-\xi*F)(\sqrt{\AbsOp})} \vvvert_{1,0,N^{-1}} \\
&\lesssim_{s,\epsilon} \vvvert \Kern_{(F-\xi*F)(\sqrt{\AbsOp})} \vvvert_{2,\beta,N^{-1}}^* \\
&\leq (1+ND)^\beta \vvvert \Kern_{(F-\xi*F)(\sqrt{\AbsOp})} \vvvert_{2,0,N^{-1}}^* \\
&\lesssim_{s,\epsilon} N^\beta \| (F-\xi*F)(N \cdot) \|_{N,q} \\
&\lesssim_{s,\epsilon} N^{-\epsilon} \|F(N \cdot)\|_{\Sob{q}{\epsilon+\beta}} \\
&\lesssim_{s,\epsilon} N^{-\epsilon} \|F(N \cdot)\|_{\Sob{q}{s}},
\end{split}
\end{equation}
where $\beta \in (\absDim/2,\min\{\absDim,s-\epsilon\})$.

Finally, observe that, if $\supp F \subseteq [0,1]$, then, by \eqref{eq:weightsint11}, H\"older's inequality and \eqref{eq:abstract_plancherel} applied with $r=N=1$,
\begin{equation}\label{eq:trivial_l1}
\begin{split}
\| F(\sqrt{\AbsOp}) \|_{1 \to 1} &= \vvvert \Kern_{F(\sqrt{\AbsOp})} \vvvert_{1,0,1} \\
&\lesssim \vvvert \Kern_{F(\sqrt{\AbsOp})} \vvvert_{2,\absDim,1}^* \\
&\leq (1+D)^\absDim \vvvert \Kern_{F(\sqrt{\AbsOp})} \vvvert_{2,0,1}^* \\
&\lesssim \| F \|_{N,q} \\
&\leq \| F \|_\infty.
\end{split}
\end{equation}

Combining \eqref{eq:mollified_l1} (applied with $t=\infty$, and $\epsilon = r=0$) and \eqref{eq:nomollified_l1} (applied with $\epsilon=0$) gives in particular that, for all $s > \absDim/2$, all  $N \in \NN \setminus \{0\}$ and all $F : \RR \to \CC$ supported in $[N/4,3N/4]$,
\begin{equation}\label{eq:l1_compact}
\| F(\sqrt{\AbsOp}) \|_{1 \to 1} \lesssim_s \| F(N \cdot) \|_{\Sob{q}{s}}.
\end{equation}
This estimate, combined with \eqref{eq:trivial_l1}, easily gives a weak version of part \ref{en:abs_compact}: namely, for all $s>\absDim/2$ and $F : \RR \to \CC$ supported in $[1/2,1]$,
\begin{equation}\label{eq:abs_compact_weak}
\sup_{t>0} \| F(t \sqrt{\AbsOp}) \|_{1 \to 1} \lesssim_s \, \|F\|_{\Sob{q}{s}}.
\end{equation}

We now prove the full version of part \ref{en:abs_compact}. Fix an even cutoff function $\chi \in C^\infty_c(\RR)$ with $\chi(0) = 1$ and $\supp \chi \subseteq [-1,1]$. Let $F : \RR \to \CC$ be supported in $[-1,1]$ and set $\tilde F = F - F(0) \chi$. Note that, for all $k \in \NN$,
\[
\|F(0) \chi(\sqrt{\cdot}) \|_{C^k} \lesssim_k \|F(0) \chi \|_{C^{2k}} \lesssim_k |F(0)| \lesssim_{s} \|F\|_{L^q_s},
\]
by Sobolev's embedding, provided $s>1/q$. In particular, from \eqref{eq:nsl1bd} it follows that
\begin{equation}\label{eq:value_in_zero}
\sup_{t>0} \|F(0) \chi(t \sqrt{\AbsOp}) \|_{1 \to 1} \lesssim_{s} \|F\|_{L^q_s}
\end{equation}
for all $s>1/q$, and moreover
\[
\|\tilde F\|_{\Sob{q}{s}} \lesssim_s \|F\|_{\Sob{q}{s}}.
\]
Let now $\xi \in C^\infty_c(\RR)$ be such that $\supp \xi \subseteq (1/2,2)$ and $\sum_{k \in \ZZ} \xi(2^k \cdot) = 1$ on $(0,\infty)$. Decompose $\tilde F = \sum_{k \in \NN} \tilde F_k(2^k \cdot)$ on $[0,\infty)$, where $\tilde F_k = \tilde F(2^{-k} \cdot) \, \xi$; since $\supp \tilde F_k \subseteq (1/2,2)$, from \eqref{eq:abs_compact_weak} we deduce that
\[
\sup_{t>0} \|\tilde F_k(t\sqrt{\AbsOp}) \|_{1 \to 1} \lesssim_{\beta} \|\tilde F_k\|_{\Sob{q}{\beta}}
\]
provided $\beta>\absDim/2$. On the other hand, arguing as in the proof of \cite[Lemma 4.8]{MAif}, one deduces that, for all $\beta\geq 0$ and $s > \max\{\beta,1/q\}$, there exists $\epsilon>0$ such that
\[
\|\tilde F_k\|_{\Sob{q}{\beta}} \lesssim_{\beta,s} \|\tilde F_k\|_\infty + 2^{-k\epsilon} \|\tilde F\|_{\Sob{q}{s}} \lesssim_s 2^{-k\epsilon} \|\tilde F\|_{\Sob{q}{s}};
\]
the latter estimate is due to the fact that $\tilde F(0) = 0$ and, by Sobolev's embedding, if $\|\tilde F\|_{\Sob{q}{s}} < \infty$ for some $s>1/q$, then $\tilde F$ is H\"older-continuous. In conclusion, for all $t>0$ and $s>\absDim/2$,
\begin{equation}\label{eq:subtracted_value_in_zero}
\|\tilde F(t\sqrt{\AbsOp}) \|_{1 \to 1} \leq \sum_{k\in \NN} \| \tilde F_k(2^k t \sqrt{\AbsOp}) \|_{1 \to 1} \lesssim_s \sum_{k \in\NN} 2^{-k\epsilon} \|\tilde F\|_{\Sob{q}{s}} \lesssim_s \|F\|_{\Sob{q}{s}};
\end{equation}
combining the estimates \eqref{eq:value_in_zero} and \eqref{eq:subtracted_value_in_zero} gives part \ref{en:abs_compact}.

As for part \ref{en:abs_mh},
since the right-hand side of \eqref{eq:main_mh_wt11} is essentially independent of the cut-off function $\eta$, we may assume that $\supp \eta \subseteq (1/4,1)$ and $\sum_{k \in \ZZ} \eta(2^k \cdot) = 1$ on $(0,\infty)$. 
Then, arguing as in \cite[proof of Theorems 3.1 and 3.2]{DOS}, by the use of the dyadic decomposition $F = \sum_{k \in \NN} \eta(2^{-k} \cdot) F$ and an application of \cite[Theorem 1]{DMc}, from \eqref{eq:mollified_l1} and \eqref{eq:nomollified_l1} we obtain that, for all $F : \RR \to \CC$ supported in $[1/2,\infty)$,
\begin{equation}\label{eq:l1weak}
\| F(\sqrt{\AbsOp}) \|_{L^1 \to L^{1,\infty}} \lesssim_s \sup_{k \in \NN} \| \eta \, F(2^k \cdot) \|_{\Sob{q}{s}}.
\end{equation}
Via a partition of unity subordinated to $\{(1/2,\infty),(-\infty,1)\}$, we can now combine \eqref{eq:l1weak} and \eqref{eq:trivial_l1} and obtain part \ref{en:abs_mh}.
\end{proof}

%%%%%%%%%%%%%%%%%%%%%%%%%%%%%%%%%%%%%%%%%%%%%%%%%%%%%%%%%%%%%%%%%%

\end{document}